\documentclass[mnsc,nonblindrev]{informs3aa} 

\OneAndAHalfSpacedXI

\usepackage{natbib}
 \bibpunct[, ]{(}{)}{,}{a}{}{,}%
 \def\bibfont{\small}%
\TheoremsNumberedThrough     

\EquationsNumberedThrough    

\MANUSCRIPTNO{\ } 
\usepackage{amsfonts, amssymb, amsmath, mathtools}
\usepackage[small]{titlesec}
\usepackage{enumitem}  
\usepackage{xcolor}
\usepackage{graphicx}
\usepackage[caption=false]{subfig}
\usepackage[hidelinks]{hyperref}
\usepackage{bbm}
\usepackage{float}
\usepackage{algpseudocode}
\usepackage{algorithm}
\usepackage{endnotes}
\let\footnote=\endnote
\usepackage{tikz}
\usetikzlibrary{decorations.pathreplacing}

\newcommand{\Var}{\operatorname{Var}}
\newcommand{\E}{\operatorname{\mathrm{E}}}
\newcommand{\Prob}{\operatorname{\mathbb{P}}}
\renewcommand{\complement}{c}

\begin{document}


 \RUNAUTHOR{Bumpensanti and Wang}

 \RUNTITLE{A Re-solving Heuristic for Network Revenue Management}

 \TITLE{A Re-solving Heuristic  with Uniformly Bounded Loss for Network Revenue Management}

\ARTICLEAUTHORS{%
\AUTHOR{Pornpawee Bumpensanti, He Wang}
\AFF{School of Industrial and Systems Engineering, Georgia Institute of Technology, Atlanta, GA 30332,\\ \EMAIL{pornpawee@gatech.edu}, \EMAIL{he.wang@isye.gatech.edu}}
} 

\ABSTRACT{%
We consider the canonical \emph{(quantity-based) network revenue management} problem, where a firm accepts or rejects incoming customer requests irrevocably in order to maximize expected revenue given limited resources. 
Due to the curse of dimensionality, the exact solution to this problem by dynamic programming is  intractable when the number of resources is large.
We study a family of re-solving heuristics that periodically re-optimize an approximation to the original problem known as the deterministic linear program (DLP), where random customer arrivals are replaced by their expectations. We find that, in general, frequently re-solving the DLP produces the same order of revenue loss as one would get without re-solving, which scales as the square root of the time horizon length and resource capacities.
By re-solving the DLP at a few selected points in time and applying thresholds to the customer acceptance probabilities, we design a new re-solving heuristic whose revenue loss is uniformly bounded by a constant that is independent of the time horizon and resource capacities. 
}%



\maketitle

%

\section{Introduction}
The \emph{network revenue management} (NRM) problem \citep{williamson1992airline,gallego1997multiproduct} is a classical model that has been extensively studied in the revenue management literature for over two decades. {\color{black} The problem is concerned with maximizing revenue given limited resource and time, and has a wide range of applications in the airline, retail, advertising, and hospitality industries \citep[see examples in][]{talluri_theory_2004}. However, the exact solution to the NRM problem is difficult to compute when the number of resources is large.  Heuristics proposed in the previous literature typically have optimality gaps, i.e., expected revenue losses compared to the optimal solution, that increase with the time horizon and the resource capacities.
In this paper, we propose a new heuristic for the NRM problem for which the revenue loss is independent of the time horizon and the resource capacities.
}

The NRM problem is stated as follows: there is a set of resources with finite capacities that are available for a finite time horizon.
Heterogeneous customers arrive sequentially over time. Customers are divided into different classes based on their consumption of resources and the prices they pay. Each class of customer may request multiple types of resources and multiple units of each resource. 
Upon a customer's arrival, a decision maker must irrevocably accept or reject the customer. If the customer is accepted and there is enough remaining capacities, she consumes the resources requested and pays a fixed price associated with her class.
Otherwise,  if the customer is rejected, no revenue is collected and no resources are used.
Unused resources at the end of the finite horizon are perishable and have no salvage value.
The decision maker's objective is to maximize the expected revenue earned during the finite horizon. 

We note that the formulation stated above is more specifically known as the ``quantity-based'' NRM problem. In another formulation referred to as the ``price-based'' NRM problem, the decision maker chooses posted prices rather than accept/reject decisions. The two formulations are different, but are equivalent in some special cases
\citep{maglaras2006dynamic}. We focus on the quantity-based formulation in this paper.

\begin{figure}[!htb]
	\centering
			\includegraphics[width=0.3\columnwidth]{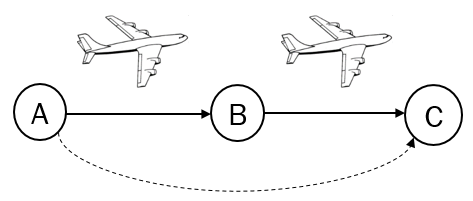}%
		\caption{A flight network of two flight legs ($A\to B, B\to C$) and three itineraries.}\label{fig:flight}
\end{figure}

A classical application of the NRM problem is in airline seat revenue management \citep{williamson1992airline, gallego1997multiproduct}. Here, the resources correspond to flight legs and the capacity corresponds to the number of seats on each flight. The resources are perishable on the date of flight departure.
Arriving customers are divided into separate classes defined by combinations of  itinerary and fare. A simple flight network of two flight legs and three itineraries is shown in Fig.~\ref{fig:flight}. 
The objective of the airline is to maximize the expected revenue earned from allocating available seats to different classes of customers. Notice that the problem cannot be decomposed for each individual flight leg, since some itineraries use multiple resources simultaneously (e.g., in Fig.~\ref{fig:flight}, customers traveling from $A$ to $C$ would request itinerary $A\to B\to C$). 
In practice, the huge size of airline networks makes solving this problem challenging.

\subsection{Deterministic LP approximation and re-solving heuristics}

In theory, the NRM problem can be solved by dynamic programming; however, since the state space grows exponentially with the number of resources, the dynamic programming formulation is often intractable. Therefore, we focus on heuristics with provable performance guarantees in this paper. We define \emph{revenue loss} as the gap between the expected revenue of a heuristic policy and that of the optimal policy.
As common in the revenue management literature, the effectiveness of heuristic polices are evaluated in an asymptotic regime where resource capacities and  customer arrivals are both scaled proportionally by a factor of $k\ (k=1,2,\ldots)$. 
Intuitively, this asymptotic regime increases market size while keeping
resource scarcity, i.e., the ratio of capacity to demand, at a constant level. 
We assume this standard asymptotic setting throughout the paper. 
 
One popular heuristic for the NRM problem that is extensively studied in the academic literature and widely used in practice is based on the deterministic linear programming (DLP) approximation, where the  customer demand distributions are replaced by their expectations. 
The solution of the DLP can then be used to construct heuristic  policies. 
Under the asymptotic scaling defined above, \citet{GallegoOptimalDynamicPricing1994,gallego1997multiproduct} have shown that the revenue loss of DLP-based static control policies is $\Theta(\sqrt{k})$ when the system size is scaled by $k$. 
The book by \citet{talluri_theory_2004} provides a comprehensive overview of different types of DLP-based control policies such as booking limit control, bid-price control, etc., and their variations. 

An apparent weakness of the DLP approximation is that it 
ignores randomness in the arrival process and fails to incorporate information acquired through time. 
To include updated information, a simple approach is to re-optimize the DLP from time to time, while replacing the initial capacity in the DLP with the remaining capacity at each re-solving  point. 
The new solution to the updated DLP is then used to adjust control policies. The re-solving approach is intuitive and widely used in practice.  We refer to this family of solution techniques as \emph{re-solving heuristics}.
One might expect that re-solving the DLP would yield better performance since it includes updated information. Surprisingly, 
%
\citet{CooperAsymptoticBehaviorAllocation2002} provides a counter-example where the performance of booking limit control deteriorates by re-solving the DLP.
Furthermore, \citet{ChenResolvingstochasticprogramming2010} give an example where re-solving the DLP worsens the performance for bid-price control. \citet{JasinAnalysisDeterministicLPBased2013} analyze the performance of re-solving  of both booking limit and bid-price controls.
They showed that when the initial capacity and customer arrival rates are both scaled by $k$, the revenue loss of re-solving heuristics is $\Omega(\sqrt{k})$, even by optimizing over the re-solving schedule or increasing re-solving frequency.

Despite those negative results, we note that there are several ways to construct control policies from the DLP, so it is possible that some control policies are suitable for applying the re-solving technique, while others are not. 
Some recent literature draws attention to a specific type of control policy called \emph{probabilistic allocation}, which seems suitable for applying the re-solving technique. Probabilistic allocation control is a randomized algorithm that accepts each arriving customer with some probability.
Using the probabilistic allocation control,
\citet{ReimanAsymptoticallyOptimalPolicy2008} propose a heuristic policy that re-solves the DLP exactly once during the horizon. In their proposed policy, the re-solving time is random and determined endogenously by the heuristic policy.
In the asymptotic setting, \citet{ReimanAsymptoticallyOptimalPolicy2008} show that the revenue loss of their policy is $o(\sqrt{k})$. This is {an} improvement over the $\Theta(\sqrt{k})$ revenue loss of DLP-based static policies.

\citet{JasinReSolvingHeuristicBounded2012} consider an algorithm that is based on probabilistic allocation control and re-solves the DLP after each unit of time. They show the algorithm has a revenue loss of $O(1)$ when the system size is scaled by $k\to\infty$. A similar $O(1)$ revenue loss is obtained by \citet{wu2015algorithms} for the case of one resource. However, both \citet{JasinReSolvingHeuristicBounded2012} and \citet{wu2015algorithms}'s results require the optimal solution to DLP (before any updating) to be nondegenerate; this assumption will be formally stated in Section~\ref{sec:FR}, which seems to be central to the hardness of the NRM problem. Moreover, \citet{wu2015algorithms} show that when the optimal solution is nondegenerate but nearly degenerate, the constant factor in $O(1)$ can become arbitrarily large.  In this paper, we aim to establish a uniform $O(1)$ loss for the general NRM problem without assuming nondegeneracy.


\subsection{Main contributions}

We propose a new re-solving heuristic that has a uniformly bounded revenue loss when the system size is scaled by $k\to\infty$. {\color{black} (Recall that the rate of revenue loss is defined for a sequence of problems indexed by $k=1,2,\ldots$, where the capacities and arrival rates are multiplied by $k$, while other parameters are treated as constants.)}
The bound is uniform in the sense that it does not depend on ratio between capacities and time. Therefore, this result does not require the nondegeneracy assumption.
Our $O(1)$ bound improves the $o(\sqrt{k})$ bound in \citet{ReimanAsymptoticallyOptimalPolicy2008}, and also improves the $O(1)$  bound in \citet{JasinReSolvingHeuristicBounded2012}, where the constant factor requires nondegeneracy assumption and depends implicitly on problem instances.
We call our new algorithm  \textsf{Infrequent Re-solving with Thresholding (IRT)}. 
The intuition behind the \textsf{IRT} algorithm is that it is not necessary to update the DLP at early stage of the horizon, as the solution to the DLP barely changes after updating. It is sufficient to re-solve the DLP at a few carefully selected time points near the \emph{end} of the horizon. In total, the \textsf{IRT} algorithm has $O(\log\log k)$ re-solving times for a system with scaling size $k$. Furthermore, a ``thresholding'' technique is applied in case that the DLP solution after re-solving is nearly degenerate.
The re-solving schedule and the thresholds of the \textsf{IRT} algorithm are designed in such a way that the accumulated random deviations before the re-solving point can be corrected after re-solving with high probability.

Then, we give a tight performance bound of the re-solving heuristic proposed by \cite{JasinReSolvingHeuristicBounded2012}, but without assuming the optimal solution to the DLP is nondegenerate. 
The heuristic in \cite{JasinReSolvingHeuristicBounded2012}, which we call \textsf{Frequent Re-solving (FR)}, re-solves the DLP after each unit of time.
One would expect that by re-solving the DLP frequently and thus constantly updating capacity information, the decision maker can improve the expected revenue. Indeed, \citet{JasinReSolvingHeuristicBounded2012} have shown that under the nondegeneracy assumption,  the revenue loss of this policy is $O(1)$ when the system size is scaled  by $k\to\infty$.
However, we find that the revenue loss of this policy is $\Theta(\sqrt{k})$ in general, which has the same order of revenue loss as DLP-based static heuristics without any re-solving \citep{gallego1997multiproduct, talluri_analysis_1998, CooperAsymptoticBehaviorAllocation2002}.
In particular, Proposition \ref{thm: ho res lb} shows that there exists a problem instance where the revenue loss of this policy is at least $\Omega(\sqrt{k})$. To analyze this instance, we  used the Berry-Esseen bound and Freedman's inequality to show that the probability of revenue loss being larger than $\Omega(\sqrt{k})$ is bounded away from 0.
This result suggests that the nondegeneracy assumption made by \citet{JasinReSolvingHeuristicBounded2012} is necessary to obtain $O(1)$ revenue loss, and explains why the $O(1)$ factor in \cite{wu2015algorithms} must be arbitrarily large when the DLP optimal solution is converging to a degenerate point. Then, Proposition~\ref{thm: DLP and resolving ub} shows that the revenue loss of this policy is bounded above by $O(\sqrt{k})$ in the general case, which also improves the $o(k)$ bound in \cite{maglaras2006dynamic}. The proof is based on a key inequality that bounds the average remaining capacity as a function of the remaining time.

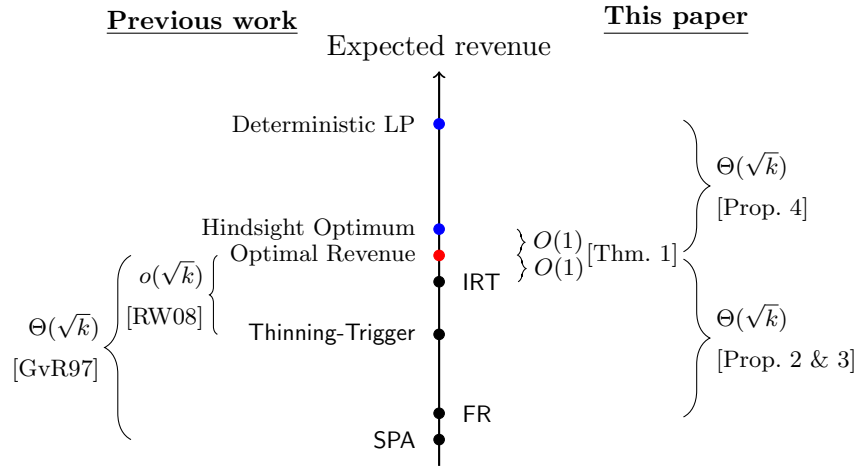
\begin{figure}[!htbp]
	\centering
	\begin{tikzpicture}[scale=0.7]
	\draw [thick,->] (0,2.5) --(0,10)node[anchor=south,black] {Expected revenue} ;
	\draw [decorate,decoration={brace,mirror,amplitude=10pt},xshift=-10pt,yshift=0pt]
	(-5.5,6.5) -- (-5.5,3) node [black,midway,xshift=-1cm,align=right,font=\footnotesize] 
	{$\Theta(\sqrt{k})$\\ {[GvR97]}};
	\draw [decorate,decoration={brace,mirror,amplitude=3pt},xshift=-6pt,yshift=0pt]
	(-4,6.5) -- (-4,5) node [black,midway,xshift=-0.7cm,align=right,font=\footnotesize] 
	{\footnotesize $o(\sqrt{k})$\\ {[RW08]}};
	\draw [decorate,decoration={brace,amplitude=3pt},xshift=-0.6pt,yshift=0pt]
	(1.5,7) -- (1.5,6.5) node [black,midway,xshift=0.55cm,align=left] 
	{\footnotesize $O(1)$};
	\draw [decorate,decoration={brace,amplitude=3pt},yshift=0pt]
	(1.5,6.5) -- (1.5,6) node [black,midway,xshift=0.55cm,align=left,font=\footnotesize] {
		$O(1)$};
	\draw [decorate,decoration={brace,amplitude=10pt},xshift=4pt,yshift=2pt]
	(4.5,9) -- (4.5,6.5) node [right,black,midway,xshift=0.3cm,align=left,font=\footnotesize] 
	{$\Theta(\sqrt{k})$\\ {[Prop.~\ref{prop: DLP and optimal lb}]}};
	\draw [decorate,decoration={brace,amplitude=10pt},xshift=4pt,yshift=-2pt]
	(4.5,6.5) -- (4.5,3.5) node [right,black,midway,xshift=0.3cm,align=left,font=\footnotesize] {
		$\Theta(\sqrt{k})$\\ {[Prop.~\ref{thm: ho res lb} \& \ref{thm: DLP and resolving ub}]}};
	\draw [fill, blue] (0,9) circle [radius=0.1]{};	
	\draw [fill, blue] (0,7) circle [radius=0.1]{};
	\draw [fill, red] (0,6.5) circle [radius=0.1]{};
	\draw [fill, black] (0,6) circle [radius=0.1]{};
	\draw [fill, black] (0,5) circle [radius=0.1]{};
	\draw [fill, black] (0,3.5) circle [radius=0.1]{};
	\draw [fill, black] (0,3) circle [radius=0.1]{};
	\node [left] (1) at (-0.25,9){\footnotesize Deterministic LP};
	\node [left] (2) at (-0.25,7){\footnotesize Hindsight Optimum};
	\node [left] (3) at (-0.25,6.5){\footnotesize Optimal Revenue};
	\node [right] (4) at (0.25,6){\footnotesize \textsf{IRT}};
	\node [left] (5) at (-0.25,5){\footnotesize \textsf{Thinning-Trigger}};
	\node [right] (6) at (0.25,3.5){\footnotesize \textsf{FR}};
	\node [left] (7) at (-0.25,3){\footnotesize \textsf{SPA}};
	\node [right] (4) at (2.6,6.5){\footnotesize [Thm. \ref{thm: TR}]};
	
	\node [above] at (-4.5,10.5) {\small \textbf{\underline{Previous work}}};
	\node [above] at (4.5,10.5) {\small \textbf{\underline{This paper}}};
	\end{tikzpicture}
	\caption{Summary of the results in the previous literature (on the left side) and our main results (on the right side).  The red node  ({\color{red}$\bullet$}) represents the expected revenue of the optimal policy (hard to compute); the blue nodes ({\color{black}$\bullet$}) represent upper bounds to the optimal revenue; and the black nodes ({\color{black}$\bullet$}) refer to revenues earned under different heuristics. The factor $k$ is the scale of both time horizon and capacities.}
	\label{fig: summary}
\end{figure}

In Fig.~\ref{fig: summary}, we summarize the performance of existing re-solving heuristics for the NRM problem. 
In this figure, the vertical axis represents the expected revenue which increases from the bottom to the top. 
We highlight the gap between different heuristics and upper bounds compared to the optimal revenue, which in principle can be obtained from dynamic programming but is hard to compute directly. The main result of the paper (Theorem~\ref{thm: TR}) simultaneously establishes an $O(1)$ upper bound of the hindsight optimum and an $O(1)$ revenue loss of the \textsf{IRT} algorithm.

\subsection{Other related work}


The re-solving heuristics defined in the NRM context is generally known as  \emph{certainty equivalent control} in dynamic programming. In certainty equivalent control, each random disturbance is fixed at a nominal value (e.g., its mean), and then an optimal control sequence {\color{black}for the certainty equivalence approximation} is found. Only the first control in the sequence is applied, the rest of them are discarded, and the same process is repeated in the next stage. 
An introduction to certainty equivalent control can be found in \citet[Section 6.1]{bertsekas2005dynamic}. \citet{secomandi2008analysis} discussed whether  certainty equivalent control guarantees performance improvement in the network revenue management setting.

The quantity-based NRM model can be generalized in several ways.
One extension assumes that the decision maker offers a set of products to each arriving customer, and customers choose some products from the offered set based on some discrete choice model \citep{talluri2004revenue, liu_choice-based_2008}.
Another stream of literature assumes that {\color{black}either the customers' arrival process or the distribution of their reservation price} is unknown, and requires the decision maker to learn the distribution exclusively from past observations \citep{besbes2012blind,jasin2015performance,ferreira2017online}. 
\cite{talluri_theory_2004,maglaras2006dynamic} discussed the case where the decision maker posts price (price-based NRM) versus the case where the decision maker chooses accept/reject (quantity-based NRM).

The NRM problem considered here is related to the online knapsack/secretary problem studied by
\citet{kleywegt1998dynamic, kleinberg2005multiple, babaioff2007knapsack}, \citet{arlotto2017uniformly}, and \citet{arlotto2018logarithmic}. 
{\color{black} In particular, \citet{arlotto2017uniformly} considers a multi-selection secretary problem, where the decision maker  sequentially selects i.i.d.\ random variables in order to maximize the expected value of the sum given a fixed budget. As such, by viewing each random variable as a customer arrival,
the multi-selection secretary problem is a special case of the NRM problem in which there is only a single resource and each customer requests exactly one unit of the resource.  \citet{arlotto2017uniformly} proposes an online policy that has a uniformly bounded regret compared to the optimal offline policy.
Their policy accepts or rejects an arriving customer by comparing the budget ratio, i.e., ratio of remaining budget to remaining arrivals, to some fixed thresholds. 
However, it is unclear whether their technique can be generalized to the general NRM setting with multiple resources, since the thresholds in their policy are specifically defined for a single resource.

Recently, \citet{vera2018bayesian} studies an online packing problem, which has the same mathematical formulation as the network revenue management problem. They propose a re-solving heuristic that achieves $O(1)$ revenue loss
without the nondegeneracy assumption and under mild assumptions on the customer arrival processes.
Unlike the \textsf{IRT} algorithm, their proposed algorithm re-solves the DLP every time there is an arrival; the algorithm then accepts that arrival if the acceptance probability from the DLP is greater than 0.5 and rejects it otherwise. 
Their proof is based on a novel argument that compensates the optimal offline algorithm and forces it to follow the decisions of their online algorithm. The design of their algorithm and their proof idea are significantly different from those in this paper.

}

%
%
%
%
%
%

{\color{black}
\subsection{Notation}\label{sec:notations}
For a positive integer $n$, let $[n]$ denote the set $\{1,\ldots,n\}$.
Given two real numbers $a \in \mathbb{R}$ and $b \in \mathbb{R}$, let $a \land b := \min\{a,b\}$, $a \lor b := \max\{a,b\}$, and $a^+ := a \lor 0$. For any real number $x$, let $\lfloor x \rfloor$ be the largest integer less than or equal to $x$, and let $\lceil x \rceil$ be the smallest integer greater than or equal to $x$.
For a set $S$, let $|S|$ denote the cardinality of $S$.
For two functions $f(T)$ and $g(T)>0$, we write $f(T)=O(g(T))$ if there exists a constant $M_1$ and a constant $T_1$ such that $f(T) \leq M_1 g(T)$ for all $T \geq T_1$; we write $f(T)=\Omega(g(T))$ if there exists a constant $M_2$ and a constant $T_2$ such that $f(T) \geq M_2 g(T)$ for all $T \geq T_2$. If $f(T)=O(g(T))$ and $f(T)=\Omega(g(T))$ both hold, we denote it by $f(T)=\Theta(g(T))$.
}

\section{Problem Formulation and Approximations}
\label{sec: Model}


Suppose there is a finite horizon with length $T$. There are $n$ classes of customers indexed by $j\in[n]$. 
The arrival process of customers in class $j$, $\{\Lambda_j(t),0\leq t\leq T\}$, follows a Poisson process of rate $\lambda_j$. We let $\Lambda_j(t_1,t_2)$ denote the number the arrivals of class $j$ customers during $(t_1,t_2]$ for $0\leq t_1 <t_2\leq T$, i.e., $\Lambda_j(t_1,t_2)=\Lambda_j(t_2)-\Lambda_j(t_1).$ Arrival processes of different classes are independent. 
Upon arrival, each customer must either be accepted or rejected.
{\color{black} Let $r_j$ denote the revenue received by accepting a class $j$ customer  and $r=[r_1,\ldots,r_n]^\top$ be the vector of such revenues. 
There are $m$ resources indexed by $l\in[m]$, where resource $l$ has initial capacity $C_l$. The vector of the initial capacities is given by $C=[C_1, \ldots,C_m]^\top$.  If a customer is accepted, $a_{l j}$ units of resource $l$ is consumed to serve a class $j$ customer; let {\color{black}$A_j=[a_{1j},\ldots,a_{mj}]^\top$} be the  column vector associated with class $j$ customers. Let $A\in \mathbb{R}^{m\times n}$ be the \emph{bill-of-materials} (BOM) matrix defined as $A=[A_1;\ldots;A_n]$}. 
If a customer is rejected, no revenue is collected and no resource is used.
Unused resources at the end of the horizon are perishable and have no salvage value.
The objective of the decision maker is to maximize the expected revenue earned during the entire horizon by deciding whether or not to accept each arriving customer. 

{\color{black}
For a control policy  $\pi$, 
let $z^{\pi}_j(t_1,t_2)$ be the number of class $j$ customers admitted during $(t_1,t_2]$ ($\forall j\in[n], 0\leq t_1 <t_2\leq T$) under that policy. We call a policy \emph{admissible} if it is non-anticipating and satisfies
\[
\sum_{j=1}^{n}A_jz_j^{\pi}(0,T)\leq C\, \text{ a.s.}, \quad\text{and}\quad
	 z_j^{\pi}(t_1,t_2) \leq \Lambda_j(t_1,t_2)\, \text{ a.s.}, \; \forall j\in[n], 0\leq t_1 < t_2 \leq T. \]
Let $\Pi$ be the set of all admissible policies. The expected revenue under policy $\pi\in \Pi$ is defined as $v^{\pi} = \E\big[\sum_{j=1}^{n}r_jz_j^{\pi}(0,T)\big]$.
We use $v^* = \sup_{\pi \in \Pi} v^{\pi}$ to denote the expected revenue under the optimal policy.
If $v^{\pi}$ is the expected revenue of a feasible policy $\pi \in \Pi$, we call $v^*-v^{\pi}$ the \emph{revenue loss} of policy $\pi$. 
}
\subsection{Asymptotic framework}
\label{subsec:performance}

The standard asymptotic framework in revenue management measures performance of heuristics when 
the capacities and customer arrivals are scaled up proportionally. Under this asymptotic scaling, we consider revenue loss of a sequence of problems, indexed by $k=1,2,\ldots$, where the capacities and arrival rates are multiplied by $k$, while all other problem parameters are treated as constants.

To avoid cumbersome notation where lots of variables and quantities are indexed by $k$, in the rest of the paper, we consider a different but equivalent asymptotic scaling, where the customer arrival rates $\lambda_j$ ($j\in[n]$) are kept as constants, the time horizon is scaled up by $T=1,2,\ldots$, and the resource capacities  are scaled up proportionally by $C_l = b_l T$ ($l\in[m]$).  Since the arrivals follow  Poisson processes, scaling up the arrival rates and scaling up the horizon length have the same effect.
We will thus express the revenue loss of heuristics in the order of $T$. Note that the horizon length $(T)$ plays the same role as the scaling factor $(k)$ in the standard asymptotic regime. For example, if we say the revenue loss of an algorithm is $O(\sqrt{T})$, it implies that revenue loss of that algorithm is $O(\sqrt{k})$ under the standard scaling regime.

\subsection{Previous work on upper bound approximations}	
\subsubsection{Deterministic linear program (DLP).} \label{subsec: DLP}

The DLP formulation is obtained by replacing all random variables with their expectations. As the expected number of arrivals of class $j$ customers during the horizon is $\lambda_jT$ for $j\in[n]$, the DLP formulation is given by
\begin{equation}
v^{\mathrm{DLP}}=\max_y \Big\{\sum_{j=1}^{n}r_jy_j\; \Big|\; \sum_{j=1}^{n}A_jy_j\leq C, \text{and } 0\leq y_j \leq \lambda_jT , \; \forall j \in[n]\Big\}. \label{form: DLP}
\end{equation}
In this formulation, decision variables $y_j$ can be viewed as the expected number of class $j$ customers to be accepted in $[0,T]$. The first constraint specifies that the expected usage of all $m$ resources  cannot exceed their initial capacities, $C =[C_1, \ldots,C_m]^\top$, and the second constraint specifies that the number of accepted customers from class $j$ cannot exceed the expected number of arrivals, $\lambda_j T$. 

Suppose $y^*$ is an optimal solution to \eqref{form: DLP}. The optimal value of DLP is given by $v^{\mathrm{DLP}}=\sum_{j=1}^nr_jy_j^*$. It can be shown that $v^{\mathrm{DLP}}$ is an upper bound of  the expected revenue of the optimal policy, $v^*$, namely $v^*\leq v^{\mathrm{DLP}}$ \citep{gallego1997multiproduct}. Intuitively, DLP is a relaxation of the original problem since it only requires the capacity constraints to be satisfied in expectation, so $v^{\mathrm{DLP}}$ is an upper bound of $v^*$.

Equivalently, we can reformulate the DLP in \eqref{form: DLP} by letting $x_j$ be the average number of class $j$ customers accepted per unit time, i.e., $x_j=y_j/T$. Then, we get
\begin{equation}
v^{\mathrm{DLP}}=\max_x\Big\{T\sum_{j=1}^{n}r_jx_j \;\Big|\; \sum_{j=1}^{n}A_jx_j\leq b, \text{ and } 0\leq x_j \leq \lambda_j , \; \forall j\in[n] \Big\}, \label{form: DLP2}
\end{equation}
where $b=[b_1,\ldots,b_m]^\top$ refers to the vector of available resources per unit time, i.e., $b_l=C_l/T, \forall l\in[m]$. 
Let $x_j^*$ for $j\in[n]$ be an optimal solution to \eqref{form: DLP2}.
The optimal value to the DLP is given by $v^{\mathrm{DLP}}=T\sum_{j=1}^nr_jx_j^*$. 

\subsubsection{Hindsight optimum.} \label{subsec: hindsight}

The hindsight optimum is the optimal revenue obtained when	the total number of arrivals is known in advance.
Recall that the random variable $\Lambda_j(T)$ represents the total arrivals of class $j$ customers in $[0,T]$. 
If the values of $\Lambda_j(T)$ are known, let $z_j$ be the number of class $j$ customers accepted in $[0,T]$; the optimal acceptance policy is given by
\begin{equation}
V^{\mathrm{HO}}=\max_y \Big\{\sum_{j=1}^{n}r_j z_j\; \Big|\; \sum_{j=1}^{n}A_j z_j\leq C, \text{ and } 0\leq z_j \leq \Lambda_j(T) , \; \forall j \in [n]\Big\}. \label{form: hindsight}
\end{equation}
Let   $V^{\mathrm{HO}}$ be the optimal objective value and $\bar{z}_j$, $j\in [n]$ be the optimal solution;
note that $V^{\mathrm{HO}}$ and $\bar{z}_j$'s are random variables that depend on $\Lambda_j(T)$.
The \emph{hindsight optimum} (HO) is defined as the expectation of the optimal objective value, i.e, $v^{\mathrm{HO}}=\E[V^{\mathrm{HO}}]=\E[\sum_{j=1}^n r_j \bar{z}_j]$.

The hindsight optimum is obviously an upper bound to the optimal revenue of the original problem, since the decision maker does not know the future arrivals at time $t=0$. In fact, it can be shown that hindsight optimum is a \emph{tighter} upper bound than the DLP, namely $v^*\leq v^{\mathrm{HO}} \leq v^{\mathrm{DLP}}$ \citep{talluri_analysis_1998}. This is easily verified since the expectation of the hindsight optimal solution, $\E[\bar{z}_j]$, is a feasible solution to the DLP.
{\color{black} 
We use the following definition throughout the paper.
\begin{definition}\label{def:regret}
Let $v^{\pi}$ be the expected revenue associated with an admissible control policy $\pi$. We refer to  $v^{\mathrm{HO}} - v^{\pi}$ as the \emph{regret} of that policy. (Note: since $v^*\leq v^{\mathrm{HO}}$, the revenue loss of the control policy, $v^* - v^{\pi}$, is upper  bounded by its regret.)
\end{definition}
}



\subsection{Static probabilistic allocation heuristic}
\label{subsec: FP}

There are various ways to construct heuristic policies using the optimal solution of DLP. An overview can be found in \citet[Ch. 2]{talluri_theory_2004}.
One intuitive approach is to interpret the solution to DLP as acceptance probabilities. 
Suppose $x^*$ is an optimal solution to DLP in \eqref{form: DLP2}. For each arriving customer, if the customer belongs to class $j$, s/he would be accepted independently with probability $x^*_j/\lambda_j$ throughout the time horizon. Since customers from each class are accepted with probabilities that are static, we call this heuristic Static Probabilistic Allocation (\textsf{SPA}). The \textsf{SPA} policy is formally stated in Algorithm~\ref{algo:FP}.

{\SingleSpacedXI
\begin{algorithm} 
	\begin{algorithmic}
		\State {\bf initialize} $x^* \leftarrow \arg\max_x \Big\{\sum_{j=1}^{n}r_jx_j \;\Big|\; \sum_{j=1}^{n}A_jx_j\leq C/T, \text{ and }\, 0\leq x_j \leq \lambda_j, \forall  j \in [n] \Big\}$; $C'\leftarrow C$
		\For {all customers arriving in $[0,T]$}
		\If {the customer belongs to class $j$ and $A_j \leq C'$ ($\forall j\in[n]$)}
		\State accept the customer with probability $x^*_j/\lambda_j$
		\State if the customer is accepted, update capacity $C'\leftarrow C' - A_j$
		\Else 
		\State reject the customer
		\EndIf
		\EndFor
	\end{algorithmic}
	\caption{\textsf{Static probabilistic allocation heuristic: \textsf{SPA} }\label{algo:FP}}
\end{algorithm}
}

The expected revenue of the \textsf{SPA} policy, denoted by $v^{\textsf{SPA}}$, can be computed as follows. 
Since the total number of arrivals from class $j$ follows a Poisson distribution with mean $\lambda_j T$,
the number of customers that the algorithm \emph{attempts} to accept from class $j$ follows a Poisson distribution with mean $(\lambda_j T)\cdot x^*_j /\lambda_j=x^*_j T=y^*_j$. Due to limited capacity, we must reject any customer from class $j$ if the remaining capacity $C'$ does not satisfy $A_j\leq C'$. It is straightforward to show that the expected number of customers who are turned away due to capacity limits is $O(\sqrt{T})$ \citep[see e.g.][]{gallego1997multiproduct,ReimanAsymptoticallyOptimalPolicy2008}. Thus, we have
\[
v^{\textsf{SPA}}=\sum_{j=1}^n r_j y^*_j - O(\sqrt{T}) = v^{\mathrm{DLP}} - O(\sqrt{T}).
\]
Recall from  \S\ref{subsec: DLP} that $v^{\mathrm{DLP}}$ is an upper bound of the expected revenue under the optimal policy, namely $v^*\leq v^{\mathrm{DLP}}$. 
Thus, the revenue of \textsf{SPA} is bounded by $v^{\textsf{SPA}}\geq v^* - O(\sqrt{T})$. 


\section{Frequent Re-solving and Degeneracy} 
\label{sec:FR}

An obvious drawback of the \textsf{SPA} policy constructed from the DLP is that it does not take into account the randomness of demand or the updated information after $t=0$.
This motivates us to consider \emph{re-solving heuristics}, which periodically re-optimize the DLP using the updated capacity information to adjust customer admission controls. 

In particular, the following re-solving heuristic, which we referred to as \textsf{Frequent Re-solving (FR)}, has been studied by \citet{JasinReSolvingHeuristicBounded2012} and \citet{wu2015algorithms}.
The \textsf{FR} policy divides the horizon into $T$ periods and  re-solves the LP at the beginning of each period.	At time $t=0,1,\ldots,T-1$, let $C_l(t)$ denote the remaining capacity of resource $l \in [m]$. 
We let $b_l (t):= \frac{C_l (t)}{T-t}$ be the average available capacity of resource $l$ in period $t$. Let $C(t)$ and $b(t)$ denote the vectors of the remaining capacities and the average remaining capacities per unit time at time $t$, respectively, for all the resources. We outline the \textsf{FR} policy in Algorithm~\ref{algo:res}.

{\SingleSpacedXI
\begin{algorithm} 
	\begin{algorithmic}
		\State {\bf initialize}: set $C(0)=C$ and $b(0)=C/T$
		\For {$t = 0, 1, \ldots, T-1$}
		\State set $x(t) \leftarrow \arg\max_x \Big\{\sum_{j=1}^{n}r_jx_j \;\Big|\; \sum_{j=1}^{n}A_jx_j\leq b(t), \text{ and }\, 0\leq x_j \leq \lambda_j, \forall  j \in [n] \Big\}$
		\State set $C' \leftarrow C(t)$
		\For {all customers arriving in $[t,t+1)$}
		\If {the customer belongs to class $j$ and $A_j \leq C'$ ($\forall j\in[n]$)}
		\State accept the customer with probability $x_j(t)/\lambda_j$
		\State if the customer is accepted, update $C'\leftarrow C' - A_j$
		\Else 
		\State reject the customer
		\EndIf
		\EndFor
		\State set $C(t+1)\leftarrow C'$ and $b(t+1)\leftarrow \frac{C(t+1)}{T-t-1}$
		\EndFor
	\end{algorithmic}
	\caption{\textsf{Frequent Re-solving Heuristic: \textsf{FR} }\label{algo:res}}
\end{algorithm}
}

\citet{JasinReSolvingHeuristicBounded2012} show that when the optimal solution to DLP \eqref{form: DLP2}
is \emph{nondegenerate}, \textsf{FR} has a revenue loss of $O(1)$, namely, the revenue loss is bounded when the problem size $k$ grows.  The optimal solution $x^*$ is nondegenerate if 
\begin{equation} \label{eq:non-degeneracy}
\left|\{j\in[n]: x^*_j =0 \text{ or } x^*_j = \lambda_j\}\right|+ \bigl|\{l\in[m]: \sum_{j=1}^n a_{lj} x^*_j = b_j\} \bigr|= n.
\end{equation}
The $O(1)$ loss is a significant improvement from the $	O(\sqrt{T})$ revenue loss of \textsf{SPA}. 
{\color{black}
However, the assumption of nondegenerate DLP solution is critical to achieve the $O(1)$ loss.
The proofs by \cite{JasinReSolvingHeuristicBounded2012} and \cite{wu2015algorithms} are built on a key observation that the ratio of remaining capacities to remaining time, $b(t)$, is a martingale (see also \citet{arlotto2017uniformly} for a discussion on this martingale property).
If the optimal solution $x^*$ is safely far from any degenerate solutions, with high probability, the adjusted solution $x(t)$ in Algorithm~\ref{algo:res} shares the same basis with $x^*$, so the revenue loss of \textsf{FR} can be bounded.
It is unclear from the analysis of \citet{JasinReSolvingHeuristicBounded2012} and 
\citet{wu2015algorithms} 
whether the nondegeneracy assumption is just an artifact of their analysis technique or something intrinsic to the performance of \textsf{FR}.
This motivates us to examine closely the role of the nondegeneracy assumption.
}

\subsection{A degenerate example}
\label{subsec:degenerate example}
{\color{black} 
We will illustrate the issue of degenerate DLP solutions using the following numerical example, while deferring the theoretical analysis of the \textsf{FR} policy to Section~\ref{sec: analysis of resolving}.
}

Suppose there are two classes of customers and one resource. Customers from each class arrive according to a Poisson process with rate 1. Customers from both classes, if accepted, consume one unit of resource, but pay different prices, $r_1$ and $r_2$.
First, we  compare the expected revenue loss of the \textsf{FR} policy and the \textsf{SPA} policy, which does not re-solve after $t=0$, to examine the effect of frequent re-solving. We simulate the \textsf{FR} policy and the \textsf{SPA} policy when the average capacity per unit time $b=1$ (so the total capacity is $T$) for two price scenarios: (a) $r_1=2$ and $r_2=1$; (b) $r_1=5$ and $r_2=1$ and for varying horizon length $T=500,\ldots,5000$. In both scenarios, the optimal solution to the DLP \eqref{form: DLP2} is $x^*_1 = 1, x^*_2 = 0$. {\color{black} From Equation~\eqref{eq:non-degeneracy}, we have
\[
\left|\{j\in[n]: x^*_j =0 \text{ or } x^*_j = \lambda_j\}\right|+ \bigl|\{l\in[m]: \sum_{j=1}^n a_{lj} x^*_j = b_j\} \bigr|= 3 > n = 2,
\]
thus the DLP solution in this example is degenerate.} 

\begin{figure}[!htp]
	\centering
	\subfloat[$r_1=2$ and $r_2=1$]{%
		\includegraphics[clip,width=0.46\columnwidth]{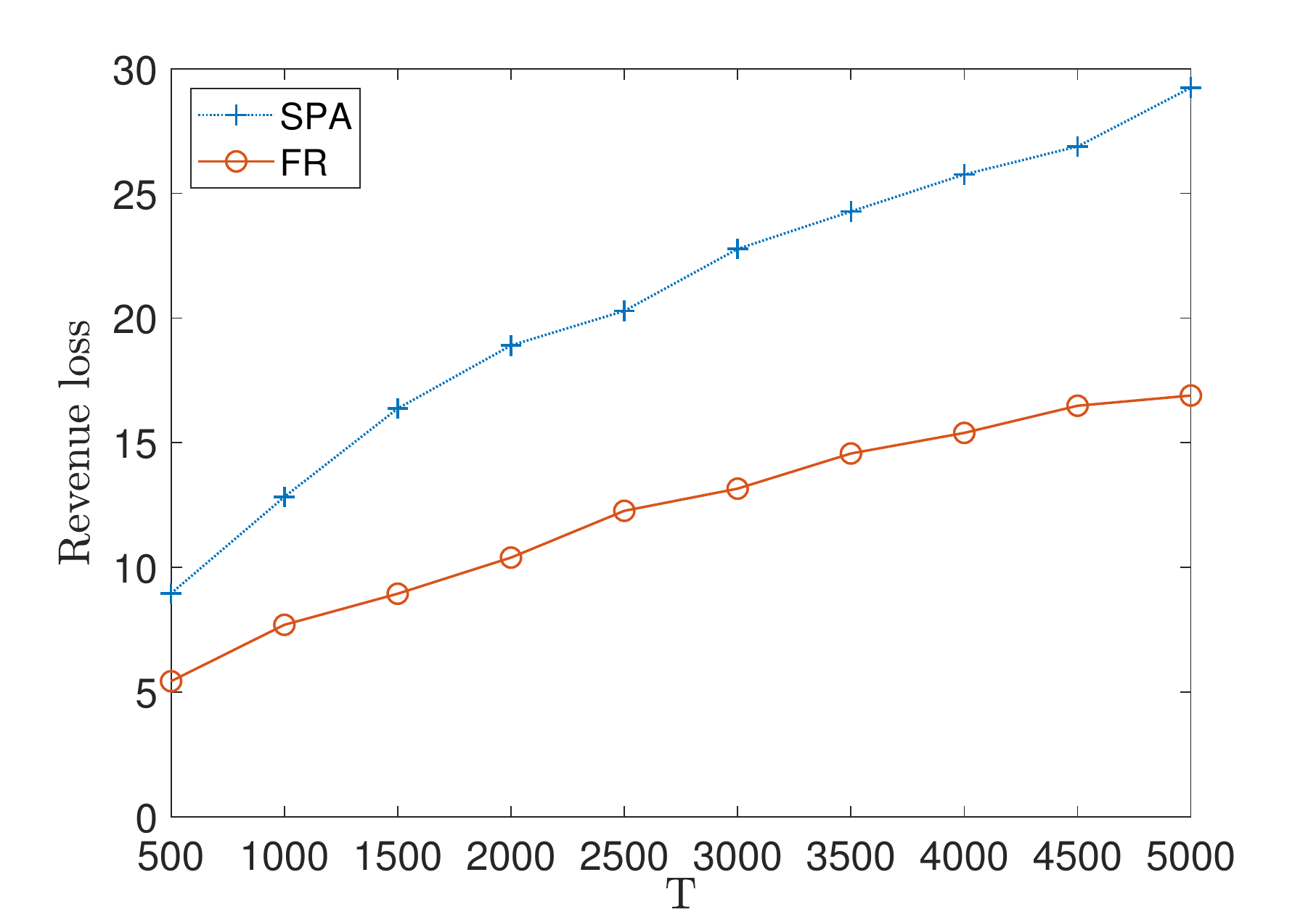}%
	}
	\subfloat[$r_1=5$ and $r_2=1$]{%
		\includegraphics[clip,width=0.46\columnwidth]{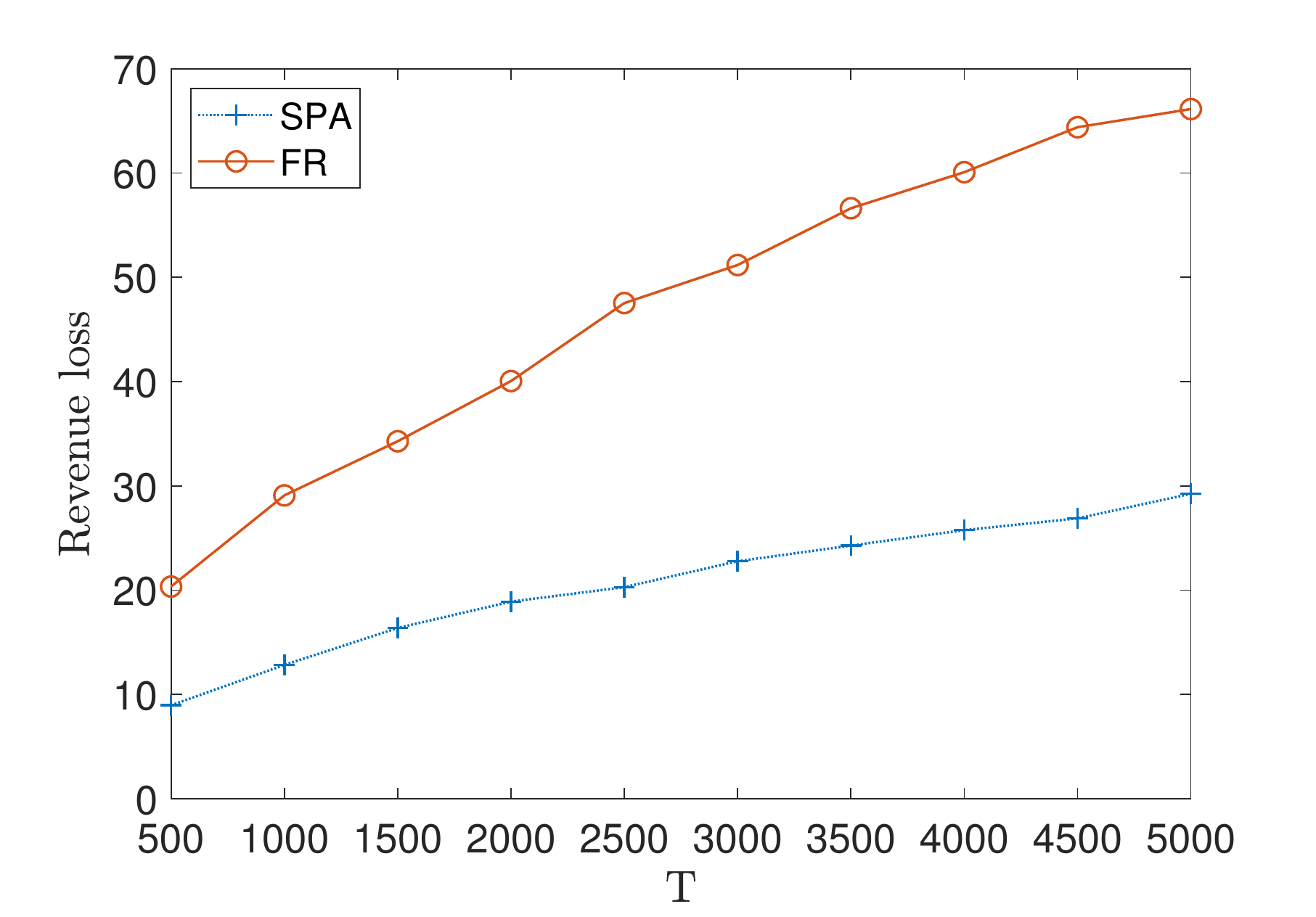}%
	}
	\caption{Regret under the \textsf{FR} policy (with re-solving) and the \textsf{SPA} policy (without re-solving).}
	\label{fig: resolve vs fp}
\end{figure}

Recall that the expected revenue loss of the \textsf{FR} policy is defined as $v^*-v^{\textsf{FR}}$. Since calculating $v^*$ requires solving dynamic programs, we use the regret $v^{\mathrm{HO}}-v^{\color{black}\mathsf{FR}}$ (see the definition in \S\ref{subsec: hindsight}) as a proxy of the expected revenue loss. In \S\ref{sec:new-heuristic}, we will show that $v^{\mathrm{HO}}-v^* = O(1)$, so this substitution does not affect the rate of revenue loss. Fig.~\ref{fig: resolve vs fp} plots the average revenue losses under the \textsf{FR} policy and the \textsf{SPA} policy over 1000 sample paths.

We make the following observations from Fig.~\ref{fig: resolve vs fp}.
First, while the revenue loss of  \textsf{FR} in scenario (a)  is lower than that obtained from applying the \textsf{SPA} policy, the relationship is reversed in scenario (b). In other words, re-solving the DLP does not always lead to better performance. The intuition behind this result is that when the ratio $r_1 / r_2$ is large, such as in scenario (b), rejecting a Class 2 customer to save the capacity for a potential future Class 1 customer is more profitable.
The \textsf{SPA} policy accepts every customer from Class 1 and rejects all customers from Class 2, since the solution to the DLP (without re-solving) is $x^*_1 = 1, x^*_2 = 0$. This static policy is indeed optimal when $r_1 / r_2 \to \infty$.
In contrast, the \textsf{FR} policy constantly adjusts accepting probabilities, and starts to accept Class 2 customers when the actual arrival of Class 1 customers falls below its average. 
Second, we observe from Fig.~\ref{fig: resolve vs fp} that the revenue losses of both \textsf{SPA} and \textsf{FR} seem to have the same growth rate as horizon length $T$ increases. (It is well-known that the revenue loss of \textsf{SPA} is of order $\Theta(\sqrt{T})$; see \S\ref{subsec: FP} and Proposition~\ref{prop: DLP and FPA} in Appendix~\ref{append:additional-result}.) This result is in contrast with
 \citet{JasinReSolvingHeuristicBounded2012}, which show that when the solution to the DLP is nondegenerate, the expected revenue loss of \textsf{FR} is $O(1)$. However, we note that the nondegeneracy assumption made by Jasin and Kumar does not hold in this example, since the DLP has a unique solution that is degenerate.

{\color{black} 
Next, we 
simulate the \textsf{FR} policy when $r_1=2$, $r_2=1$ and $T=5000$ for varying average capacity per unit time $b=0.5,\ldots,2$. Note that when $b=1$ and $b=2$, the optimal solutions to the DLP \eqref{form: DLP2} are $x^*_1 = 1, x^*_2 = 0$ and $x^*_1 = 1, x^*_2 = 1$, respectively, which are degenerate according to Equation \eqref{eq:non-degeneracy}. When $b\neq 1$ and $b\neq 2$, the solution to the DLP is nondegenerate.  Therefore, by changing the value of $b$, we can evaluate the performance of \textsf{FR} with either degenerate or nondegenerate DLP solutions.
Fig.~\ref{fig:single fr vary c} shows the average revenue loss under the \textsf{FR} policy over 1000 sample paths.

\begin{figure}[!htp]
	\centering
	\includegraphics[clip,width=0.7\columnwidth]{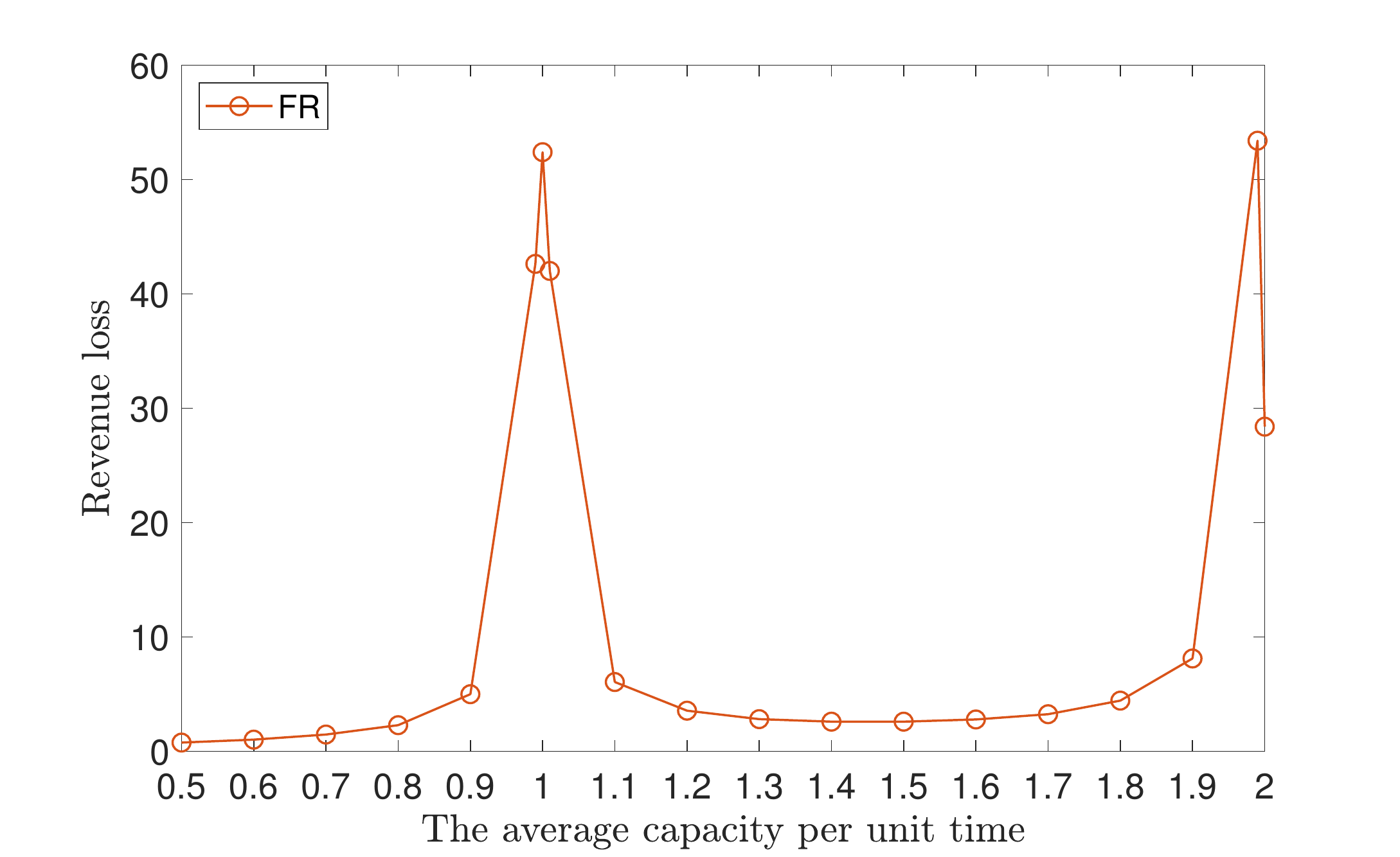}
	\caption{Regret under the \textsf{FR} policy  for $r_1=2, r_2=1$ and $T=5000$.}
	\label{fig:single fr vary c}
\end{figure}

The simulation result from Fig.~\ref{fig:single fr vary c} shows that the expected revenue loss under the \textsf{FR} policy is sensitive to the value of capacity rate $b$. When $b$ is far away from the degenerate points (i.e., $b=1$ and $b=2$),  \textsf{FR} performs well and has small revenue loss. However, the revenue loss increases significantly when the optimal DLP solution is close to degenerate (e.g., $b=0.95$).

We notice that the observation from  Fig.~\ref{fig:single fr vary c} is consistent with the analysis by \cite{JasinReSolvingHeuristicBounded2012}.
Even though \cite{JasinReSolvingHeuristicBounded2012} proves that the revenue loss of \textsf{FR} is bounded by a constant whenever the DLP solution is nondegenerate, their analysis \emph{does not imply the constant is  uniform over all $b$'s}. Rather, the constant bound from their analysis critically depends on the distance between $b$ and its nearest degenerate point. When the optimal DLP solution is close to degenerate, the bound in \cite{JasinReSolvingHeuristicBounded2012} can be arbitrarily large.  Fig.~\ref{fig:single fr vary c} shows that this phenomenon is not merely a consequence of the analysis technique from \cite{JasinReSolvingHeuristicBounded2012}, but reflects the actual performance of the \textsf{FR} policy.
}

\section{A Re-solving Heuristic with Uniformly Bounded Loss}
\label{sec:new-heuristic}

In this section, we propose a new re-solving algorithm. 
The main result of this section is to show that this algorithm has uniformly bounded revenue loss given any horizon length $T$ and starting capacity $C$, without requiring the nondegeneracy assumption.

\subsection{Definition of the \textsf{IRT} algorithm}\label{sub: TR}

We propose an algorithm called   \textsf{Infrequent Re-solving and Thresholding (IRT)}. 
The \textsf{IRT} policy has two distinct features compared to the \textsf{FR} policy: 1) the DLP is \emph{not} re-solved in every period; 2) customers acceptance probabilities are adjusted by some thresholds. 

Unlike the $\mathsf{FR}$ policy, the \textsf{IRT} policy re-solves the DLP for only $O(\log\log T)$ times during a horizon of length $T$.
The re-solving schedule is defined as follows.
Given horizon length $T$, we set 
$
K=\left\lceil\frac{\log \log T}{\log (6/5)}\right\rceil.
$
Let $\{t^*_u, \forall u\in [K] \}$ denote a sequence of re-solving times, where $\tau_u = T^{(5/6)^u}$ and $t^*_u = T - \tau_u$ for all $u \in [K]$. {\color{black} In addition, let $t^*_{K+1}=T$. Thus, the re-solving times divide the entire horizon into $K+1$ epochs:
$[0, t^*_1)$, $[t^*_1, t^*_2)$, $\cdots$, $[t^*_K, t^*_{K+1}]$. }
Fig.~\ref{fig:recursive} illustrates the re-solving schedule of the \textsf{IRT} policy.
{
\begin{figure}[!htb]
\centering
	\includegraphics[width=0.7\columnwidth]{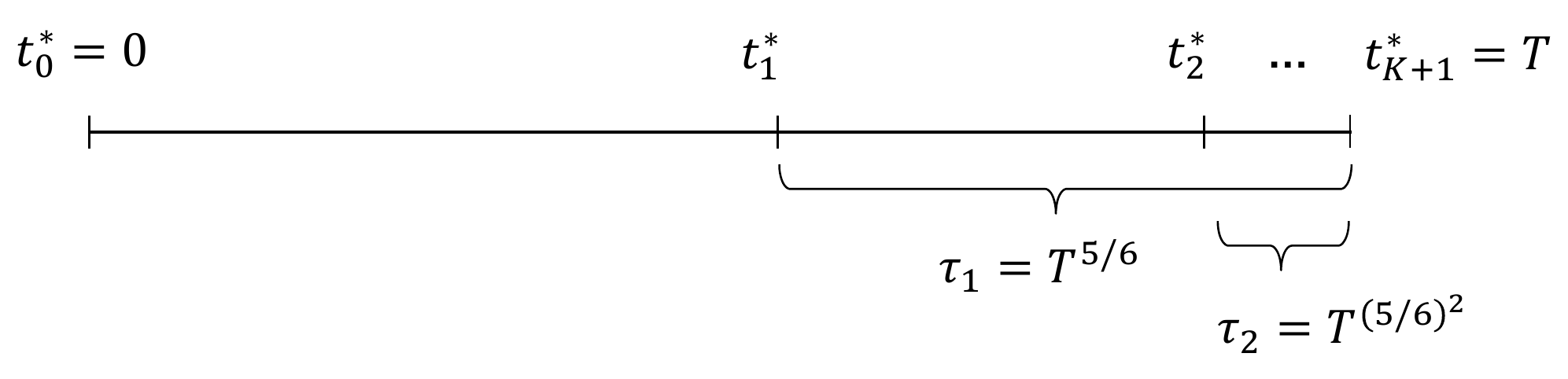}%
	\caption{The re-solving times of the \textsf{IRT} policy is constructed recursively.}\label{fig:recursive}
\end{figure}
}

At the beginning of each epoch $u$ $(0 \leq u \leq K)$, the algorithm solves an LP approximation to the dynamic programming problem---this LP is identical to the LP used in the \textsf{FR} algorithm, which uses information about remaining capacities and the mean of remaining customer arrivals.
The optimal solution of the LP is then used to construct a probabilistic allocation control policy. {\color{black} The  \textsf{IRT} policy applies thresholds to the allocation probabilities. In particular, in epoch $u \in \{0\}\cup [K-1]$ (except for the last epoch), the allocation probability for each class is rounded down to 0 if it is less than $\tau^{-1/4}_u$, or  rounded up to 1 if it is larger than $1-\tau^{-1/4}_u$.} 
The complete definition of \textsf{IRT} is given in Algorithm~\ref{algo:less-is-more}.

{\SingleSpacedXI
\begin{algorithm} 
	\begin{algorithmic}
		\State {\bf initialize}: set $\tau_u=T^{(5/6)^u}$ and $t^*_u = T - \tau_u$ for all $u\in\{0\}\cup [K]$, where 
		$
		K=\left\lceil\frac{\log \log T}{\log (6/5)}\right\rceil
		$
		\For {$u = 0, 1, \ldots, K$}
		\State set $x^u \leftarrow \arg\max_x \Big\{\sum_{j=1}^{n}r_jx_j \;\Big|\; \sum_{j=1}^{n}A_jx_j\leq C(t^*_k)/{\tau_k}, \text{ and }\, 0\leq x_j \leq \lambda_j, \forall  j \in [n] \Big\}$
		\If {$u<K$}
		\For {$j \in [n]$}
		\If {$x^u_j < \lambda_j \tau_u^{-1/4}$}
		\State set $p^u_j \leftarrow 0$
		\ElsIf {$x^u_j >\lambda_j (1 - \tau_u^{-1/4})$}
		\State set $p^u_j \leftarrow 1$
		\Else
		\State set $p^u_j \leftarrow x^u_j / \lambda_j$ 
		\EndIf
		\EndFor
		\Else
		\State set $p^u_j \leftarrow x^u_j / \lambda_j$ for all $j\in[n]$
		\EndIf
		\State set $C' \leftarrow C(t^*_u)$
		\For {$t\in [t_u^*,t_{u+1}^*)$}
		\State observe requests from all arrival of customers
		\If {an arriving customer belongs to class $j$ and $A_j \leq C'$ ($\forall j\in[n]$)}
		\State accept the customer with probability $p^u_j$
		\State if accepted, update $C'\leftarrow C' - A_j$ 
		\Else 
		\State reject the customer
		\EndIf
		\EndFor
		\State set $C(t^*_{u+1}) \leftarrow C'$
		\EndFor
	\end{algorithmic}
	\caption{\textsf{Infrequent Re-solving with Thresholding (IRT)}\label{algo:less-is-more}}
\end{algorithm}
}

{\color{black}
Before we present the formal analysis of the \textsf{IRT} algorithm,
it might be helpful to discuss the intuition behind the design of this algorithm. We start with the choice of the first re-solving time, $t^*_1$. The analysis by \cite{ReimanAsymptoticallyOptimalPolicy2008} shows that by setting $t^*_1 \approx T - O(\sqrt{T})$, one re-solving of DLP is sufficient to reduce the regret to $O(T^{1/4})$. But we note that if $t^*_1$ is defined as in \cite{ReimanAsymptoticallyOptimalPolicy2008}, additional re-optimizations after $t^*_1$ cannot improve the regret rate.
In the \textsf{IRT} algorithm, we choose the first re-solving time to be $t^*_1 = T - T^{5/6}$, which is earlier than the re-solving time in \cite{ReimanAsymptoticallyOptimalPolicy2008}. If no further re-solving is used, this policy leads to a regret rate of $O(T^{5/12})$ (Proposition~\ref{prop: resolve once}). 
Even though the $O(T^{5/12})$ rate is worse than the $O(T^{1/4})$ rate in \cite{ReimanAsymptoticallyOptimalPolicy2008}, as we choose an earlier re-solving time, more time is left for making further adjustments. Once we establish the  $O(T^{5/12})$ regret rate with the first re-solving, we then use induction to prove that subsequent re-optimizations of the DLP can further reduce the regret, eventually reducing it to a constant. 
By definition, $\tau_u$, the length of epoch $u$ satisfies the recursive relationship $\tau_{u+1} = \tau_{u}^{5/6}$, $\forall u \in [K]$. This enables us to apply the induction hypothesis to epochs $u \geq 1$. 

The $\tau_u^{-1/4}$ thresholds in the algorithm are critical to bounding the regret. As we have seen from the numerical example in \S\ref{subsec:degenerate example}, large losses can occur when the DLP solution is nearly degenerate. If we use a nearly degenerate solution to construct probabilistic allocation controls, some customer classes would have acceptance probabilities that are either very close to 0 or very close 1. As a result, the mean number of accepted or rejected customers is dominated by its standard deviation, making the control policy ineffective. 
More specifically, if the acceptance probability of class $j$ customers is $\epsilon \to 0$, the coefficient of variation of the number of customer accepted in one unit time is $1/\sqrt{\lambda_j \epsilon} \to +\infty$. Therefore, if the acceptance probability of a customer class is almost 0, we might as well reject all customers from that class in the current epoch, as long as there is sufficient time left to accept customers in the next epoch. Similarly, if the acceptance probability of a customer class is almost 1,  
we might as well accept all customers from that class in the current epoch. This is the intuition behind adding thresholds to the acceptance probabilities in the \textsf{IRT} policy.




}

%

\subsection{Analysis of the \textsf{IRT} policy}
We now formally analyze the revenue loss (regret) of the \textsf{IRT} policy. The main result of this section is the following.
{\color{black}
\begin{theorem} \label{thm: TR}
The regret of \textsf{IRT} policy define in Algorithm~\ref{algo:less-is-more} 
is bounded by $
	v^{\mathrm{HO}} - v^{\textup{\textsf{IRT}}} =  O(1).$
	The constant factor depends on the customer arrival rate $\lambda_j$ ($\forall j\in[n]$), the revenues per customer $r_j$ ($\forall j\in[n]$), and the BOM matrix $A$; however, this constant is
	independent of the time horizon $T$ and the capacity vector $C$.
\end{theorem}
}
Theorem~\ref{thm: TR} states that the regret of \textsf{IRT} policy is $O(1)$. {\color{black} Moreover, this constant is independent of time horizon and capacities, so the performance of \textsf{IRT} is uniformly bounded when the capacity ratio $C/T$ varies. Because degenerate DLP solution occurs only for some specific capacity ratios, the result in Theorem~\ref{thm: TR} does not require the nondegeneracy assumption in \cite{JasinReSolvingHeuristicBounded2012}. 
}
Since the hindsight optimum $v^{\mathrm{HO}}$ is an upper bound of the expected revenue of the optimal policy $v^*$, we immediately get a bound on its revenue loss:
$
v^* -  v^{\textsf{IRT}} \leq v^{\mathrm{HO}} - v^{\textsf{IRT}} = O(1).
$
Moreover, Theorem~\ref{thm: TR} implied that hindsight optimum is a tight upper bound, satisfying
$
v^{\mathrm{HO}} -  v^* \leq v^{\mathrm{HO}} - v^{\textsf{IRT}}  = O(1).
$

The complete proof of Theorem~\ref{thm: TR} can be found in Appendix~\S\ref{subsec:proof-TR}.
We outline the main idea of the proof here.
In the proof, we define a sequence of auxiliary re-solving policies with increasing re-solving frequency. Recall that $K=\left\lceil\frac{\log \log T }{\log (6/5)}\right\rceil
$ is the number of re-optimizations made by the \textsf{IRT} algorithm. 
For any $u \in [K]$, we define a policy that follows the \textsf{IRT} heuristic exactly in $[0, t_u^*)$, but then applies static allocation control in $[t_u^*, T]$.
We refer to such a policy as
$\mathsf{IRT}^u$. 
{\color{black} Notice that when $u=K$, $\mathsf{IRT}^u$ coincides with \textsf{IRT}. }
Similarly, we define $\mathsf{HO}^u$ as a policy that  is exactly the same as \textsf{IRT} in $[0,t^*_u)$ but applies the hindsight optimal policy in $[t^*_u,T]$. 
Our proof of Theorem~\ref{thm: TR} depends on the following proposition, proved in Appendix \S\ref{sub: proof_proposition}.
\begin{proposition}\label{prop: resolve once}
	Given horizon length $T$, suppose the first re-solving time is  $t_1^*=T-T^{5/6}$, then
	\begin{enumerate}
		\item the regret of $\mathsf{HO}^1$ is $O(T e^{-\kappa T^{1/6}})$;
		\item the regret of $\mathsf{IRT}^1$ is $O(T e^{-\kappa T^{1/6}}) + O(T^{5/12})$.
	\end{enumerate}
Here, we define $\kappa = 	 \frac{\lambda_{min}}{27(\alpha|J_\lambda|+1)^2}$, where $J_\lambda：=\{j: x^*_j=\lambda_j\}$ (recall that $x^*$ is the solution to DLP),
	 $\lambda_{\min}：=\min_{j\in[n]}\lambda_j$ , and $\alpha$ is a positive constant that depends on the BOM matrix $A$.
\end{proposition}
Notice that $\mathsf{IRT}^1$ is a non-anticipating and admissible policy, and its regret of $O(T^{5/12})$ is an improvement over the $O(\sqrt{T})$ bound of \textsf{SPA}. The policy $\mathsf{HO}^1$ is not non-anticipating since it requires access to future arrival information; thus it is not practical and its sole purpose is to bound the performance of  $\mathsf{IRT}^1$ in the proof.

We then use Proposition~\ref{prop: resolve once} to prove Theorem~\ref{thm: TR} by induction. We illustrate the induction step using $\mathsf{IRT}^2$, a policy that re-solves at $t^*_1 = T-T^{5/6}$ and  again at $t_2^*=T-T^{(5/6)^2}$.
	The regret of $\mathsf{IRT}^2$ can be written as
\[
\E[V^{\mathrm{HO}}-V^{\mathsf{IRT}^2}]\nonumber=\underbrace{\E[V^{\mathrm{HO}}-V^{\mathsf{HO}^1}]}_{(*)}+\underbrace{\E[V^{\mathsf{HO}^1}-V^{\mathsf{HO}^2}]}_{(**)}+\underbrace{\E[V^{\mathsf{HO}^2}-V^{\mathsf{IRT}^2}]}_{(***)}.
\]
	The term $(*)$ is bounded  by $O(Te^{-\kappa T^{1/6}})$ according to Proposition~\ref{prop: resolve once}. For the term $(**)$, the policies $\mathsf{HO}^1$ and $\mathsf{HO}^2$ are identical up to time $t^*_1$. So applying 
part (1) of Proposition \ref{prop: resolve once} to the subproblem in $(t^*_1,T]$, we get
$
	\E[V^{\mathsf{HO}^1}-V^{\mathsf{HO}^2}]
	=O(T^{5/6} e^{-\kappa (T^{5/6})^{1/6}})
	= O(T^{5/6} e^{-\kappa T^{5/36}}).
$
For the last term $(***)$, using the well-known result that static probabilistic allocation has a squared root regret, we have
$\E[V^{\mathsf{HO}^2}-V^{\mathsf{IRT}^2}]=
\E[V^{\mathrm{HO}}(t_2^*,T)-V^{\mathsf{SPA}}(t_2^*,T)]
	= O(\sqrt{T-t_2^*})
	= O(T^{(5/6)^2/2}).$
	Combining these three terms, we get
\[
	v^{\mathrm{HO}}-v^{\mathsf{IRT}^2}=O(Te^{-\kappa T^{1/6}})+O( T^{5/6} e^{-\kappa T^{5/36}})+O(T^{25/72})= O(T^{25/72}).
\]

By induction, we show that if the decision maker re-solves for $K \geq 1$ times, where the $u$-th ($u=1,\cdots,K$) re-solving time is $t^*_u = T - T^{(5/6)^u}$, the regret is given by
\[
	v^{\mathrm{HO}}-v^{\mathsf{IRT}^K}=\sum_{u=0}^{K-1} 
	O\left((T^{(5/6)^u} \exp\left(-\kappa T^{(5/6)^u/6}\right)\right) 
	+ O( T^{(5/6)^K/2}).
\]
When
$K=\left\lceil\frac{\log \log T}{\log (6/5)}\right\rceil,$
the right-hand side of the above equation is bounded by a constant. In additional, the policy $\mathsf{IRT}^K$ is the same as $\mathsf{IRT}$, so we
 prove that the regret of \textsf{IRT} is $v^{\mathrm{HO}}-v^{\mathsf{IRT}} =O(1)$. 

\subsection{Revisiting the degenerate example in Section~\ref{subsec:degenerate example}}

{\color{black}
In Section~\ref{subsec:degenerate example}, we considered a numerical example with two classes and one resource.
We 
simulated the \textsf{FR} policy when $r_1=2$, $r_2=1$ and $T=5000$ for varying average capacity  $b=0.5,\ldots,2$, and showed that \textsf{FR} has poor performance when the DLP solution is either degenerate (i.e., $b=1$ or $b=2$) or nearly degenerate. 
We now test the  \textsf{IRT} policy using the same example and compare it to the \textsf{FR} policy.  Fig.~\ref{fig:single fr and irt vary c} plots the average 
regret under \textsf{FR} and \textsf{IRT} over 1000 sample paths.
	\begin{figure}[!htp]
		\centering
		\includegraphics[clip,width=0.7\columnwidth]{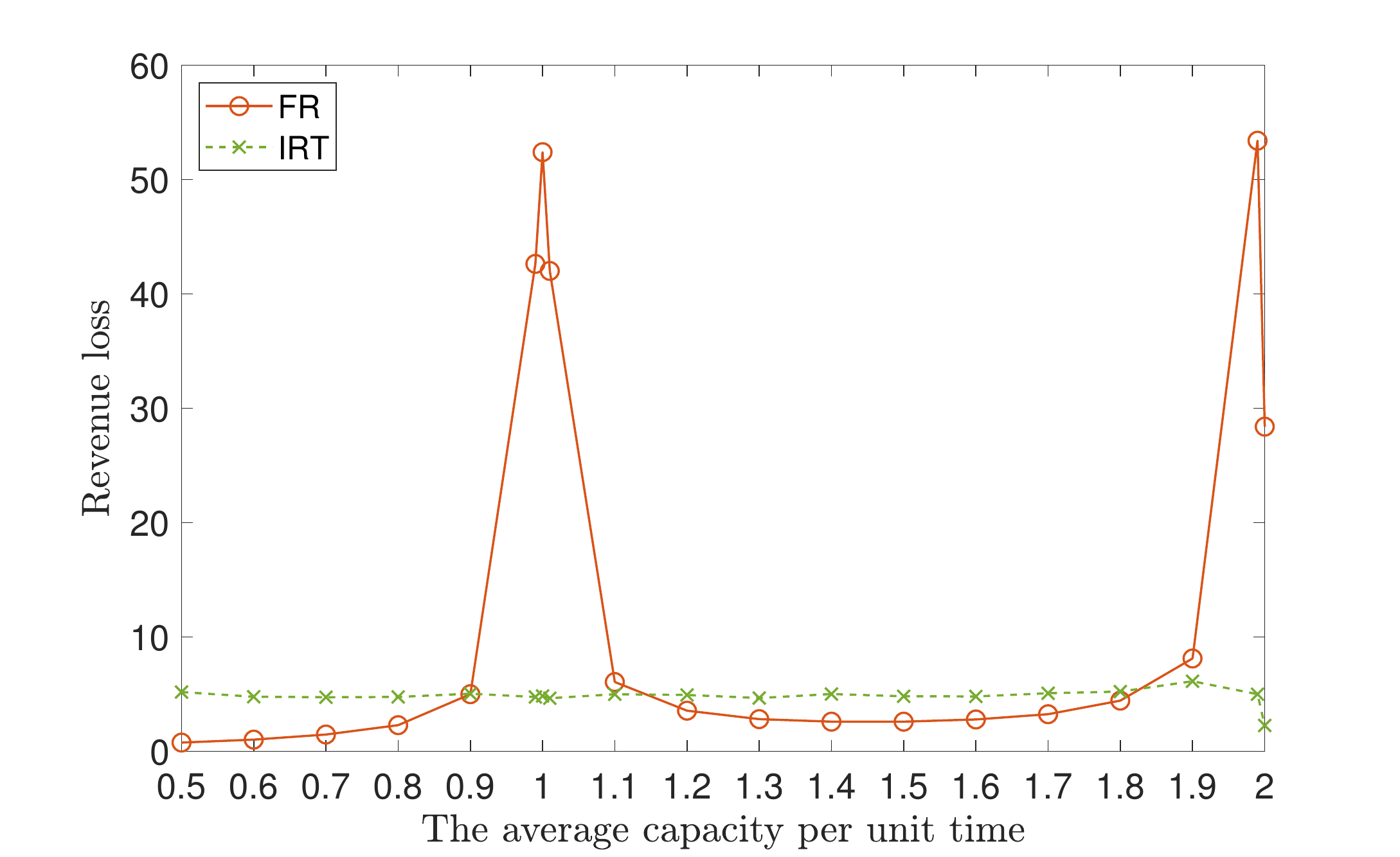}
		\caption{Regret under the \textsf{FR} policy and the \textsf{IRT} policy for $r_1=2, r_2=1$ and $T=5000$.}
		\label{fig:single fr and irt vary c}
	\end{figure}
	
	It can be observed from Fig.~\ref{fig:single fr and irt vary c} that  the regret under the proposed \textsf{IRT} policy is not sensitive to the average capacity per unit time. 
This result verifies Theorem~\ref{thm: TR} in that the regret of \textsf{IRT} is uniformly bounded with respect to the ratio between capacity and time.
	In contrast, the regret under the \textsf{FR} policy has two spikes that are associated with the two degenerate points ($b=1$ and $b=2$). 

}

\section{Analysis of the Frequent Resolving Policy}
\label{sec: analysis of resolving}

\subsection{Lower bound of the revenue loss of  \textsf{FR}}
\label{subsec:HO-FR}

The simulation in Section~\ref{subsec:degenerate example} inspires us to analyze the performance of \textsf{FR}  without the nondegeneracy assumption in order to gain a better understanding of the effect of frequent re-solving.
First, we show that the regret under the \textsf{FR} policy is bounded below by $\Omega(\sqrt{T})$. 
%
%

\begin{proposition}\label{thm: ho res lb}
	There exists a problem instance for which the regret of the \textsf{FR} policy defined in Algorithm~\ref{algo:res} is bounded below by
	\[
	v^{\mathrm{HO}}-v^\mathsf{FR}=\Omega(\sqrt{T}).
	\]
\end{proposition}

Proposition~\ref{thm: ho res lb} implies that the expected revenue loss under \textsf{FR} policy is bounded below by $\Omega(\sqrt{T})$ as well, because the revenue gap between the hindsight optimum ($v^{\mathrm{HO}}$) and the optimal revenue ($v^*$) is $O(1)$ (Theorem~\ref{thm: TR}).
That is, we have
\begin{align*}
v^*-v^\textsf{FR}=-(v^{\mathrm{HO}}-v^*)+v^{\mathrm{HO}}-v^\textsf{FR}= -O(1)+\Omega(\sqrt{T})=\Omega(\sqrt{T}).
\end{align*}

To prove Proposition~\ref{thm: ho res lb}, we consider a problem instance with two classes of customers and one resource. We assume that customers from each class arrive according to a Poisson process with rate 1; the arrivals from two classes are independent.
The initial resource capacity is $T$. 
Customers from both classes, if accepted, consume one unit of the resource, but pay different prices, $r_1>r_2$. 
We consider the event when the number of class 1 customers that arrive during $T$ period is more than $T$. If this event happens, the hindsight optimum will accept $T$ of class 1 customers and none of class 2 customers.
Conditional on that event,  we use Freedman's inequality \citep{freedman1975} to show that with positive probability, the \textsf{FR} policy accepts  $\Omega(\sqrt{T})$ of class 2 customers, and thus at most $T-\Omega(\sqrt{T})$ of class 1 customers.
So the revenue of \textsf{FR} is at least $\Omega(\sqrt{T})$ less than the hindsight optimum. The complete proof can be found in Appendix~\S\ref{subsec: proof of HO-FR}.




\subsection{Upper bound of the revenue loss of \textsf{FR}}
\label{subsec:res-ub}
In this section, we provide an upper bound of the expected revenue loss of the \textsf{FR} policy. 
\begin{proposition}\label{thm: DLP and resolving ub}
	The gap between the expected revenue of the \textsf{FR} policy defined in Algorithm~\ref{algo:res} and the optimal value of the $\mathrm{DLP}$ is bounded by
	\[
	v^{\mathrm{DLP}}-v^\mathsf{FR}= O(\sqrt{T}).
	\]
	The constant pre-factor depends on the customer arrival rate $\lambda_j$ ($\forall j\in[n]$), the revenues per customer $r_j$ ($\forall j\in[n]$), and the BOM matrix $A$; however, it does not depend on the starting capacity $C_l$ ($\forall l\in[m]$).
\end{proposition}

Since $v^{\mathrm{DLP}}$ is an upper bound of the expected revenue under the optimal policy (see Section~\ref{subsec: DLP}), Proposition~\ref{thm: DLP and resolving ub} immediately implies that the expected revenue loss of the \textsf{FR} policy when compared with the optimal revenue is bounded by $O(\sqrt{T})$. That is,
$
v^* - v^\mathsf{FR} \leq v^{\mathrm{DLP}} - v^\mathsf{FR} = O(\sqrt{T}).
$
Combining Propositions~\ref{thm: ho res lb} and \ref{thm: DLP and resolving ub} gives
$
v^*-v^\mathsf{FR}=\Theta(\sqrt{T}).
$

The proof of Proposition \ref{thm: DLP and resolving ub} can be found in Appendix~\S\ref{subsec: proof DLP-FR}. The proof is based on the following idea.
Since the LP solved under the \textsf{FR} policy and the DLP \eqref{form: DLP2} only differ in the right hand side of the capacity constraints, $b(t)$ and $b$, the expected revenue loss of the \textsf{FR} policy when compared to the optimal value of the DLP can be expressed in terms of $b(t)$ and $b$. More specifically, we show that the expected revenue loss during $[t,t+1)$ can be expressed as $O(\E[(b_l- b_l(t))^+])$ for each resource $l\in[m]$. 
Then, using the relationship between the average remaining capacity, $b(t)$, and the number of accepted customers up to time $t$, we prove that
$
O(\E[(b_l-b_l(t))^+]) =  O(\frac{1}{\sqrt{T-t}}).
$
This completes the proof since $\sum_{t=0}^{T-1} O(\frac{1}{\sqrt{T-t}}) = O(\sqrt{T})$.


Although the $O(\sqrt{T})$ bound in Theorem~\ref{thm: DLP and resolving ub} is looser than the $O(1)$ bound of \textsf{FR} in \citet{JasinReSolvingHeuristicBounded2012}, it does not require the additional condition that the optimal solution to the DLP is nondegenerate. 
Given that the expected revenue loss of \textsf{SPA} is also $O(\sqrt{T})$ (see Appendix \S\ref{app:spa}), we conclude that re-solving at least guarantees the same order of revenue loss compared to no re-solving.


\section{Numerical Experiment}

In this section, we evaluate the numerical performance of five different heuristics, which include
	\begin{enumerate}
		\item \textsf{SPA}:  static probabilistic allocation heuristic (Algorithm~\ref{algo:FP})
		\item \textsf{FR}:  frequent re-solving heuristic (Algorithm~\ref{algo:res})
		\item \textsf{IRT}: infrequent re-solving with threhoslding (Algorithm~\ref{algo:less-is-more})
        \item \textsf{IR}: this algorithm uses the same re-solving schedule as \textsf{IRT} but without applying thresholding; i.e., the acceptance probability at iteration $u$ is always set to $p^u_j \leftarrow x^u_j / \lambda_j$
			\item \textsf{FRT}: this algorithm is motivated by \textsf{IRT}. We apply the same $\tau^{-1/4}$ thresholds from  \textsf{IRT} to  the frequent re-solving algorithm; see the complete description in Algorithm~\ref{algo:TFR}.
	\end{enumerate}
{\color{black}
Recall that \textsf{IRT} has two distinct features compared to \textsf{FR}: it uses an infrequent re-solving schedule and adds thresholds for acceptance probabilities.
The motivation to include \textsf{IR} and \textsf{FRT} in this test is to evaluate which of the two features plays a more important role.
}

{\SingleSpacedXI
\begin{algorithm} 
	\begin{algorithmic}
		\State {\bf initialize}: set $C(0)=C$ and $b(0)=C/T$
		\For {$t = 0, 1, \ldots, T-1$}
		\State set $x(t) \leftarrow \arg\max_x \Big\{\sum_{j=1}^{n}r_jx_j \;\Big|\; \sum_{j=1}^{n}A_jx_j\leq b(t), \text{ and }\, 0\leq x_j \leq \lambda_j, \forall  j \in [n] \Big\}$
		\State set $C' \leftarrow C(t)$
		\For {all customers arriving in $[t,t+1)$}
		\If {the customer belongs to class $j$ and $A_j \leq C'$ ($\forall j\in[n]$)}
		\If {$x_j(t) < \lambda_j (T-t)^{-1/4}$}
		\State reject the customer
		\ElsIf {$x_j(t) >\lambda_j (1 - (T-t)^{-1/4})$}
		\State accept the customer
		\Else
		\State accept the customer with probability $x_j(t) / \lambda_j$ 
		\EndIf
		\State if the customer is accepted, update $C'\leftarrow C' - A_j$
		\Else 
		\State reject the customer
		\EndIf
		\EndFor
		\State set $C(t+1)\leftarrow C'$ and $b(t+1)\leftarrow \frac{C(t+1)}{T-t-1}$
		\EndFor
	\end{algorithmic}
	\caption{\textsf{Frequent Re-solving with Thresholding: \textsf{FRT} }\label{algo:TFR}}
\end{algorithm}
}

	\subsection{Single resource}
	We consider a revenue management problem with a single resource and two classes of customers. We assume that customers from each class arrive according to a Poisson process with rate 1. The arrivals of two classes are independent.  Customers from both classes, if accepted, consume one unit of resource, but pay different prices, $r_1$ and $r_2$.	
We consider two cases: 1) $r_1=2, r_2=1$ and 2) $r_1=5,r_2=1$.
	We also test three settings for the average capacity per unit time: $b=$1, 1.1 and 1.5. When the average capacity per unit time is 1, the solution to the DLP is degenerate. The scenario where the average capacity is 1.1 represents a setting where the DLP solution is ``nearly degenerate,'' and the scenario of 1.5 represents a setting where the DLP solution is far away from any degenerate point.
	We simulate the heuristics for two price and three average capacity per unit time scenarios defined above and for varying horizon length $T=500, 1000,\ldots,5000$.
	
\begin{figure}[!htp]
		\centering
		\subfloat[$b=1, r_1=2$ and $r_2=1$]{%
			\includegraphics[clip,width=0.5\columnwidth]{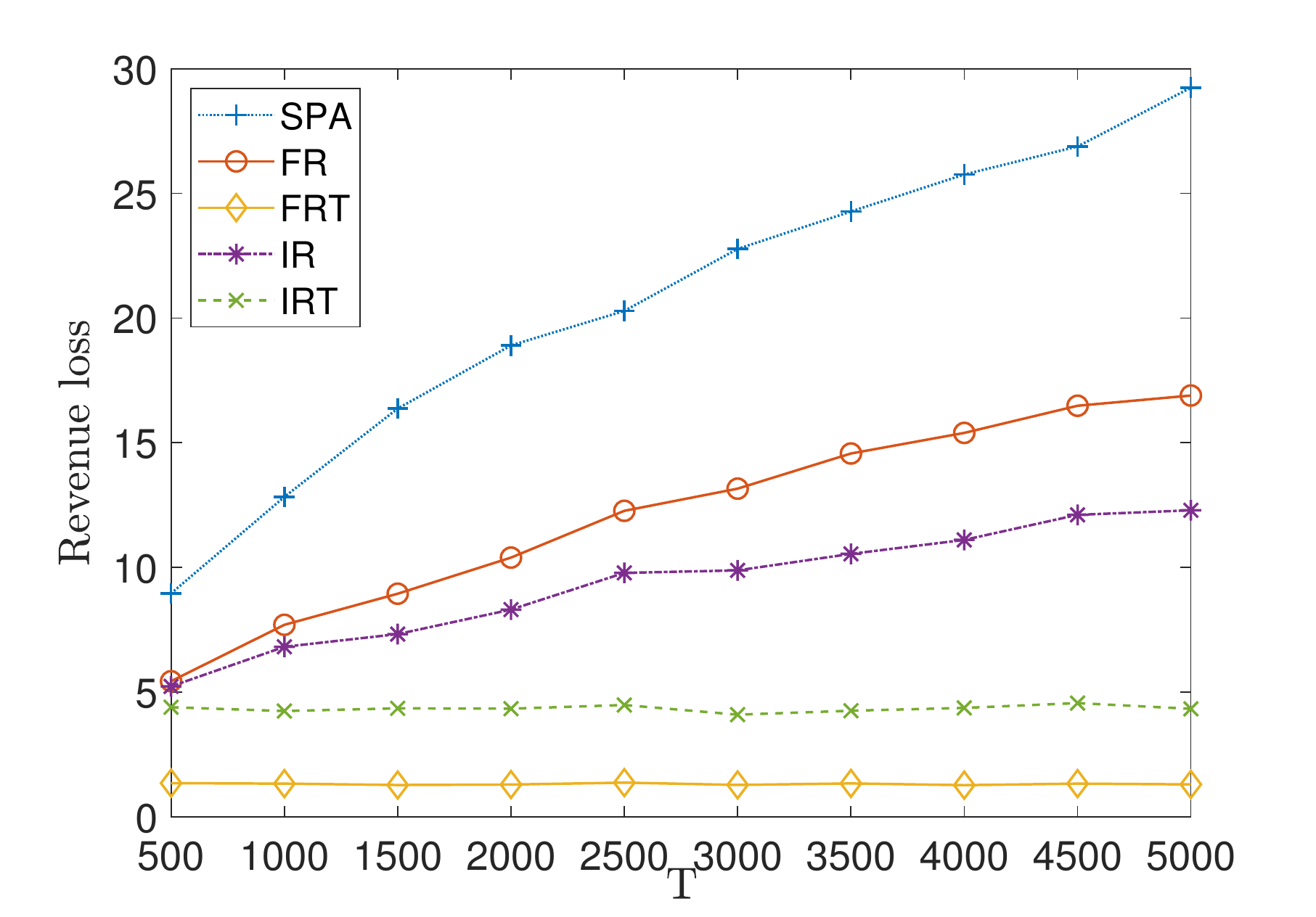}%
		}
	\subfloat[$b=1, r_1=5$ and $r_2=1$]{%
		\includegraphics[clip,width=0.5\columnwidth]{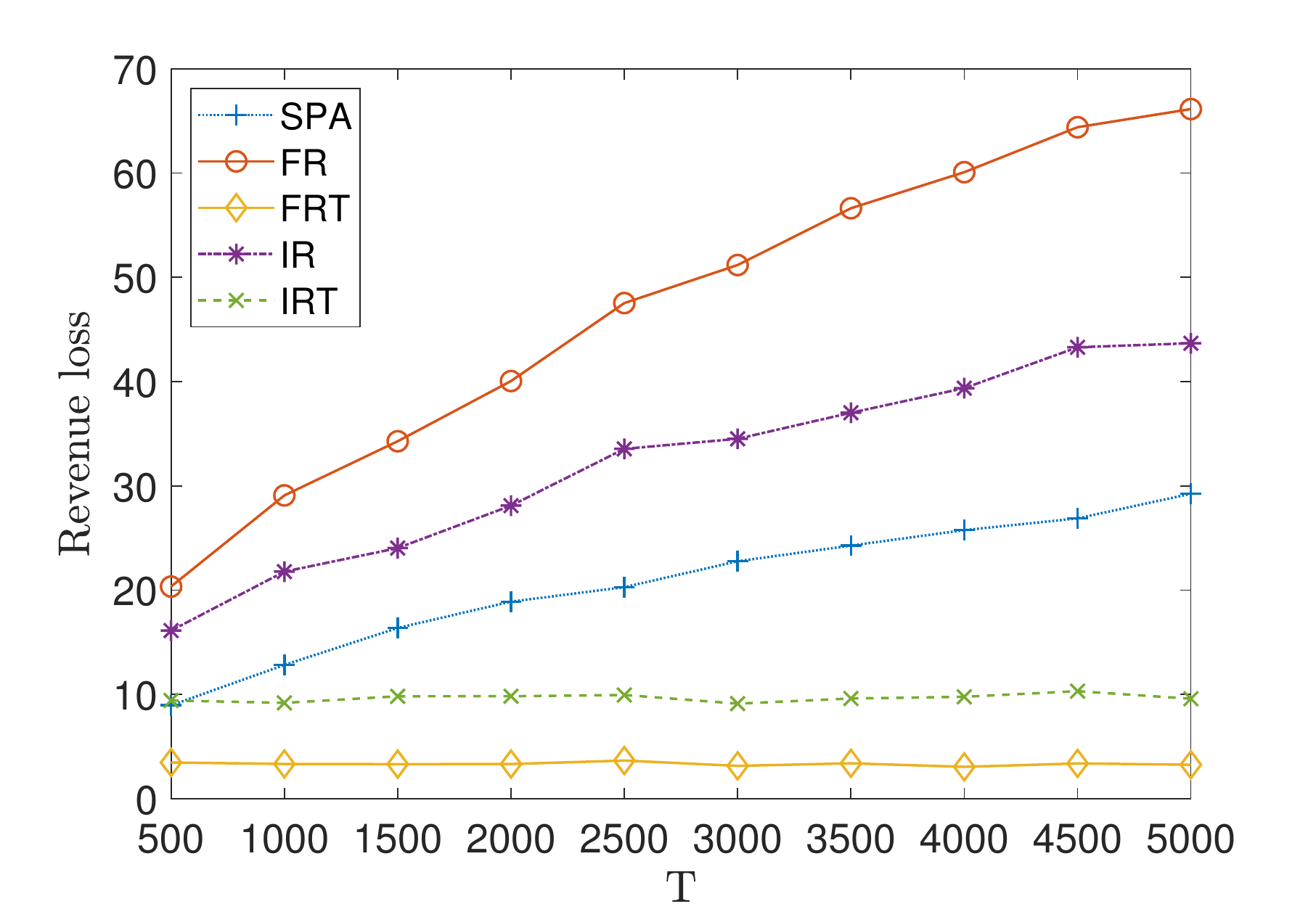}%
	}
		\\
		\subfloat[$b=1.1, r_1=2$ and $r_2=1$]{%
			\includegraphics[clip,width=0.5\columnwidth]{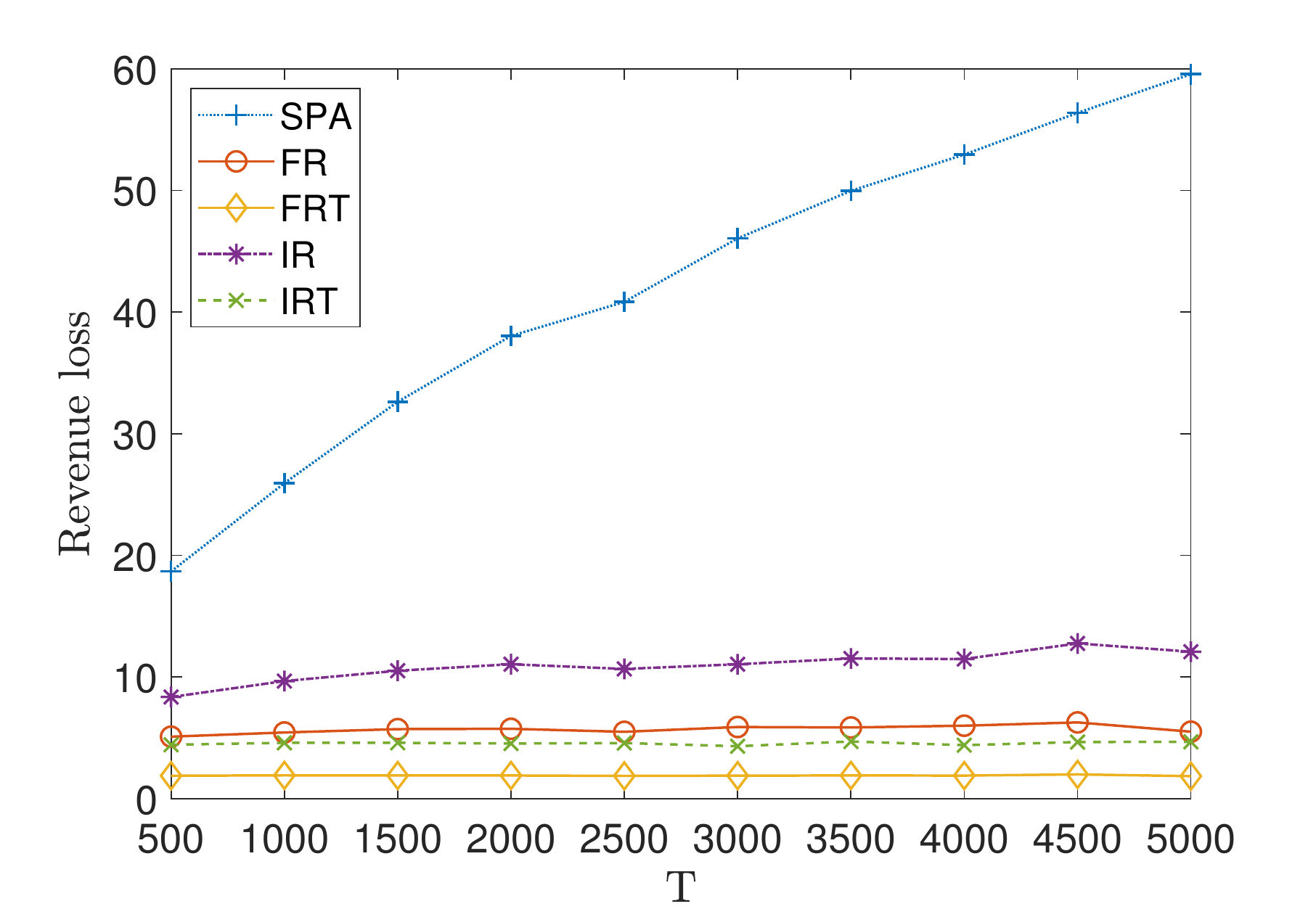}%
		}
	\subfloat[$b=1.1, r_1=5$ and $r_2=1$]{%
		\includegraphics[clip,width=0.5\columnwidth]{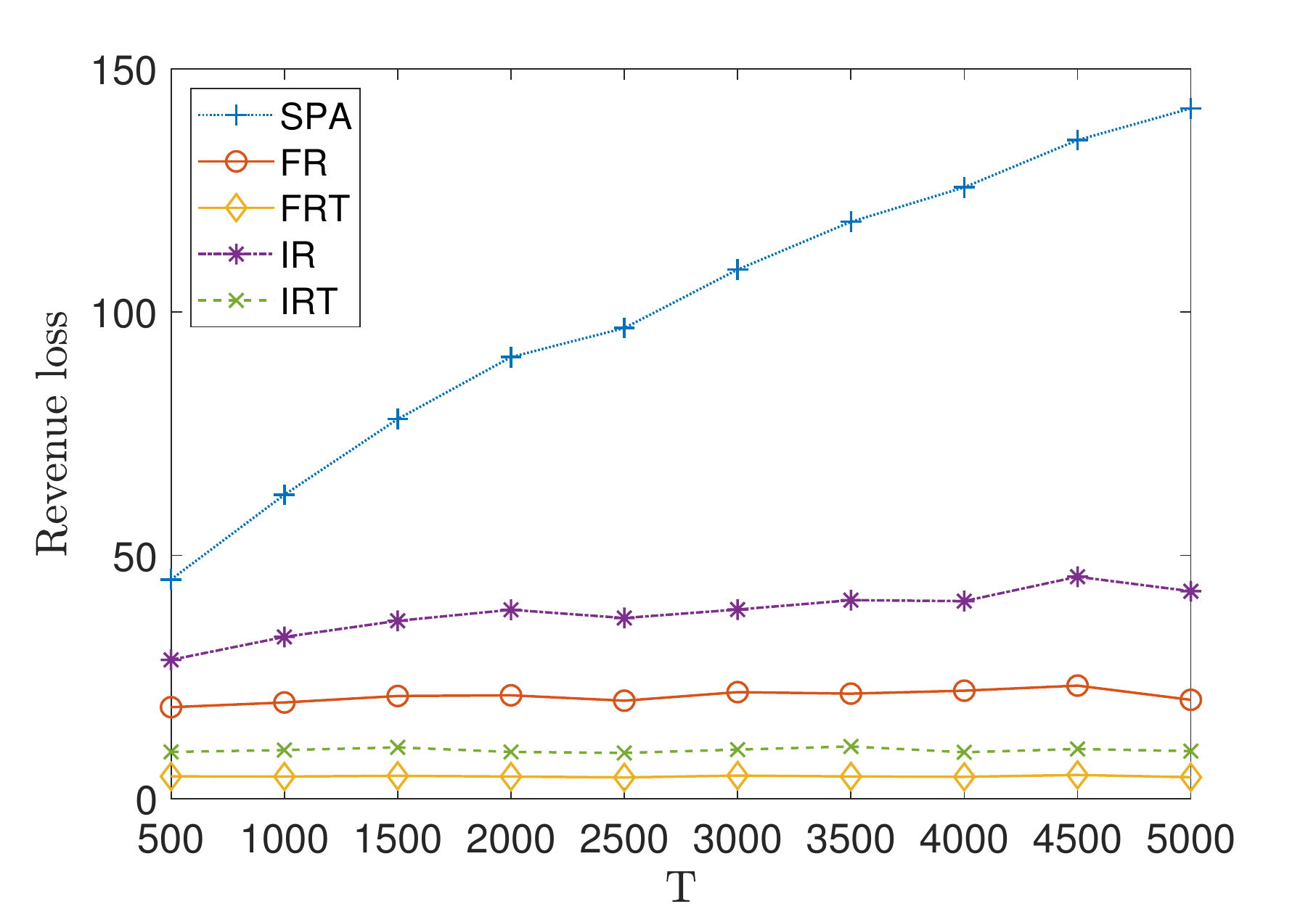}%
	}\\
			\subfloat[$b=1.5, r_1=2$ and $r_2=1$]{%
		\includegraphics[clip,width=0.5\columnwidth]{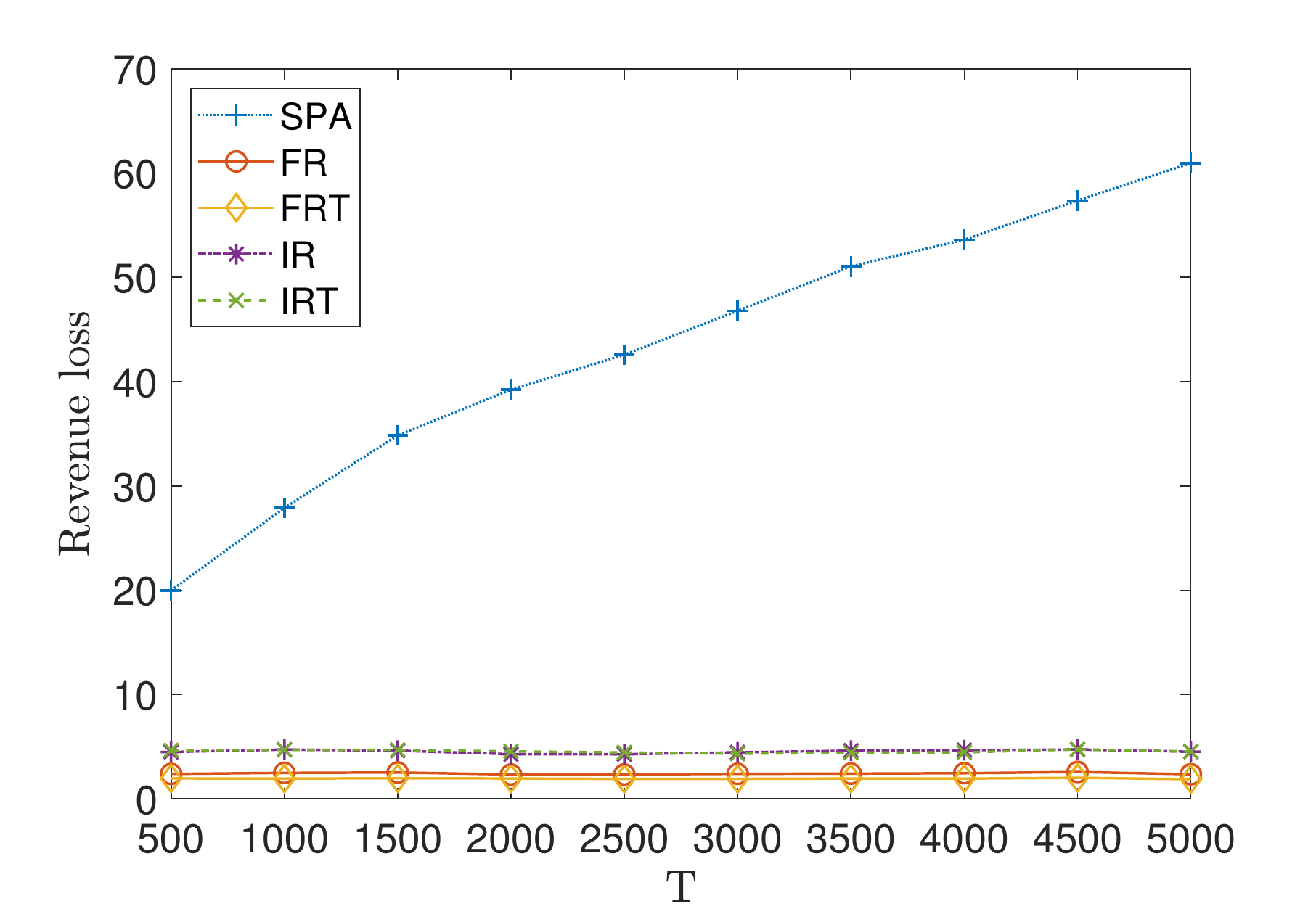}%
	}
\subfloat[$b=1.5, r_1=5$ and $r_2=1$]{%
	\includegraphics[clip,width=0.5\columnwidth]{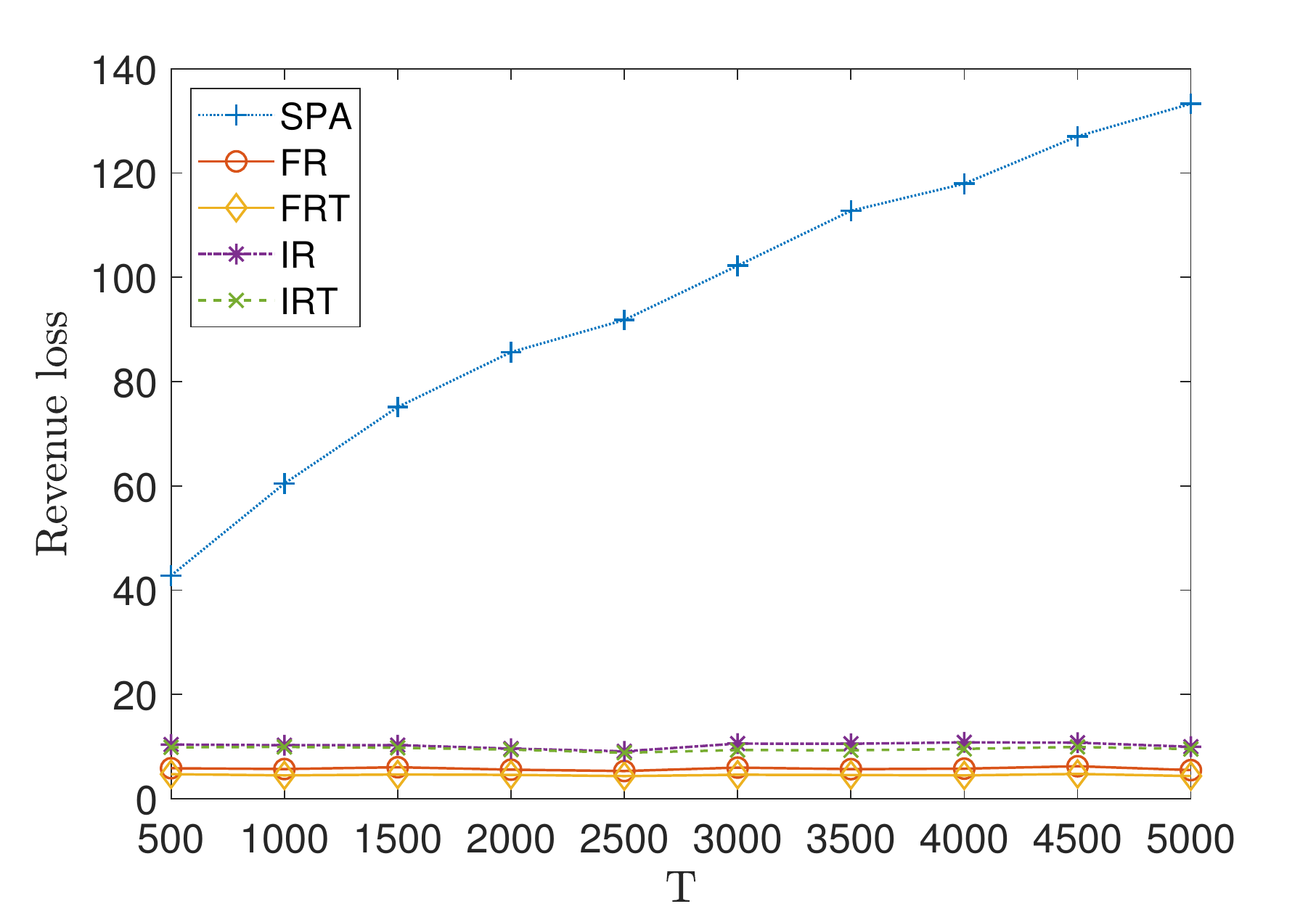}%
}
		\caption{Regret under the \textsf{SPA}, the \textsf{FR}, the \textsf{FRT}, the \textsf{IR}, and the \textsf{IRT} policies for $T=500, 1000,\dots,5000$.}
		\label{fig:numerical results single}
	\end{figure}

	Fig.~\ref{fig:numerical results single} plots the average revenue losses  compared to the hindsight optimum upper bound under  \textsf{SPA},  \textsf{FR},  \textsf{FRT},  \textsf{IR}, and  \textsf{IRT}  over 1000 sample paths. The first column shows the case when $r_1=2$ and $r_2=1$, while the second column shows the case when $r_1=5$ and $r_2=1$. The first, the second and the third rows illustrate the case when $b=1$, $b=1.1$, and $b=1.5$  respectively.
	We make the following observations:
	\begin{enumerate}
		\item When $r_1=2$ and $r_2=1$, the expected revenue loss under  \textsf{SPA} is the largest for all average capacity per unit time and horizon length. This does not hold when $r_1=5$, $r_2=1$ and $b=1$, where \textsf{SPA} is better than either frequent re-solving (\textsf{FR}) or infrequent re-solving (\textsf{IR}).
		\item The expected revenue loss under \textsf{IR} is higher than the expected revenue loss under  \textsf{FR} except when the problem is degenerate ($b=1$). We conclude that choosing an infrequent re-solving schedule alone is not enough to achieve $O(1)$ loss.
		\item The expected revenue losses under  \textsf{IRT} and  \textsf{FRT} remain constant for all cases as the horizon length increases. Moreover, although we don't have theoretical guarantee for \textsf{FRT},  the expected revenue loss under  \textsf{FRT} often appears smaller than the expected revenue loss under  \textsf{IRT}. {\color{black} This implies that appropriate thresholding is the main factor that leads to  uniformly bounded regret for re-solving heuristics.}
	\end{enumerate}

	\subsection{Multiple resources}
	Next, we consider a network revenue management problem with multiple resources. We consider the problem when there are five classes of customers and four types of resources. We assume that customers from each class arrive according to a Poisson process with rate 1; the arrivals of different classes are independent. The vector of the average capacities per unit time is given by
$
	b=[1,1,1,1]^\top.
$
	The vector of the revenue earned by accepting customers is given by
$
	r=[10, 3, 6, 1, 2]^\top.
$
	The bill-of-materials matrix is given by
	\begin{align*}
		A=
		\begin{bmatrix}
		1 &0 &1 &0 &0\\
		0 &1 &0 &1 &1\\
		1 &1 &0 &0 &0\\
		0 &0 &0 &0 &1
		\end{bmatrix}.
	\end{align*}
	We simulate the heuristics for varying horizon length $T=500, 1000,\ldots,5000$. Notice that in this example, the optimal solution to the DLP is degenerate.
	\begin{figure}[!htp]
		\centering
		\includegraphics[clip,width=0.7\columnwidth]{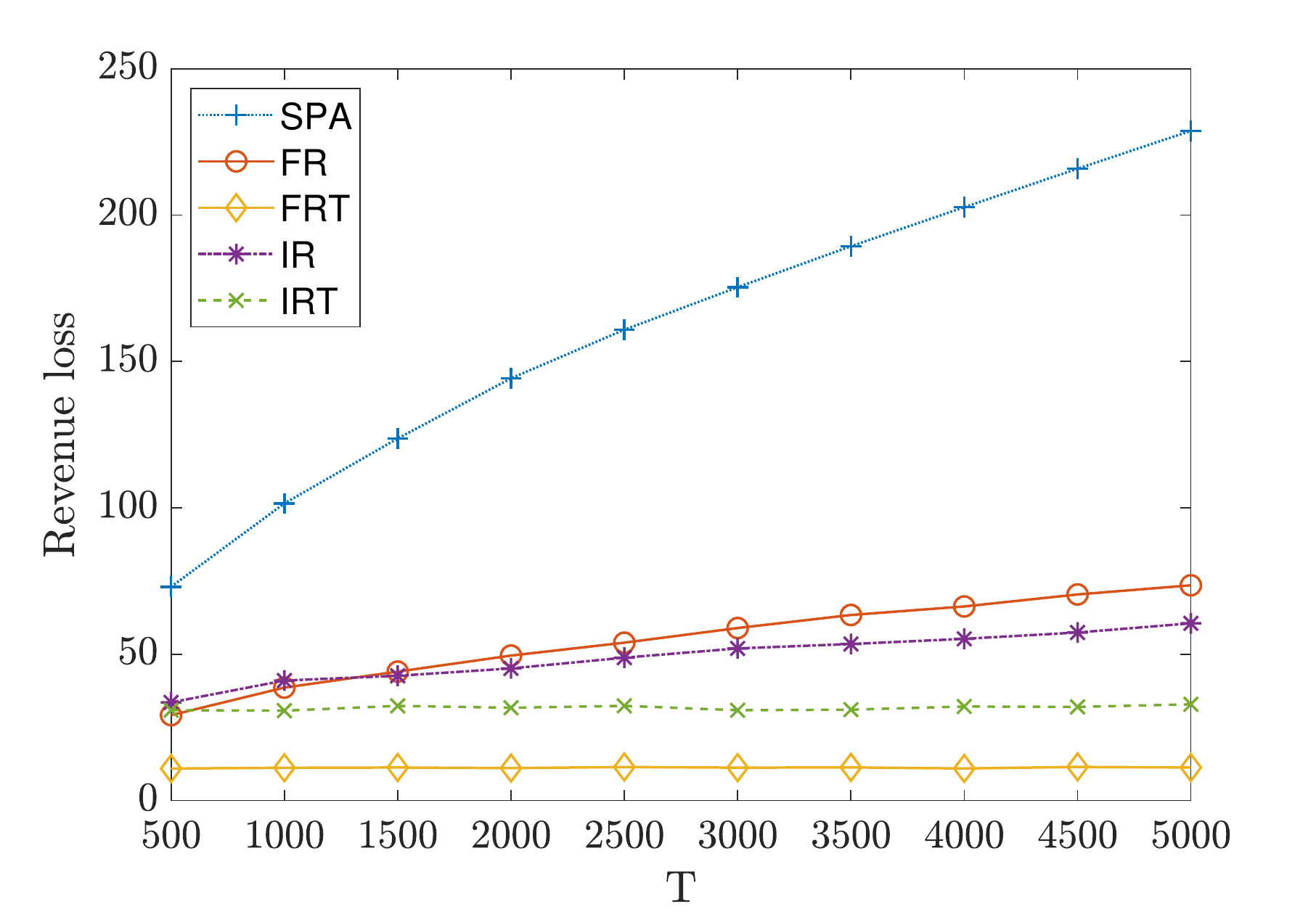}
		\caption{The expected revenue losses (regret) under the \textsf{SPA}, the \textsf{FR}, the \textsf{FRT}, the \textsf{IR}, and the \textsf{IRT} when compared with the hindsight optimal for $T=500, 1000,\dots,5000$.}
		\label{fig:numerical results multi}
	\end{figure}

Fig.~\ref{fig:numerical results multi} plots the average revenue losses under  \textsf{SPA},  \textsf{FR},  \textsf{FRT},  \textsf{IR}, and  \textsf{IRT} over 1000 sample paths. The result shows that the revenue losses of  \textsf{SPA} scales poorly with horizon length $T$. In comparison, the revenue losses of \textsf{FR} and  \textsf{IR} increase more slowly when $T$ increases, and \textsf{IR} seems to perform slightly better for large $T$. 
The revenue losses of  \textsf{FRT} and  \textsf{IRT} remain constant as $T$ grows. Moreover, the expected revenue loss under the \textsf{IRT} is higher than the expected revenue loss under the \textsf{FRT}. {\color{black} Again, this result implies that among the two factors, infrequent re-solving and thresholding, the latter plays a more important role.}

\section{Conclusion and Discussion}

{\color{black}
We study re-solving heuristics for the network revenue management (NRM) problem. A re-solving heuristic periodically re-optimizes a deterministic LP approximation of the original NRM problem. The main question considered in this paper is:  can we find a simple and computationally efficient re-solving heuristic, whose expected revenue loss compared to the optimal policy is bounded by a constant even when both the time horizon and the resource capacities scale up?

We answer the above question in the affirmative by proposing a re-solving heuristic called \textsf{Infrequent Re-solving with Thresholding (IRT)}, whose revenue loss is bounded by a constant independent of time horizon and resource capacities. This finding improves a previous result by \cite{JasinReSolvingHeuristicBounded2012}, showing that \textsf{Frequent Re-solving (FR)}, 
an algorithm that re-solves the DLP after each unit of time, has $O(1)$ revenue loss, but requires the optimal solution to the DLP to be nondegenerate.
Moreover, we show that when both time horizon and resource capacities scale up by $k=1,2,\ldots$,
\textsf{Frequent Re-solving (FR)}  has a revenue loss of $\Theta(\sqrt{k})$. This is a negative result, as most DLP-based heuristics can achieve the same revenue loss rate without using re-solving at all. 

Our simulation results show that when the controls from  \textsf{FR}  are adjusted by some thresholds, the resulting algorithm \textsf{FRT} has very promising numerical performance and seems to have a bounded revenue loss as well.  So far, we are not able to prove this result, mainly because the induction based proof we developed for \textsf{IRT} breaks down when the DLP is re-optimized every period. Recently, \cite{vera2018bayesian} propose a different re-solving heuristic for the NRM problem, where the DLP is re-optimized every period and a fixed acceptance probability threshold of 0.5 is applied to every class for all periods. They show their heuristic also achieves $O(1)$ regret. Although the fixed threshold used by \cite{vera2018bayesian} is different from the time-varying thresholds we proposed in the \textsf{FRT} algorithm, we think their analysis technique may be helpful to establish the revenue loss bound of \textsf{FRT}.


}


\bibliographystyle{ormsv080} 
\bibliography{propbib}

\begin{thebibliography}{28}
\expandafter\ifx\csname natexlab\endcsname\relax\def\natexlab#1{#1}\fi
\expandafter\ifx\csname url\endcsname\relax
  \def\url#1{{\tt #1}}\fi
\expandafter\ifx\csname urlprefix\endcsname\relax\def\urlprefix{URL }\fi
\expandafter\ifx\csname urlstyle\endcsname\relax
  \expandafter\ifx\csname doi\endcsname\relax
  \def\doi#1{doi:\discretionary{}{}{}#1}\fi \else
  \expandafter\ifx\csname doi\endcsname\relax
  \def\doi{doi:\discretionary{}{}{}\begingroup \urlstyle{rm}\Url}\fi \fi

\bibitem[{Arlotto and Gurvich(2017)}]{arlotto2017uniformly}
Arlotto, Alessandro, Itai Gurvich. 2017.
\newblock Uniformly bounded regret in the multi-secretary problem.
\newblock ArXiv preprint arXiv:1710.07719.

\bibitem[{Arlotto and Xie(2018)}]{arlotto2018logarithmic}
Arlotto, Alessandro, Xinchang Xie. 2018.
\newblock Logarithmic regret in the dynamic and stochastic knapsack problem.
\newblock {\it arXiv preprint arXiv:1809.02016\/} .

\bibitem[{Babaioff et~al.(2007)Babaioff, Immorlica, Kempe, and
  Kleinberg}]{babaioff2007knapsack}
Babaioff, Moshe, Nicole Immorlica, David Kempe, Robert Kleinberg. 2007.
\newblock A knapsack secretary problem with applications.
\newblock {\it Approximation, randomization, and combinatorial optimization.
  Algorithms and techniques\/}. Springer, Princeton, NJ, 16--28.

\bibitem[{Bertsekas(2005)}]{bertsekas2005dynamic}
Bertsekas, Dimitri~P. 2005.
\newblock {\it Dynamic programming and optimal control (3rd edition), volume
  I\/}.
\newblock Belmont, MA: Athena Scientific.

\bibitem[{Besbes and Zeevi(2012)}]{besbes2012blind}
Besbes, Omar, Assaf Zeevi. 2012.
\newblock Blind network revenue management.
\newblock {\it Operations Research\/} {\bf 60}(6) 1537--1550.

\bibitem[{Chen and Homem-de
  Mello(2010)}]{ChenResolvingstochasticprogramming2010}
Chen, Lijian, Tito Homem-de Mello. 2010.
\newblock Re-solving stochastic programming models for airline revenue
  management.
\newblock {\it Annals of Operations Research\/} {\bf 177}(1) 91--114.

\bibitem[{Cooper(2002)}]{CooperAsymptoticBehaviorAllocation2002}
Cooper, William~L. 2002.
\newblock Asymptotic behavior of an allocation policy for revenue management.
\newblock {\it Operations Research\/} {\bf 50}(4) 720--727.

\bibitem[{Ferreira et~al.(2017)Ferreira, Simchi-Levi, and
  Wang}]{ferreira2017online}
Ferreira, Kris, David Simchi-Levi, He~Wang. 2017.
\newblock Online network revenue management using {Thompson} sampling.
\newblock {\it Operations Research\/} (forthcoming).

\bibitem[{Freedman(1975)}]{freedman1975}
Freedman, David~A. 1975.
\newblock On tail probabilities for martingales.
\newblock {\it Annals of Probabability\/} {\bf 3}(1) 100--118.

\bibitem[{Gallego and van Ryzin(1994)}]{GallegoOptimalDynamicPricing1994}
Gallego, Guillermo, Garrett van Ryzin. 1994.
\newblock Optimal {{Dynamic Pricing}} of {{Inventories}} with {{Stochastic
  Demand}} over {{Finite Horizons}}.
\newblock {\it Management Science\/} {\bf 40}(8) 999.

\bibitem[{Gallego and van Ryzin(1997)}]{gallego1997multiproduct}
Gallego, Guillermo, Garrett van Ryzin. 1997.
\newblock A multiproduct dynamic pricing problem and its applications to
  network yield management.
\newblock {\it Operations Research\/} {\bf 45}(1) 24--41.

\bibitem[{Jasin(2015)}]{jasin2015performance}
Jasin, Stefanus. 2015.
\newblock Performance of an {LP}-based control for revenue management with
  unknown demand parameters.
\newblock {\it Operations Research\/} {\bf 63}(4) 909--915.

\bibitem[{Jasin and Kumar(2012)}]{JasinReSolvingHeuristicBounded2012}
Jasin, Stefanus, Sunil Kumar. 2012.
\newblock A re-solving heuristic with bounded revenue loss for network revenue
  management with customer choice.
\newblock {\it Mathematics of Operations Research\/} {\bf 37}(2) 313--345.

\bibitem[{Jasin and Kumar(2013)}]{JasinAnalysisDeterministicLPBased2013}
Jasin, Stefanus, Sunil Kumar. 2013.
\newblock Analysis of deterministic {LP}-based booking limit and bid price
  controls for revenue management.
\newblock {\it Operations Research\/} {\bf 61}(6) 1312--1320.

\bibitem[{Kleinberg(2005)}]{kleinberg2005multiple}
Kleinberg, Robert. 2005.
\newblock A multiple-choice secretary algorithm with applications to online
  auctions.
\newblock {\it Proceedings of the sixteenth annual ACM-SIAM symposium on
  Discrete algorithms\/}. Society for Industrial and Applied Mathematics, PA,
  USA, 630--631.

\bibitem[{Kleywegt and Papastavrou(1998)}]{kleywegt1998dynamic}
Kleywegt, Anton~J, Jason~D Papastavrou. 1998.
\newblock The dynamic and stochastic knapsack problem.
\newblock {\it Operations Research\/} {\bf 46}(1) 17--35.

\bibitem[{Liu and {van Ryzin}(2008)}]{liu_choice-based_2008}
Liu, Qian, Garrett {van Ryzin}. 2008.
\newblock On the choice-based linear programming model for network revenue
  management.
\newblock {\it Manufacturing \& Service Operations Management\/} {\bf 10}(2)
  288--310.

\bibitem[{Maglaras and Meissner(2006)}]{maglaras2006dynamic}
Maglaras, Constantinos, Joern Meissner. 2006.
\newblock Dynamic pricing strategies for multiproduct revenue management
  problems.
\newblock {\it Manufacturing \& Service Operations Management\/} {\bf 8}(2)
  136--148.

\bibitem[{Pollard(2015)}]{pollard2015inequality}
Pollard, David. 2015.
\newblock A few good inequalities.
\newblock Http://www.stat.yale.edu/~pollard/Books/Mini/Basic.pdf.

\bibitem[{Reiman and Wang(2008)}]{ReimanAsymptoticallyOptimalPolicy2008}
Reiman, Martin~I., Qiong Wang. 2008.
\newblock An asymptotically optimal policy for a quantity-based network revenue
  management problem.
\newblock {\it Mathematics of Operations Research\/} {\bf 33}(2) 257--282.

\bibitem[{Secomandi(2008)}]{secomandi2008analysis}
Secomandi, Nicola. 2008.
\newblock An analysis of the control-algorithm re-solving issue in inventory
  and revenue management.
\newblock {\it Manufacturing \& Service Operations Management\/} {\bf 10}(3)
  468--483.

\bibitem[{Shevtsova(2011)}]{shevtsova2011absolute}
Shevtsova, Irina. 2011.
\newblock On the absolute constants in the berry-esseen type inequalities for
  identically distributed summands.
\newblock {\it arXiv preprint arXiv:1111.6554\/} .

\bibitem[{Talluri and van Ryzin(1998)}]{talluri_analysis_1998}
Talluri, Kalyan, Garrett van Ryzin. 1998.
\newblock An analysis of bid-price controls for network revenue management.
\newblock {\it Management Science\/} {\bf 44}(11) 1577--1593.

\bibitem[{Talluri and van Ryzin(2004{\natexlab{a}})}]{talluri2004revenue}
Talluri, Kalyan, Garrett van Ryzin. 2004{\natexlab{a}}.
\newblock Revenue management under a general discrete choice model of consumer
  behavior.
\newblock {\it Management Science\/} {\bf 50}(1) 15--33.

\bibitem[{Talluri and van Ryzin(2004{\natexlab{b}})}]{talluri_theory_2004}
Talluri, Kalyan~T., Garrett van Ryzin. 2004{\natexlab{b}}.
\newblock {\it The Theory and Practice of Revenue Management\/}.
\newblock No.~68 in International series in operations research \& management
  science, {Kluwer Academic Publishers}, Boston, Mass.

\bibitem[{Vera and Banerjee(2018)}]{vera2018bayesian}
Vera, Alberto, Siddhartha Banerjee. 2018.
\newblock The bayesian prophet: A low-regret framework for online decision
  making.
\newblock {\it Available at SSRN: https://ssrn.com/abstract=3158062\/} .

\bibitem[{Williamson(1992)}]{williamson1992airline}
Williamson, Elizabeth~Louise. 1992.
\newblock Airline network seat inventory control: Methodologies and revenue
  impacts.
\newblock Ph.D. thesis, Massachusetts Institute of Technology.

\bibitem[{Wu et~al.(2015)Wu, Srikant, Liu, and Jiang}]{wu2015algorithms}
Wu, Huasen, R~Srikant, Xin Liu, Chong Jiang. 2015.
\newblock Algorithms with logarithmic or sublinear regret for constrained
  contextual bandits.
\newblock {\it Advances in Neural Information Processing Systems\/}. 433--441.

\end{thebibliography}
\newpage
\clearpage
\begin{APPENDIX}{\ }
\section{Additional Results}
\label{append:additional-result}

\subsection{A note on the DLP upper bound}

In the paper, we establish the revenue loss of heuristics by comparing their revenues to the hindsight optimum upper bound $v^{\mathrm{HO}}$. This bound is tighter than the DLP upper bound $v^{\mathrm{DLP}}$. The following result suggests that $v^{\mathrm{DLP}}$ is not an appropriate benchmark to prove $O(1)$ revenue loss, because even the gap between the optimal policy $v^*$ and $v^{\mathrm{DLP}}$ is $\Omega(\sqrt{T})$.

\begin{proposition}\label{prop: DLP and optimal lb}
	The gap between the optimal value of the DLP and the optimal value obtained by dynamic programming is bounded below by
	\[
	v^{\mathrm{DLP}}-v^*= \Omega(\sqrt{T}).
	\]
\end{proposition}
\proof{Proof of Proposition \ref{prop: DLP and optimal lb}.}
	To prove Proposition \ref{prop: DLP and optimal lb}, we consider the following instance. In this instance, there is only one class of customer and one type of resource. So for simplicity, we will suppress the subscriptions. Suppose the expected number of arrivals in one period is Poisson process with rate $\lambda$. The revenue earned by accepting a customer is 1. The resource has the capacity $\lambda T$ and the amount of the resource used to serve one customer is 1. Therefore, the DLP formulation is given by
	\[
	\max_x\Big\{T x \;\Big|\; x\leq \lambda T/T =\lambda,\; 0\leq x \leq \lambda\Big\}.
	\]	 
	It easily verified that, we have $x^*=\lambda$ and thus $v^{\mathrm{DLP}}=\lambda T$.

	On the other hand, it is obvious that the optimal policy is to admit all customers in $[0,T]$ subject to the capacity constraint. Specifically,  the optimal number of the admitted customers is either the number of the arriving customers in $[0,T$] or the capacity level, whichever is lower. Therefore, 
	the optimal revenue of the above problem instance is given by
	\begin{align}
	v^*&=\E[\min(\Lambda(T),\lambda T)]=\lambda T -\E[\max(\lambda T-\Lambda(T),0)]\nonumber=\lambda T-\sqrt{\lambda T}\E\left[\max\left(\frac{\lambda T-\Lambda(T)}{\sqrt{\lambda T}},0\right)\right].\nonumber
	\end{align}
	From the Markov's inequality, it follows that
	\begin{align}
	\E\left[\max\left(\frac{\lambda T-\Lambda(T)}{\sqrt{\lambda T}},0\right)\right]
	&\geq \Prob\left(\max\left(\frac{\lambda T-\Lambda(T)}{\sqrt{\lambda T}},0\right) \geq 1\right).\nonumber
	\end{align}
	Consequently, we can write
	\begin{align}
	v^*
	&\leq \lambda T-\sqrt{\lambda T}\Prob\left(\max\left(\frac{\lambda T-\Lambda(T)}{\sqrt{\lambda T}},0\right) \geq 1\right)\nonumber\\
	&=\lambda T-\sqrt{\lambda T}\Prob\left(\frac{\lambda T-\Lambda(T)}{\sqrt{\lambda T}} \geq 1\right)\nonumber\\
	&=\lambda T-\sqrt{\lambda T}(1-F_T(1)),\label{eq:berry1}
	\end{align}
	where $F_T$ is the cumulative distribution function of $\frac{\lambda T-\Lambda(T)}{\sqrt{\lambda T}}$. Let $\Phi$ denote the cumulative distribution function of standard normal distribution. Since
	\begin{align}
	1-F_T(1)
	&= 1-\Phi(1)+\Phi(1)-F_T(1)  \geq 1-\Phi(1)-|F_T(1)-\Phi(1)|.\label{eq:berry2}
	\end{align}
	{\color{black} We use the Berry-Esseen theorem (Lemma~\ref{lem:berry-esseen} in Appendix~\ref{app:lemmas}) to bound $|F_T(1)-\Phi(1)|$.
		Let $X_i=\lambda-(\Lambda(i)-\Lambda(i-1))$ for $i=1,\ldots,T$. From the stationary and independent increment properties of Poisson processes, we observe that $X_i$ are i.i.d. with $\E[X_1]=0,\; \E[X^2]=\lambda$ and $\E[|X_1^3|]=\lambda$. Since
\[
\frac{X_1+\ldots+X_T}{\sqrt{\lambda T}} =\frac{\lambda T - \Lambda(T)}{\sqrt{\lambda T}},
\]
from Lemma \ref{lem:berry-esseen} we have
	\begin{align}\label{eq:berry3}
|F_T(1)-\Phi(1)|
& \leq \frac{0.4748 \lambda}{ \sqrt{\lambda^3 T}}=\frac{0.4748}{\sqrt{\lambda T}}.
\end{align}
}
	Combining \eqref{eq:berry1}, \eqref{eq:berry2} and \eqref{eq:berry3}, we can write
	\begin{align}
	v^*
	&\leq \lambda T-\sqrt{\lambda T}\left(1-\Phi(1) - \frac{0.4748}{\sqrt{\lambda T}}\right)
	= \lambda T-\sqrt{\lambda T}(1-\Phi(1))+0.4748.\nonumber
	\end{align}
	Since $v^{\mathrm{DLP}}=\lambda T$ and $1-\Phi(1)=0.1587$, it follows that
	\begin{align}
	v^{\mathrm{DLP}}-v^*
	&\geq 0.1587\sqrt{\lambda T}-0.4748 =\Omega(\sqrt{T}).\nonumber \Halmos	
	\end{align}
\endproof

\subsection{Revenue loss of static probabilistic allocation}
\label{app:spa}

The following result is a well-known in the revenue management literature by
\citet{GallegoOptimalDynamicPricing1994, gallego1997multiproduct}; see also \citet{CooperAsymptoticBehaviorAllocation2002,ReimanAsymptoticallyOptimalPolicy2008}. 

\begin{proposition}\label{prop: DLP and FPA}
	The gap between the optimal value of the DLP and the optimal value obtained by the static probabilistic allocation (\textsf{SPA}) heuristic is bounded above by
	\[
	v^{\mathrm{DLP}}-v^\mathsf{SPA}= O(\sqrt{T}).
	\]
	The constant pre-factor depends on the customer arrival rate $\lambda_j$ ($\forall j\in[n]$), the revenues per customer $r_j$ ($\forall j\in[n]$), and the BOM matrix $A$; however, it does not depend on the starting capacity $C_l$ ($\forall l\in[m]$).
\end{proposition}

Since the revenue of the optimal policy, $v^*$, and the hindsight optimum, $v^{\mathrm{HO}}$, satisfies $v^*\leq v^{\mathrm{HO}} \leq v^{\mathrm{DLP}}$, a corollary of the result is $v^* - v^{\textsf{SPA}}= O(\sqrt{T})$ and $v^{\mathrm{HO}} - v^{\textsf{SPA}}= O(\sqrt{T})$.

In the revenue management literature, this result is often proved under the additional assumption that resource capacities and customer arrivals are both scaled up at the same rate. However, the result in fact holds for arbitrary capacity levels. 
We need this fact in the proof of Theorem~\ref{thm: TR}. To make the proof of Theorem~\ref{thm: TR} self-contained, we include a proof of the proposition below.

\proof{Proof of Proposition \ref{prop: DLP and FPA}.}
By Eq~\ref{form: DLP2}, we have $v^{\mathrm{DLP}}=T\sum_{j\in[n]}r_jx_j^*$.
We bound $v^{\textsf{SPA}}$ in two steps. First, consider a hypothetical setting where remaining capacities are allowed to become negative. Since \textsf{SPA} accepts each class $j$ customer with probability $x^*_j/\lambda_j$ if capacity constraints are ignored, {\color{black}the expected revenue of  \textsf{SPA} under this hypothetical  setting} is
\[
\E\left[\sum_{j\in[n]}\int_{0}^T r_j\lambda_j \frac{x^*_j}{\lambda_j}dt\right]=\sum_{j\in[n]}r_jx_j^*T=v^\mathrm{DLP}.
\]
In reality, remaining capacity is always nonnegative, and customers must be rejected if there is insufficient capacity. So to correct the revenue calculation of the hypothetical setting, we must subtract the revenue associated with customers who are rejected due to insufficient capacity. This part of revenue is bounded by
\begin{equation}\label{eq:lost-sales2}
\sum_{l=1}^m r_{max}^l \E\bigl[ \bigl(\sum_{j=1}^n a_{lj} X_j - C_l\bigr)^+ \bigr],
\end{equation}
	where $X_j$  is the number of class $j$ customers that would have been accepted by \textsf{SPA} without capacity limits, which follow a Poisson distribution with mean $x^*_j T$, and  $r_{max}^l$ is the largest possible revenue gain by increasing the capacity of resource $l$ by one unit, i.e., $r^l_{max}=\max_{j\in[n]}\{
	r_j \mathrm{I}(a_{lj}>0)/a_{lj}\}$.
	{\color{black}We have 
\begin{align*}
\E\bigl[ \bigl(\sum_{j=1}^n a_{lj} X_j - C_l\bigr)^+ \bigr] 
&\leq \E\bigl[ \bigl(\sum_{j=1}^n a_{lj} X_j - \sum_{j=1}^n a_{lj}x^*_jT\bigr)^+ \bigr]\\
&\leq  \E\bigl[ \bigl|\sum_{j=1}^n a_{lj} X_j -  \sum_{j=1}^n a_{lj}x^*_jT\bigr| \bigr]\\
&\leq \sqrt{\E\bigl[\bigl( \sum_{j=1}^n a_{lj} X_j -\sum_{j=1}^n a_{lj}x^*_jT\bigr)^2 \bigr]},
\end{align*}
where the first inequality follows the capacity constraints in DLP \eqref{form: DLP2}, and the last inequality follows from Cauchy-Schwarz inequality.
}
Note that $X_j$'s mean and variance are both equal to $x^*_j T$, and $X_j$'s are independent for all $j\in[n]$. So
\[
\E\bigl[\bigl( \sum_{j=1}^n a_{lj} X_j -\sum_{j=1}^n a_{lj}x^*_jT\bigr)^2 \bigr] = \mathrm{Var}\bigl(\sum_{j=1}^n a_{lj} X_j\bigr)= \sum_{j=1}^n a^2_{lj} \mathrm{Var}\bigl( X_j\bigr)= \sum_{j=1}^n a^2_{lj} x^*_j T
\leq \sum_{j=1}^n a^2_{lj} \lambda_j T,
\]
where the last inequality follows the demand constraints $x^*_j\leq \lambda_j$ in DLP \eqref{form: DLP2}.  Substituting this result to Eq~\eqref{eq:lost-sales2}, we have
\[
v^{\mathrm{DLP}} - v^{\textsf{SPA}} \leq \sum_{l=1}^m r_{max}^l \E\bigl[ \bigl(\sum_{j=1}^n a_{lj} X_j - C_l\bigr)^+ \bigr] \leq \sum_{l=1}^m  r_{max}^l \sqrt{\sum_{j=1}^na^2_{lj} \lambda_j T}. \quad \Halmos
\]
\endproof

\section{Proofs for Section~\ref{sec:new-heuristic}}

In this section, we provide complete proofs for the results on the \textsf{IRT} policy in Section~\ref{sec:new-heuristic}.

\subsection{Proof of Theorem~\ref{thm: TR}}
\label{subsec:proof-TR}

\proof{Proof of Theorem~\ref{thm: TR}.}
	Given remaining capacity $C(t_1)$ at time $t_1\in[0,T]$, let $x(t_1)$ be an optimal solution to the following LP
	\[
	\max_x \Big\{\sum_{j=1}^{n}r_jx_j \;\Big|\; \sum_{j=1}^{n}A_jx_j\leq C(t_1)/(T-t_1), \text{ and }\, 0\leq x_j \leq \lambda_j, \forall  j \in [n] \Big\}.
	\]
	For any $t_2\in (t_1,T]$,  let $V^{\mathsf{SPA}}(t_1,t_2)$ denote the revenue earned in $[t_1,t_2)$ under a static probabilistic allocation policy, where class $j$ customers are accepted with probability $x_j(t_1)/\lambda_j$. Let $V^{\mathsf{SPA}'}(t_1,t_2)$ be the revenue earned in $[t_1,t_2)$, where a class $j$ customer is accepted with the following probability:
	\begin{itemize}
		\item 0, if $x_j(t_1) < \lambda_j (T-t_1)^{-1/4}$
		\item 1, if $x_j(t_1) >\lambda_j (1 - (T-t_1)^{-1/4})$
		\item $x_j(t_1)/\lambda_j$, otherwise.
	\end{itemize}
	
	Let $V^{\mathrm{HO}}(t_1,T)$ denote the revenue earned from solving the hindsight optimum in $[t_1,T]$. That is,  $V^{\mathrm{HO}}(t_1,T)$ is the optimal revenue given the remaining capacity at $t_1$ and a sample path of demand in $(t_1,T]$, given by
	\[
	V^{\mathrm{HO}}(t_1,T) = \max_y \Big\{\sum_{j=1}^{n}r_jy_j \;\Big|\; \sum_{j=1}^{n}A_j y_j\leq C(t_1), \text{ and }\, 0\leq y_j \leq \Lambda_j(T)-\Lambda_j(t_1), \forall  j \in [n] \Big\}.
	\]
	
	Consider policy $\mathsf{IRT}^2$, which re-solves at $t^*_1 = T-T^{5/6}$, and re-solves again at $t_2^*=T-T^{(5/6)^2}=T-T^{25/36}$. 
	Let $v^{\mathsf{IRT}^2}$ be the expected revenue of $\mathsf{IRT}^2$.	
	The regret of this policy can be decomposed as
	\begin{align}
	v^{\mathrm{HO}}-v^{\mathsf{IRT}^2}
	&=\E[V^{\mathrm{HO}}-V^{\mathsf{IRT}^2}]\nonumber\\
	&=\E[V^{\mathrm{HO}}-V^{\mathsf{HO}^1}]+\E[V^{\mathsf{HO}^1}-V^{\mathsf{HO}^2}]+\E[V^{\mathsf{HO}^2}-V^{\mathsf{IRT}^2}].\label{eq:hoandres2_decomp}
	\end{align}
	The first term of Eq~\eqref{eq:hoandres2_decomp} is bounded by $O(Te^{-\kappa T^{1/6}})$ as stated in Proposition~\ref{prop: resolve once}. For the second term of Eq~\eqref{eq:hoandres2_decomp}, the revenue of  policy $\mathsf{HO}^1$ and $\mathsf{HO}^2$ are
\[
V^{\mathsf{HO}^1} = V^{\mathsf{SPA}'}(0,t^*_1) + V^{\mathrm{HO}}(t^*_1,T), \quad V^{\mathsf{HO}^2} = V^{\mathsf{SPA}'}(0,t^*_1) +V^{\mathsf{SPA}'}(t_1^*,t_2^*)+V^{\mathrm{HO}}(t_2^*,T).
\]
Note that 	policy $\mathsf{HO}^1$ and $\mathsf{HO}^2$ are exactly the same during time $t\in[0,t^*_1)$,
	so we have	
	\begin{align}
	V^{\mathsf{HO}^1}-V^{\mathsf{HO}^2}
	&=V^{\mathrm{HO}}(t_1^*,T)-(V^{\mathsf{SPA}'}(t_1^*,t_2^*)+V^{\mathrm{HO}}(t_2^*,T))\nonumber.
	\end{align}
	Applying part (1) of Proposition \ref{prop: resolve once} to the remaining problem in $(t^*_1,T]$, which has a horizon length $T-t^*_1=\tau_1=T^{5/6}$,  we get
	\begin{align}
	\E[V^{\mathrm{HO}}(t_1^*,T)-(V^{\mathsf{SPA}'}(t_1^*,t_2^*)+V^{\mathrm{HO}}(t_2^*,T))]=O(T^{5/6} e^{-\kappa (T^{5/6})^{1/6}})
	= O(T^{5/6} e^{-\kappa T^{5/36}}).\label{eq:ho1andho2 result}
	\end{align}
	Because the revenue of $\mathsf{IRT}^2$ can be decomposed as
	\begin{align}
	V^{\mathsf{IRT}^2}&=V^{\mathsf{SPA}'}(0,t_1^*)+V^{\mathsf{SPA}'}(t_1^*,t_2^*)+V^{\mathsf{SPA}}(t_2^*,T),\nonumber
	\end{align}
	the last term of \eqref{eq:hoandres2_decomp} can be bounded by
	\begin{align}
	\E[V^{\mathsf{HO}^2}-V^{\mathsf{IRT}^2}]
	&=\E[V^{\mathrm{HO}}(t_2^*,T)-V^{\mathsf{SPA}}(t_2^*,T)]\nonumber\\
	&= O(\sqrt{T-t_2^*})\nonumber\\
	&= O(T^{(5/6)^2/2}).\label{eq:ho2andres2-result}
	\end{align}
Eq \eqref{eq:ho2andres2-result} follows the well-known result that static probabilistic allocation has a regret of $O(\sqrt{k})$ for a problem with horizon length $k$ (see Appendix~\ref{app:spa}). 
	
	Combining \eqref{eq:ho1andho2 result} and \eqref{eq:ho2andres2-result},  Eq~\eqref{eq:hoandres2_decomp} is bounded by
	\begin{align}
	v^{\mathrm{HO}}-v^{\mathsf{IRT}^2}=O(Te^{-\kappa T^{1/6}})+O( T^{5/6} e^{-\kappa T^{5/36}})+O(T^{25/72}).\label{eq:hoandres2 result}
	\end{align}
	
	Now, consider policy $\mathsf{IRT}^3$, which follows $\mathsf{IRT}^2$ during $t\in[0,t^*_3)$, but re-solves again at time $t_3^*=T- T^{(5/6)^3}$. By the same decomposition argument, the expected regret is given by
	\begin{align}
	v^{\mathrm{HO}}-v^{\mathsf{IRT}^3}
	&=\E[V^{\mathrm{HO}}-V^{\mathsf{HO}^1}]+\E[V^{\mathsf{HO}^1}-V^{\mathsf{HO}^2}]+\E[V^{\mathsf{HO}^2}-V^{\mathsf{HO}^3}]+\E[V^{\mathsf{HO}^3}-V^{\mathsf{IRT}^3}]\nonumber\\
	&=O(T e^{-\kappa T^{1/6}})+O(T^{5/6} e^{-\kappa T^{5/36}})+O(T^{25/36} e^{-\kappa T^{25/216}})+O(T^{(5/6)^3/2}).\nonumber
	\end{align}

{\color{black}	
Let
$
	K=\left\lceil\frac{\log \log T}{\log (6/5)}\right\rceil.
$
Note that the policy $\mathsf{IRT}^K$ coincides with \textsf{IRT}.
	By induction, if the decision maker re-solves $K$ times, where the $u$th re-solving time is $t^*_u = T - T^{(5/6)^u}$, the regret is given by
	{\color{black} 
	\begin{align}
	v^{\mathrm{HO}}-v^{\mathsf{IRT}} = v^{\mathrm{HO}}-v^{\mathsf{IRT}^K}=\sum_{u=0}^{K-1} 
	O\left((T^{(5/6)^u} \exp\left(-\kappa T^{(5/6)^u/6}\right)\right) 
	+ O( T^{(5/6)^K/2}).\label{eq:induction-bound}
	\end{align}}	
For the first term of the right hand side of Equation \eqref{eq:induction-bound}, using the fact that $T^{(5/6)^K} \leq e$,
we have
	\begin{align*}
	 \sum_{u=0}^{K-1} 	
	T^{(5/6)^u} \exp\left(-\kappa T^{(5/6)^u/6}\right) 
	= & \sum_{\ell=1}^{K} 	
	T^{(5/6)^{K-\ell}} \exp\left(-\kappa T^{(5/6)^{K-\ell}/6}\right) \\
	\leq &  \sum_{\ell=1}^{K} e^{(6/5)^{\ell}} \exp\left(-\kappa e^{(6/5)^{\ell}/6}\right) \\
	\leq &  \sum_{\ell=1}^{\infty} e^{(6/5)^{\ell}} \exp\left(-\kappa e^{(6/5)^{\ell}/6}\right) \\
 \leq & \int_{0}^{\infty} x \exp\left(-\kappa x^{1/6}\right) = O(1).
	\end{align*}
	Thus, the first term of the right hand side of Equation \eqref{eq:induction-bound} is $O(1)$. 
By the definition of constant $K$, we have $T^{(5/6)^K/2}\leq e^{1/2}$.
	Thus, the second term in \eqref{eq:induction-bound} is also $O(1)$.  Therefore, we have $	v^{\mathrm{HO}}-v^{\mathsf{IRT}} = O(1)$.
	In addition, this constant factor is
	independent of the time horizon $T$ and the capacity vector $C$.}
\Halmos	
\endproof

\subsection{Proof of Proposition \ref{prop: resolve once}} \label{sub: proof_proposition}
Throughout this subsection, we focus on policies $\textsf{IRT}^1$ and $\textsf{HO}^1$ with only one re-solving at $t^*_1 = T - T^{5/6}$. We write $t^*:=t^*_1$ for simplicity.
Let $\Gamma(T)=\alpha\sum_{j:x^*_j=\lambda_j}|\Lambda_j(T)-\lambda_jT|$, where $\alpha$ is a constant whose value is determined by the BOM matrix $A=(a_{lj})_{l\in[m], j\in[n]}$. 
More specifically, $\alpha$
is the maximum absolute value of the elements in the inverses of all invertible submatrices of the BOM matrix $A$.
In a special case when all entries of $A$ are either 0 or 1, we have $\alpha \leq \max\{1, m\wedge n-1\}$.
We let $\Delta_j(t)$ be the deviation of the number of arrivals of class $j$ customer from its mean in $(t,T]$, i.e., $\Delta_j(t)=\Lambda_j(T)-\Lambda_j(t)-\lambda_j(T-t)$. 
Define $\tilde{z}_j(t)$  as the number of class $j$ customers accepted up to time $t$ if the algorithm were allowed to go over the capacity limits.  
For all $j\in[n]$, we define the following events:
\begin{align}
E_{1,j} &=\left\{ (T-t^*)x^*_j-\tilde{z}_j(t^*)+t^*x^*_j \geq \Gamma(T)\right\},\label{event: E1}\\
E_{2,j} &=\left\{ (T-t^*)(\lambda_j-x^*_j)+\tilde{z}_j(t^*)-t^*x^*_j \geq \Gamma(T)+|\Delta_j(t^*)|\right\}. \label{event: E2}
\end{align}
The event E is defined as
\begin{align}
E=\left(\bigcap_{j:x^*_j\geq\lambda_j T^{-1/4}} E_{1,j} \right) \cap \left( \bigcap_{j:x^*_j\leq \lambda_j(1-T^{-1/4})} E_{2,j} \right) .\label{event: E}
\end{align}

Now we will prove Proposition \ref{prop: resolve once}.
\proof{Proof of Proposition \ref{prop: resolve once}.}
For all $j: x^*_j< \lambda_j T^{-1/4}$, $\tilde{z}_j(t^*) =0\leq Tx^*_j$.
For all $j: x^*_j\geq\lambda_j T^{-1/4}$,	event $E$ in \eqref{event: E} implies
$
\tilde{z}_j(t^*) \leq  (T-t^*)x^*_j + t^*x^*_j = Tx^*_j.
$
So, suppose  event $E$ holds,  the capacity constraints for all resources are satisfied up to period $t^*$, and we have $z_j(t^*)=\tilde{z}_j(t^*)$. 	

If $x^*_j<\lambda_j T^{-1/4}$, we have $\bar{z}_j-z_j(t^*)=\bar{z}_j-0\geq0$. (Recall that $\bar{z}_j$ is the solution to the hindsight optimum; see Section~\ref{subsec: hindsight}.) Otherwise,
suppose  event $E$ holds, by Lemma~\ref{lem: range hindsight} in Appendix~\ref{app:lemmas},
we have
\begin{align}
\bar{z}_j-z_j(t^*)&\geq Tx^*_j-\Gamma(T)-z_j(t^*)+t^*x^*_j-t^*x_j^* = (T-t^*)x^*_j-\Gamma(T)-(z_j(t^*)-t^*x_j^*) \geq 0, \label{eq: def of E1}
\end{align}
where \eqref{eq: def of E1} follows from the condition \eqref{event: E1} and the fact that $z_j(t^*)=\tilde{z}_j(t^*)$.

Similarly, if $x^*_j>\lambda_j(1-T^{-1/4}) $, we  have $z_j(t^*)+\Lambda_j(T)-\Lambda_j(t^*)-\bar{z}_j \geq \Lambda_j(t^*) +\Lambda_j(T)-\Lambda_j(t^*)  - \Lambda_j(T) =0$. Otherwise, suppose  event $E$ holds, by Lemma~\ref{lem: range hindsight}, we have
\begin{align}
z_j(t^*)+\Lambda_j(T)-\Lambda_j(t^*)-\bar{z}_j&=z_j(t^*)+(T-t^*)\lambda_j+\Delta_j(t^*)-\bar{z}_j\nonumber\\
&\geq z_j(t^*)+(T-t^*)\lambda_j+\Delta_j(t^*)-Tx^*_j-\Gamma(T)\nonumber\\
&=(z_j(t^*)-t^*x^*_j)+(T-t^*)(\lambda_j-x^*_j)+\Delta_j(t^*)-\Gamma(T)\nonumber\\
&\geq (T-t^*)(\lambda_j-x^*_j)+(z_j(t^*)-t^*x^*_j)-|\Delta_j(t^*)|-\Gamma(T) \nonumber\\
&\geq 0,  \label{eq: def of E2}
\end{align}
where \eqref{eq: def of E2} follows from the condition \eqref{event: E2}  and the fact that $z_j(t^*)=\tilde{z}_j(t^*)$.  

Therefore, combining \eqref{eq: def of E1} and \eqref{eq: def of E2}, we have $z_j(t^*)\leq \bar{z}_j \leq z_j(t^*)+\Lambda_j(T)-\Lambda_j(t^*)$. In other words, the decision maker would still be able to achieve the hindsight optimum if she uses probabilistic allocation up to $t^*$, and then 
gets perfect information from $t^*$ onwards, because she can accept $\bar{z}_j -z_j(t^*) $ of class $j$ customers. If the decision maker re-solves once at $t^*$, then the regret can be written as
\begin{align}
v^{\mathrm{HO}}-v^{\mathsf{IRT}^1}
=\E[V^{\mathrm{HO}}-V^{\mathsf{IRT}^1}] =\E[V^{\mathrm{HO}}-V^{\mathsf{HO}^1}]+\E[V^{\mathsf{HO}^1}-V^{\mathsf{IRT}^1}].\label{eq:hoandres1_decomp}
\end{align}
where we can decompose the revenue earned under the policy $\mathsf{HO}^1$ and $\mathsf{IRT}^1$ as
\begin{align}
V^{\mathsf{HO}^1}=V^{\mathsf{SPA}'}(0,t^*)+ V^{\mathrm{HO}}(t^*,T),\quad V^{\mathsf{IRT}^1}=V^{\mathsf{SPA}'}(0,t^*)+V^{\mathsf{SPA}}(t^*,T).\nonumber
\end{align}
Consequently, we get
\begin{align}
\E[V^{\mathsf{HO}^1}-V^{\mathsf{IRT}^1}]
=\E[V^{\mathrm{HO}}(t^*,T)-V^{\mathsf{SPA}}(t^*,T)] = O(\sqrt{T-t^*})=O( T^{5/12}).\label{eq:ho1andres1 result}
\end{align}
Eq  \eqref{eq:ho1andres1 result} follows the well-known result that static probabilistic allocation without re-solving has a regret rate of $O(\sqrt{k})$ for a problem with horizon length $k$, where the constant factor does not depend on the capacity vector $C$ \citep[see e.g.][]{ReimanAsymptoticallyOptimalPolicy2008}. For completeness, we give a proof of this result in Appendix~\S\ref{append:additional-result}.

Recall that if the event $E$ happens, the hindsight optimal is still attainable starting from $t^*$. In other words, conditioned on $E$, the regret of $\mathsf{HO}^1$ is $V^{\mathrm{HO}}-V^{\mathsf{HO}^1}=0$. Therefore, the first term of \eqref{eq:hoandres1_decomp} is given by
\begin{align}
\E[V^{\mathrm{HO}}-V^{\mathsf{HO}^1}]
&=\E[V^{\mathrm{HO}}-V^{\mathsf{HO}^1} \mid E]\Prob(E)+\E[V^{\mathrm{HO}}-V^{\mathsf{HO}^1} \mid E^\complement]\Prob(E^\complement)\nonumber\\
&=\E[V^{\mathrm{HO}}-V^{\mathsf{HO}^1} \mid E^\complement]\Prob(E^\complement)\nonumber\\
&\leq \E[V^{\mathrm{HO}} \mid E^\complement]\Prob(E^\complement). \label{eq:hoandho1}
\end{align}
{\color{black}
Note that the hindsight optimum is bounded almost surely by
$
V^{\mathrm{HO}} \leq \sum_{j=1}^n r_j \Lambda_j(T),
$
where $\Lambda_j(T)$ is the total number of arrivals from class $j$. Moreover, $\Lambda_j(T)$ follows Poisson distribution with mean $\lambda_j T$. 
By the Poisson tail bound (see Lemma~\ref{lem: poisson bound} in Appendix~\ref{app:lemmas}), we have
\begin{align}
	\E[ ( \Lambda_j(T) - 2 \lambda_j T)^+ ] &= \int_{0}^{\infty} \Prob (\Lambda_j(T) - 2\lambda_j T \geq x ) dx 
	\leq \int_{0}^{\infty}2\exp\left(-\frac{(x+\lambda_jT)^2}{3\lambda_jT}\right)dx\nonumber\\
	&\leq \int_{0}^{\infty} 2\exp\left(-\frac{(x+\lambda_j T)\lambda_j T}{3\lambda_jT}\right)dx = 6\exp\left(-\frac{\lambda_j T}{3}\right).\nonumber
\end{align}
Combining the above inequality with Eq~\eqref{eq:hoandho1}, we have
\begin{align}
\E[V^{\mathrm{HO}}-V^{\mathsf{HO}^1}]
&\leq \E[V^{\mathrm{HO}} \mid E^\complement]\Prob(E^\complement) \nonumber\\
& \leq \E\left[\sum_{j=1}^n r_j \Lambda_j(T) \Big| E^\complement \right]\Prob(E^\complement) \nonumber \\
& \leq \E\left[\sum_{j=1}^n r_j  \left( 2 \lambda_j T +(  \Lambda_j(T) - 2 \lambda_j T)^+ \right) \Big| E^\complement \right]\Prob(E^\complement)  \nonumber\\
& \leq  \sum_{j=1}^n r_j \left( 2  \lambda_j T \Prob(E^\complement)  
 +   \E [ ( \Lambda_j(T) - 2 \lambda_j T)^+ | E^\complement] \Prob(E^\complement) \right)  \nonumber\\
 & \leq  \sum_{j=1}^n r_j \left( 2  \lambda_j T \Prob(E^\complement)  
 +   \E [ ( \Lambda_j(T) - 2 \lambda_j T)^+] \right)  \nonumber\\
 & \leq \sum_{j=1}^n  r_j  \left( 2\lambda_j T  \Prob(E^\complement)  + 6 e^{-\lambda_jT/3} \right).   \label{eq:hobound}
\end{align}
}
Using the result from Lemma~\ref{lem: prob E'} in Appendix~\ref{app:lemmas} and \eqref{eq:hoandho1}--\eqref{eq:hobound}, we have
\begin{align}
\E[V^{\mathrm{HO}}-V^{\mathsf{HO}^1}]
\leq \sum_{j=1}^n  r_j  \left( 2\lambda_j T \cdot O(e^{-\kappa T^{1/6}})  + 6 e^{-\lambda_jT/3} \right)
=O(Te^{-\kappa T^{1/6}}). \label{eq:hoandho1 result}
\end{align}
We can conclude from \eqref{eq:hoandres1_decomp}, \eqref{eq:ho1andres1 result} and \eqref{eq:hoandho1 result} that
\begin{align}
v^{\mathrm{HO}}-v^{\mathsf{IRT}^1} =O(Te^{-\kappa T^{1/6}})+O(T^{5/12}),\nonumber 
\end{align}
where the big O notation hides constants that are
	independent of the time horizon $T$ and the capacity vector $C$. \Halmos	
\endproof

\section{Proofs for Results in Section~\ref{sec: analysis of resolving}}

In this section, we provide complete proofs for the results on the \textsf{FR} policy in Section~\ref{sec: analysis of resolving}.

\subsection{Proof of Proposition~\ref{thm: ho res lb}}\label{subsec: proof of HO-FR}

\proof{Proof of Proposition~\ref{thm: ho res lb}.}
	Let $V^{\mathsf{FR}}$ denote the revenue earned under the $\mathsf{FR}$ policy.
	We know that
$
	v^{\mathrm{HO}}-v^{\mathsf{FR}}=\E[V^{\mathrm{HO}}]-\E[V^{\mathsf{FR}}]\nonumber =\E[V^{\mathrm{HO}}-V^{\mathsf{FR}}].
$
	From the law of total expectation, it equals to
	\begin{align}	\E[V^{\mathrm{HO}}-V^{\mathsf{FR}}|Q]\Prob(Q)+\E[V^{\mathrm{HO}}-V^{\mathsf{FR}}|Q^\complement]\Prob(Q^\complement),\label{eq:cond_E}
	\end{align}
	for any event $Q$. Since $V^{\mathrm{HO}}$ is an upper bound of $V^{\mathsf{FR}}$, i.e., $V^{\mathrm{HO}}\geq V^{\mathsf{FR}}$ a.s. and the probability of any measurable event is nonnegative, the second term of \eqref{eq:cond_E} is nonnegative. Consequently,
	\begin{align}
	v^{\mathrm{HO}}-v^{\mathsf{FR}}
	&\geq \E[V^{\mathrm{HO}}-V^{\mathsf{FR}}|Q]\Prob(Q).\nonumber
	\end{align}
	That is, to complete the proof, we want to show that $\E[V^{\mathrm{HO}}-V^{\mathsf{FR}}|Q]\Prob(Q)\geq \Omega(\sqrt{T})$ for some appropriately chosen event $Q$.
	
	{\color{black} Recall that we consider a problem instance with two classes of customers and one resource. Customers from each class arrive according to a Poisson process with rate 1. The arrivals from two classes are assumed to be independent. The initial resource capacity is $T$. Customers from both classes, if accepted, consume one unit of the resource, but pay different prices, $r_1>r_2$. To proof the result,}
	we will consider the situation when the number of arrivals of class 1 customers in $[0,T]$ is above its mean which is $T$. Since the initial capacity of the resource is $T$, the hindsight optimal policy is to accept only class 1 customer. More specifically, we should accept $T$ of class 1 customer and none of class 2 customer. Failing to do so will result in  a positive regret.
	
	We will partition time in the interval $[0,T]$ into 3 phases of equal length. Let $T'$ and $T''$ denote the beginning of phase 2 and phase 3 respectively. In other words, $T'=T/3$ and $T''=2T/3$. We will define the events of the number of the arrivals of class 1 customer in each period as follows.
	\begin{align}
	Q_1&=\{T'-{\color{black}4}\sqrt{T'}\leq \Lambda_1(0,T')\leq T'-{\color{black}3}\sqrt{T'}\},\label{event:A1}\\
	Q_2&=\{(t-T')-{\color{black}2}\sqrt{T'}\leq \Lambda_1(T',t)\leq (t-T')+{\color{black}2}\sqrt{T'},\; \forall t \in (T',T'']\},\label{event:A2}\\
	Q_3&=\{T'+{\color{black}6}\sqrt{T'}\leq \Lambda_1(T'',T)\leq T'+{\color{black}7}\sqrt{T'}\}\label{event:A3},
	\end{align}
	where $Q_i$ restricts the number of the arrivals of class 1 customer in phase $i$.
	 Let $z_j(t_1,t_2)$ denote the actual number of class $j$ customers admitted in $(t_1,t_2]$. We will further define the events of the number of accepted customers as
	\begin{align}
	B&=\{z_2(0,T)\geq \frac{1}{6}\sqrt{T'}\},\label{event:B}\\
	B_1&=\{z_1(0,T')< T'-{\color{black}10}\sqrt{T'}\}.\label{event:B1}
	\end{align}	
	
	If the event $B$ happens, the decision maker will admit at least $\frac{1}{6}\sqrt{T'}$ of class 2 customers which leads to a regret of at least $(r_1-r_2)\frac{1}{6}\sqrt{T'}$ from the hindsight optimal. On the other hand, if the event $B_1$ happens, the decision maker will admit less than $T'-{\color{black}10}\sqrt{T'}$ of class 1 customer; this means that even if the decision maker admit all arrivals of class 1 customer in the second and the third phase, the total number of admitted class 1 customer is
	\begin{align}
	z_1(0,T)<T'-{\color{black}10}\sqrt{T'}+T'+{\color{black}2}\sqrt{T'}+T'+{\color{black}7}\sqrt{T'}=T-\sqrt{T'},\nonumber
	\end{align}
	which results in a regret of at least $(r_1-r_2)\sqrt{T'}$ from the hindsight optimal. Therefore, if the event $B$ or $B_1$ happens with probability that is bounded away from zero, the incurred a regret of is at least $(r_1-r_2)\frac{1}{6}\sqrt{T'}$ from the hindsight optimal. If we can show that this event happens with positive probability, then we are done. That is, we want to show that
	\begin{align}
	\Prob((B\cup B_1)\cap Q_1 \cap Q_2 \cap Q_3)>0.\nonumber
	\end{align}
	This probability can be written as
	\begin{align}
	\Prob((B\cup B_1)\cap Q_1 \cap Q_2 \cap Q_3) &= \Prob(Q_1 \cap Q_2 \cap Q_3)-\Prob(B^\complement \cap B_1^\complement\cap  Q_1 \cap Q_2 \cap Q_3)\nonumber\\
	&= \Prob(Q_1)\Prob(Q_2)\Prob(Q_3)-\Prob(B^\complement \cap B_1^\complement\cap  Q_1 \cap Q_2 \cap Q_3)\label{eq:Aindep}\\
	&\geq \Prob(Q_1)\Prob(Q_2)\Prob(Q_3)-\Prob(B^\complement \cap B_1^\complement\cap  Q_1 \cap Q_2),\label{eq:subset}
	\end{align}
	where the first term of \eqref{eq:Aindep} follows since the arrivals of the customers in disjoint interval are independent, i.e., $Q_1, Q_2$ and $Q_3$ are independent, and the inequality \eqref{eq:subset} follows because the event $B^\complement \cap B_1^\complement \cap  Q_1 \cap Q_2  \cap Q_3$ is a subset of the event $B^\complement \cap B_1^\complement\cap  Q_1 \cap Q_2$. The remaining part of the proof relies on the results of Lemma~\ref{lemma: prob events in HO-FR} in Appendix~\ref{app:lemmas}.
	Combining the result from Lemma~\ref{lemma: prob events in HO-FR} and \eqref{eq:subset}, we get
	\begin{multline}
	\mathbb{P}((B\cup B_1)\cap  Q_1 \cap Q_2\cap Q_3)
	\geq \left({\color{black}0.0013}-\frac{0.9496}{\sqrt{T'}}\right)\left(\color{black}0.5\right)\left({\color{black}9.8531 \times 10^{-10}}-\frac{0.9496}{\sqrt{T'}}\right)-{\color{black}e^{-0.0026\sqrt{T'}}}.\nonumber
	\end{multline}
	Therefore, the regret of \textsf{FR} is given by
	\begin{align*}
	v^{\mathrm{HO}}-v^{\color{black}\mathsf{FR}} &\geq {\color{black}(r_1-r_2)}\frac{1}{6}\sqrt{T'}\left(\left({\color{black}0.0013}-\frac{0.9496}{\sqrt{T'}}\right)\left(\color{black}0.5\right)\left({\color{black}9.8531 \times 10^{-10}}-\frac{0.9496}{\sqrt{T'}}\right)-{\color{black}e^{-0.0026\sqrt{T'}}}\right)\nonumber\\
= &\Omega (\sqrt{T}).\nonumber \Halmos	
	\end{align*}
\endproof

\subsection{Proof of Proposition \ref{thm: DLP and resolving ub}}\label{subsec: proof DLP-FR}
\proof{Proof Proposition \ref{thm: DLP and resolving ub}.}

The \textsf{FR} algorithm divides the horizon $[0,T]$ into $T$ periods: $[0,1) \cup [1,2) \cup \cdots \cup [T-1,T]$. In each period $[t,t+1)$, the algorithm attempts to accept class $j$ customers with probability $x_j(t)/\lambda_j$.
If we ignore capacity constraints, the algorithm on average accepts $x_j(t)$ customers from class $j$.
However, the decision maker can potentially reject customers due to capacity constraints; if that happens, the expected number of admitted class $j$ customers in period $t$ is less than $x_j(t)$ per period.

More specifically, let us consider the LP solved at time $t$:
\[
\max_y \Big\{\sum_{j=1}^{n}r_j y_j \;\Big|\; \sum_{j=1}^{n}A_j y_j\leq C(t), \text{ and }\, 0\leq y_j \leq \lambda_j(T-t), \forall  j \in [n] \Big\}.
\]
Let $y(t)$ be the optimal solution.
The algorithm accepts customers in class $j$ with probability {\color{black}$\frac{x_j(t)}{\lambda_j}=\frac{y_j(t)}{\lambda_j(T-t)}$}.
{\color{black} If there is insufficient capacity for class $j$ in period $[t,t+1)$, two cases can happen: case 1) $A_j \nleq C(t)$, i.e., there exists $l\in [m]$ such that $a_{lj}>C_l(t)$, so there is insufficient capacity for class $j$ when this period starts; 
case 2) $A_j \leq C(t)$, namely there is sufficient capacity when this period starts, but during $[t,t+1)$ the capacity of certain resource runs out, and a class $j$ customer that arrives at time $s\in[t,t+1)$ finds $A_j \nleq C(s)$. For case 1), $A_j \nleq C(t)$ implies $y_j(t) < 1$. So, the expected number of class $j$ customers that $\mathsf{FR}$ would accept, but that are not accepted because of the capacity constraint}
in period $[t,t+1)$ is less than $\lambda_j \cdot \frac{1}{\lambda_j(T-t)}=\frac{1}{T-t}$. 
The revenue loss from that group of customers over the entire horizon is bounded by
\[
\sum_{j=1}^n \sum_{t=0}^{T-1} r_j \frac{1}{T-t} \leq  \sum_{j=1}^n r_j\left(\log T+1\right),
\]
where we use the fact that $\sum_{i=2}^T \frac{1}{i}\leq
\int_{x=1}^{T} \frac{1}{x} dx=\log (T)$.
{\color{black}For case 2), we note that this situation  can happen \emph{at most once} during the entire horizon for each class of customers. Since the expected number of class $j$ customers that $\mathsf{FR}$ would accept, but that are not accepted because of the capacity constraint in one period is bounded above by the expected number of class $j$ arrivals in one period which is $\lambda_j$, the revenue loss caused by case 2) is bounded by $\sum_{j=1}^nr_j\lambda_j.$ 
}
We can write the expected revenue of the re-solving heuristic as
\begin{align*}
v^{\mathsf{FR}}\geq\E\left[\sum_{t=0}^{T-1}\sum_{j=1}^n r_jx_j(t)\right]-\sum_{j=1}^n r_j \left(\log T+1\right){\color{black}-\sum_{j=1}^nr_j\lambda_j},
\end{align*}
{\color{black}where the last two terms account for the lost sales of case 1) and case 2) respectively.}
Therefore, the expected revenue loss of the re-solving heuristic can be bounded by
\begin{align}
v^{\mathrm{DLP}}-v^{\mathsf{FR}}\leq T\sum_{j=1}^nr_jx_j^*-\mathbb{E}\left[\sum_{t=0}^{T-1}\sum_{j=1}^nr_jx_j(t)\right]+\sum_{j=1}^n r_j\left(\log T+1\right){\color{black}+\sum_{j=1}^nr_j\lambda_j}. \label{eq:revloss_dlpres}
\end{align}
{\color{black}Since the solutions to the $\mathrm{DLP}(x^*)$ and the LP solved under the $\mathsf{FR}$ policy$(x(t))$ only differ in the right hand side of the capacity constraints ($b$ and $b(t)$), we can write,}
for any time period $[t,t+1)$,
\begin{align}
\sum_{j=1}^nr_jx_j^*-\sum_{j=1}^nr_jx_j(t) \leq \sum_{l=1}^{L}r_{max}^l(b_l-b_l(t))^+,\label{eq:diff DLP and FR}
\end{align}
where $r_{max}^l$ is the largest possible revenue gain by increasing the capacity of resource $l$ by one unit, i.e., $r^l_{max}=\max_{j\in[n]}\{
r_j \mathrm{I}(a_{lj}>0)/a_{lj}\}$. {\color{black} Equation~\eqref{eq:diff DLP and FR} holds because $r_{max}^l$ is an upper bound of the dual price for resource $l\in[m]$.} From the definition of $b_l(i)$, we can write
\begin{align}
(b_l-b_l(t))^+&=\left[\sum_{i=0}^{t-1}(b_l(i)-b_l(i+1))\right]^+ =\left[\sum_{i=0}^{t-1}\left(\frac{C_l(i)}{T-i}-\frac{C_l(i+1)}{T-i-1}\right)\right]^+. \nonumber
\end{align}
Note that
$
\frac{1}{T-i}=\frac{1}{T-i-1}-\frac{1}{(T-i-1)(T-i)},
$
so it follows that
\begin{align}
(b_l-b_l(t))^+
&=\left[\sum_{i=0}^{t-1}\left(\frac{C_l(i)}{T-i-1}-\frac{C_l(i)}{(T-i-1)(T-i)}-\frac{C_l(i+1)}{T-i-1}\right)\right]^+ \nonumber\\
&=\left[\sum_{i=0}^{t-1}\left(\frac{C_l(i)-C_l(i+1)}{T-i-1}-\frac{C_l(i)}{(T-i-1)(T-i)}\right)\right]^+. \label{eq:capacity-bound}
\end{align}
Let $z_j(t)$ be the actual number of class $j$ customers admitted in $[0,t]$. The change in the capacity of resource $l$ is given by
$
C_l(i)-C_l(i+1)=\sum_{j=1}^n a_{lj}(z_j(i+1)-z_j(i)).
$
Because of the capacity constraint, the decision maker may fail to accept some customers. Therefore, the actual number of the admitted customers in any period is bounded above by the number of the customers admitted in that period by ignoring the capacity constraint. 
{\color{black}
More specifically, we define stochastic processes $\{\tilde{z}_j(t), t\geq 1, j\in[n]\}$, such that $\tilde{z}_j(t+1)-\tilde{z}_j(t)$ follows Poisson distribution with mean $x_j(t)$ (the solution to the LP at period $t$). Therefore,  $\tilde{z}_j(t+1)-\tilde{z}_j(t)$ is the the number of class $j$ customers that the algorithm could have admitted if there were no capacity constraint in $[t,t+1)$.
Since the number of customers who are actually admitted, $z_j(t+1)-z_j(t)$, follows the same Poisson distribution with additional rejections due to capacity constraints, we always have
$
z_j(t+1)-z_j(t)\leq \tilde{z}_j(t+1)-\tilde{z}_j(t), 
$
and therefore
\[C_l(i)-C_l(i+1)=\sum_{j=1}^n a_{lj}(z_j(i+1)-z_j(i))\leq \sum_{j=1}^n a_{lj}(\tilde{z}_j(i+1)-\tilde{z}_j(i)).\]
}
We can now bound \eqref{eq:capacity-bound} by
\begin{align}
(b_l-b_l(t))^+
&\leq\left[\sum_{i=0}^{t-1}\left(\frac{\sum_{j=1}^n a_{lj}(\tilde{z}_j(i+1)-\tilde{z}_j(i))}{T-i-1}-\frac{b_l(i)}{T-i-1}\right)\right]^+ \label{eq:lost sales}\\
&\leq\left[\sum_{i=0}^{t-1}\left(\frac{\sum_{j=1}^n a_{lj}(\tilde{z}_j(i+1)-\tilde{z}_j(i))}{T-i-1}-\frac{\sum_{j=1}^na_{lj}x_j(i)}{T-i-1}\right)\right]^+, \label{eq:def of resolving}
\end{align}
where the second term in the RHS of \eqref{eq:lost sales} follows from the definition of $b_l(i)= C_l(i)/(T-i)$ and 
the inequality of \eqref{eq:def of resolving} follows from the definition of the re-solving LP in Algorithm~\ref{algo:res} which is $\sum_{j=1}^n a_{lj}x_j(i)\leq b_l(i)$.

We will use the result from Lemma~\ref{lem: expectation b} in Appendix~\ref{app:lemmas}
to finish the proof.
{\color{black}If we sum the inequality \eqref{eq:diff DLP and FR} over $t\in\{0\}\cup[T-1]$ and apply Lemma~\ref{lem: expectation b}, it follows that}
\begin{align}
T\sum_{j=1}^nr_jx_j^*-\E\left[\sum_{t=0}^{T-1}\sum_{j=1}^nr_jx_j(t)\right] &\leq \E\left[\sum_{t=0}^{T-1}\sum_{l=1}^{L}r_{max}^l(b_l-b_l(t))^+\right]\nonumber \\
&\leq \sum_{l=1}^L r_{max}^lK_l\sum_{t=0}^{T-1}\sqrt{\sum_{i=0}^{t-1}\frac{1}{(T-i-1)^2}} \nonumber\\
&= \sum_{l=1}^L r_{max}^lK_l\sum_{t=1}^{T}\sqrt{\sum_{i=1}^{t-1}\frac{1}{(T-i)^2}} \nonumber\\
&=\sum_{l=1}^L r_{max}^lK_l\left[\sum_{t=1}^{T-1}\sqrt{\sum_{i=1}^{t-1}\frac{1}{(T-i)^2}}+\sqrt{\sum_{i=1}^{T-1}\frac{1}{(T-i)^2}}\right]\nonumber\\
&\leq \sum_{l=1}^L r_{max}^lK_l \left[\sum_{t=1}^{T-1}\sqrt{\int_{1}^{t}\frac{1}{(T-s)^2}ds}+\sqrt{\int_{1}^{T-1}\frac{1}{(T-s)^2}ds+1}\right]\nonumber\\
&\leq \sum_{l=1}^L r_{max}^lK_l\left[\sum_{t=1}^{T-1}\sqrt{\frac{1}{T-t}}+\sqrt{1-\frac{1}{T-1}+1}\right]\nonumber\\
&\leq \sum_{l=1}^L r_{max}^lK_l(2\sqrt{T}+\sqrt{2}).\nonumber
\end{align}
Combining the result to \eqref{eq:revloss_dlpres}, we can conclude that
\begin{align}
v^{\mathrm{DLP}}-v^{\mathsf{FR}}
&{\color{black}\leq\sum_{l=1}^L r_{max}^lK_l(2\sqrt{T}+\sqrt{2})+n r_{max} \left(\log T+1\right)+ n r_{max}\lambda_{max}} =O(\sqrt{T}),\nonumber
\end{align}
{\color{black}where $r_{max}=\max_{j\in[n]}r_j$ and $\lambda_{max}=\max_{j\in[n]}\lambda_j.$}  \Halmos	
\endproof

{
\section{Lemmas}
\label{app:lemmas}
{\color{black}
	\begin{lemma}[Berry-Esseen theorem, Corollary 1 in \citet{shevtsova2011absolute}]\label{lem:berry-esseen}
	Let $X_1,X_2,\ldots$ be independent and identically distributed random variables (i.i.d.) with $\E[X_1]=0,\;\E[X_1^2]=\sigma^2>0$ and $\E[|X_1^3|]=\rho< \infty$. Let $F_n$ be the cumulative distribution function of $\frac{X_1+\ldots+X_n}{\sigma\sqrt{n}}$ and $\Phi$ the cumulative distribution function of standard normal distribution. For any $x$ and $n$,
	\[
	|F_n(x)-\Phi(x)|\leq \frac{0.4748\rho}{\sigma^3\sqrt{n}}.
	\]
\end{lemma}
}

{\color{black}
\begin{lemma}[Doob's maximal inequality]\label{lem:doob}
	Suppose $M_t$ is a martingale with paths that are right continuous with left limits. Then, for any constant $a>0$,
	\[
		\Prob(\sup_{s\leq t}|M_s|\geq a)\leq \frac{\E[|M_t|]}{a}.
	\]
\end{lemma}

\begin{lemma}[\citet{freedman1975}]\label{lem:freedman}
Given a sequence of real-valued supermartingale differences $(\xi_i,\mathcal{F}_i)_{i\in\{0\}\cup [n]}$ with $\xi_0=0.$ Set $S_k=\sum_{i=0}^k\xi_i$ for $k \in [n]$. Then $S=(S_k,\mathcal{F}_k)_{k \in [n]}$ is a supermartingale. Let $\langle S\rangle_k=\sum_{i=1}^k\E[\xi_i^2|\mathcal{F}_{i-1}]$. Suppose $\xi_i\leq 1$. Then, for all $x,v>0$,
	\begin{align*}
	\mathbb{P}(\exists k \in [n]: S_k\geq x, \langle S\rangle_k \leq v^2)\leq \exp\left(-\frac{x^2}{2(v^2+x)}\right).
	\end{align*}
\end{lemma}

\begin{lemma}[Poisson Tail Bound, \citet{pollard2015inequality}]\label{lem: poisson bound}
	If random variable $X$ follows a Poisson distribution with parameter $\lambda>0$, for any constant $x>0$, we have
	\[ 
	\Prob(|X-\lambda|\geq x) \leq 2\exp\left(-\frac{x^2}{3 \lambda }\right).
	\]
\end{lemma} 
}

\begin{lemma} \label{lem: range hindsight}
	There exists an optimal solution to the hindsight LP defined in \eqref{form: hindsight}, $\bar{z}=(\bar{z}_1,\ldots,\bar{z}_n)$, such that
	\begin{equation}\label{eq:hindsight-interval}
	\bar{z}_j \in [\max\{Tx^*_j-\Gamma(T), 0\}, \min \{Tx^*_j+\Gamma(T), \Lambda_j(T)\}] \quad\text{for all }j\in[n].
	\end{equation}
\end{lemma}
\proof{Proof of Lemma~\ref{lem: range hindsight}.}
Theorem 4.2 in \citet{ReimanAsymptoticallyOptimalPolicy2008} shows that
there exists an optimal solution to the hindsight LP defined in \eqref{form: hindsight}, $\bar{z}=(\bar{z}_1,\ldots,\bar{z}_n)$, such that
\[
\bar{z}_j \in [Tx^*_j-\Gamma(T), Tx^*_j+\Gamma(T)] \quad\text{for all }j\in[n].
\]
The lemma immediately follows since any feasible solution to the hindsight LP satisfies $0\leq \bar{z}_j \leq \Lambda_j(T)$ for all $j\in [n]$.
\Halmos
\endproof

\begin{lemma}\label{lem: prob E'} The probability of event $E$ defined in \eqref{event: E}  satisfies
	$
	\Prob(E^\complement)=O\left(n |J_\lambda| \exp \bigl(	-\kappa T^{1/6}	\bigr)\right),
	$
	where the constant $\kappa$ is given by
	{\color{black}$
	\kappa= \frac{\lambda_{min}}{27(\alpha|J_\lambda|+1)^2}.
	$}
	The set $J_\lambda=\{j: x^*_j=\lambda_j\}$, and $\alpha$ is a positive constant that depends on the BOM matrix $A$.
\end{lemma}

\proof{Proof of Lemma \ref{lem: prob E'}.}
We can write
\begin{align}
\Prob(E^\complement)
&=\Prob(\bigcup_{j:x^*_j\geq\lambda_j T^{-1/4}} E_{1,j}^\complement \cup \bigcup_{j:x^*_j\leq\lambda_j(1-T^{-1/4})}E_{2,j}^\complement)\nonumber\\
&\leq \sum_{j:x^*_j\geq\lambda_j T^{-1/4}} \Prob(E_{1,j}^\complement)+\sum_{j:x^*_j\leq\lambda_j(1-T^{-1/4})}\Prob(E_{2,j}^\complement).\label{eq:p_E'}
\end{align}
First, we will bound
\begin{align}
\Prob(E_{1,j}^\complement)
&=\Prob(\tilde{z}_j(t^*)-t^*x^*_j+\Gamma(T)>(T-t)x^*_j).\nonumber
\end{align}	
Observe that if the event $E_{1,j}^\complement$ happens, at least one of the following two events must happen:
\begin{align}
\left\{\tilde{z}_j(t^*)-t^*x^*_j>\frac{(T-t^*)x^*_j}{2}\right\} ,\;\left\{\Gamma(T)>\frac{(T-t^*)x^*_j}{2}\right\}.\nonumber \nonumber
\end{align}
Thus, we can apply the union bound and write
\begin{align}
\Prob(E_{1,j}^\complement)
&\leq \Prob\left(\tilde{z}_j(t^*)-t^*x^*_j>\frac{(T-t^*)x^*_j}{2}\right) + \Prob\left(\Gamma(T)>\frac{(T-t^*)x^*_j}{2}\right).\label{eq: union bound E1'}
\end{align}

Let us consider three cases: 1) $x^*_j < \lambda_j T^{-1/4}$, 2) $\lambda_j T^{-1/4} \leq x^*_j \leq \lambda_j (1 - T^{-1/4})$, 3) $x^*_j > \lambda_j (1 - T^{-1/4})$. 
Case 1) is already eliminated in definition of event $E$, so we focus on case 2) and 3).
In Case 2), $\tilde{z}_j(t^*)$ is Poisson random variable with parameter $t^*x^*_j$.
We use the Poisson tail bound (Lemma \ref{lem: poisson bound}) to bound such events.

It follows from Lemma \ref{lem: poisson bound} that
\begin{align}
\Prob\left(\tilde{z}_j(t^*)-t^*x^*_j>\frac{(T-t^*)x^*_j}{2}\right)
&\leq 2 \exp \left( -\frac{1}{3}\left(\frac{(T-t^*)x^*_j}{2}\right)^2\frac{1}{t^*x^*_j}\right)\nonumber\\
&\leq 2 \exp \left(-\frac{(T-t^*)^2 x_j^{*}}{12 T}\right)\nonumber\\
&\leq 2 \exp \left(-\frac{(T^{5/6})^2\cdot T^{-1/4}\lambda_j}{12 T}\right)\nonumber\\
&\leq 2 \exp \left( -\frac{\lambda_{min}}{12}T^{5/12}\right),\label{eq: poisson bound 1}
\end{align}
where $\lambda_{min}=\min_{j\in[n]}\lambda_j$. 
In case 3), $\tilde{z}_j(t^*)$ is Poisson random variable with parameter $t^*{\lambda}_j$. We have
\begin{align*}
\Prob\left(\tilde{z}_j(t^*)-t^*x^*_j>\frac{(T-t^*)x^*_j}{2}\right)=&\Prob\left(\tilde{z}_j(t^*)-t^*{\lambda}_j>\frac{(T-t^*)x^*_j}{2}-t^*(\lambda_j-x^*_j)\right)\\
\leq & \Prob\left(\tilde{z}_j(t^*)-t^*{\lambda}_j>\frac{T^{5/6}{\lambda}_j/2}{2}-{\lambda}_jT\cdot T^{-1/4}\right).
\end{align*}
Choose a constant $T_0$ such that $T^{5/6}/4\geq 2 T^{3/4}$ for $T\geq T_0$.
It follows from the Poisson tail bound (Lemma \ref{lem: poisson bound}) that
\begin{align}
\Prob\left(\tilde{z}_j(t^*)-t^*{\lambda}_j>\frac{T^{5/6}{\lambda}_j/2}{2}-{\lambda}_jT\cdot T^{-1/4}\right)
&\leq 2 \exp \left( -\frac{1}{3}\left(T^{3/4}\lambda_j\right)^2\frac{1}{t^*\lambda_j }\right)\nonumber\\
&\leq 2 \exp \left(-\frac{1}{3}\frac{\lambda_j T^{3/2}}{T}\right)\nonumber\\
&\leq 2 \exp \left( -\frac{\lambda_{min}}{3}T^{1/2}\right).\label{eq: poisson bound 1-3}
\end{align}

Let $J_\lambda=\{j: x^*_j=\lambda_j\}$.
The second term can be written as
\begin{align}
\Prob\left(\Gamma(T)>\frac{(T-t^*)x^*_j}{2}\right)
&=\Prob\left(\alpha \sum_{j\in J_\lambda}|\Lambda_j(T)-\lambda_jT|>\frac{(T-t^*)x^*_j}{2}\right) \nonumber\\
&\leq \sum_{j\in J_\lambda}\Prob\left(\alpha|\Lambda_j(T)-\lambda_jT|>\frac{(T-t^*)x^*_j}{2|J_\lambda|}\cdot I(|J_\lambda|>0)\right)\nonumber\\
&= \sum_{j\in J_\lambda}\Prob\left(|\Lambda_j(T)-\lambda_jT|>\frac{(T-t^*)x^*_j}{2\alpha|J_\lambda|+1}\right).\nonumber
\end{align}
Applying the Poisson tail bound (Lemma \ref{lem: poisson bound}) again, we have
\begin{align}
\Prob\left(|\Lambda_j(T)-\lambda_jT|>\frac{(T-t^*)x^*_j}{2\alpha|J_\lambda|+1}\right)
&\leq 2\exp\left(-\frac{1}{3}\left(\frac{(T-t^*)x^*_j}{2\alpha|J_\lambda|+1}\right)^2\frac{1}{\lambda_j T}\right)\nonumber\\
&\leq 2 \exp\left( -\frac{ (T^{5/6}\cdot \lambda_jT^{-1/4})^2}{3(2\alpha|J_\lambda|+1)^2\lambda_j T}\right)\nonumber\\
&\leq 2 \exp\left( -\frac{\lambda_{min}}{3(2\alpha|J_\lambda|+1)^2}T^{1/6}\right)\nonumber.
\end{align}
Therefore, we have
\begin{align}
\Prob\left(\Gamma(T)>\frac{(T-t^*)x^*_j}{\color{black}2}\right) \leq 2|J_{\lambda}| \exp\left(  -\frac{\lambda_{min}}{3(2\alpha|J_\lambda|+1)^2}T^{1/6}\right).\label{eq: poisson bound 2}
\end{align}
%
%
From \eqref{eq: union bound E1'}, \eqref{eq: poisson bound 1},\eqref{eq: poisson bound 1-3} and \eqref{eq: poisson bound 2}, it follows that
\begin{align}
\Prob(E_{1,j}^\complement)	\leq O\left(\exp\bigl(-\frac{\lambda_{min}}{3(2\alpha|J_\lambda|+1)^2}T^{1/6}\bigr)\right).\label{eq:p_E1'}
\end{align}

We can also apply the similar argument to bound $\Prob(E_{2,j}^\complement)$. 
That is, if the event $E_{2,j}^\complement$ happens, at least one of the following three events must happen:
\begin{align}
\left\{t^*x^*_j-\tilde{z}_j(t^*)>\tfrac{(T-t^*)(\lambda_j-x^*_j)}{3}\right\} ,\;\left\{\Gamma(T)>\tfrac{(T-t^*)(\lambda_j-x^*_j)}{3}\right\} ,\;\left\{|\Delta_j(t^*)|>\tfrac{(T-t^*)(\lambda_j-x^*_j)}{3} \right\}. \nonumber
\end{align}
We can apply the Poisson tail bound (Lemma \ref{lem: poisson bound}) to these three events.
{\color{black} We consider only case 1) and case 2) here, because case 3) is eliminated in definition of event E. For the first event, in case 1) we have $\tilde{z}_j(t^*)=0$. It follows that
	\begin{align}
	\Prob\left(t^*x^*_j-\tilde{z}_j(t^*)>\frac{(T-t^*)(\lambda_j-x^*_j)}{3}\right)&=I\left(t^*x^*_j>\frac{(T-t^*)(\lambda_j-x^*_j)}{3}\right)\nonumber\\
	&=I\left(\frac{(T+2t^*)x^*_j}{3}>\frac{(T-t^*)\lambda_j}{3}\right)\nonumber\\
	&\leq I\left(\frac{3T\lambda_jT^{-1/4}}{3}>\frac{T^{5/6}\lambda_j}{3}\right)\nonumber\\
	&=I\left(T^{3/4}\lambda_j>\frac{T^{5/6}\lambda_j}{3}\right).\nonumber
	\end{align}}
	Choosing a constant $T_0$ such that $T^{5/6}/3\geq T^{3/4}$ for $T\geq T_0$, we have
	\begin{align}
	\Prob\left(t^*x^*_j-\tilde{z}_j(t^*)>\frac{(T-t^*)(\lambda_j-x^*_j)}{3}\right)&\leq 0. \label{eq:poisson bound2 1.1}
	\end{align}	
	In case 2), $\tilde{z}_j(t^*)$ is Poisson random variable with {\color{black}parameter} $t^*x^*_j$, thus we have
	\begin{align}
	\Prob\left(t^*x^*_j-\tilde{z}_j(t^*)>\frac{(T-t^*)(\lambda_j-x^*_j)}{3}\right)
	&\leq 2 \exp \left( -\frac{1}{3}\left(\frac{(T-t^*)(\lambda_j-x^*_j)}{3}\right)^2\frac{1}{t^*x^*_j}\right)\nonumber\\
	&\leq 2 \exp \left(-\frac{(T-t^*)^2 (\lambda_j-x^*_j)^2}{27 T\lambda_j}\right)\nonumber\\
	&\leq 2 \exp \left(-\frac{(T^{5/6})^2 (\lambda_jT^{-1/4})^2}{27 T\lambda_j}\right)\nonumber\\
	&\leq 2 \exp \left( -\frac{\lambda_{min}}{27}T^{1/6}\right).\label{eq:poisson bound2 1.2}
	\end{align}
	For the second term, we have
	\begin{align}
	\Prob\left(\Gamma(T)>\frac{(T-t^*)(\lambda_j-x^*_j)}{3}\right)
	&=\Prob\left(\alpha \sum_{j\in J_\lambda}|\Lambda_j(T)-\lambda_jT|>\frac{(T-t^*)(\lambda_j-x^*_j)}{3}\right) \nonumber\\
	&\leq \sum_{j\in J_\lambda}\Prob\left(|\Lambda_j(T)-\lambda_jT|>\frac{(T-t^*)(\lambda_j-x^*_j)}{3\alpha|J_\lambda|+1}\right)\nonumber\\
	&\leq \sum_{j\in J_\lambda} 2\exp\left(-\frac{1}{3}\left(\frac{(T-t^*)(\lambda_j-x^*_j)}{3\alpha|J_\lambda|+1}\right)^2\frac{1}{\lambda_j T}\right)\nonumber\\
	&\leq \sum_{j\in J_\lambda} 2 \exp\left( -\frac{ (T^{5/6}\cdot \lambda_jT^{-1/4})^2}{3(3\alpha|J_\lambda|+1)^2\lambda_j T}\right)\nonumber\\
	&\leq 2|J_{\lambda}|\exp\left( -\frac{\lambda_{min}}{3(3\alpha|J_\lambda|+1)^2}T^{1/6}
	\right).\label{eq:poisson bound2 2}
	\end{align}
	Moreover, for the last term we have
	\begin{align}
	\Prob\left(|\Delta_j(t^*)|>\frac{(T-t^*)(\lambda_j-x^*_j)}{3}\right)
	&=\Prob\left( |\Lambda_j(T)-\Lambda_j(t^*)-\lambda_j(T-t^*)|>\frac{(T-t^*)(\lambda_j-x^*_j)}{3} \right)\nonumber\\
	&\leq 2 \exp \left( -\frac{1}{3}\left(\frac{(T-t^*)(\lambda_j-x^*_j)}{3}\right)^2\frac{1}{\lambda_j(T-t^*)}\right)\nonumber\\
	&\leq 2 \exp \left( -\frac{ T^{5/6}(\lambda_j T^{-1/4})^2}{27\lambda_j }\right)\nonumber\\
	&\leq 2 \exp \left( -\frac{\lambda_{min}}{27}T^{1/3}\right).\label{eq:poisson bound3}
	\end{align}
	From union bound, it follows from \eqref{eq:poisson bound2 1.1}, \eqref{eq:poisson bound2 1.2}, \eqref{eq:poisson bound2 2} and \eqref{eq:poisson bound3} that
	\begin{align}
	\Prob(E_{2,j}^\complement)\leq O\left(\exp\left(-\frac{\lambda_{min}}{27(\alpha|J_\lambda|+1)^2}T^{1/6}\right)\right).\label{eq:p_E2'}
	\end{align}
	Let
	\[
	\kappa = \frac{\lambda_{min}}{27(\alpha|J_\lambda|+1)^2},
	\]
	The results from \eqref{eq:p_E'}, \eqref{eq:p_E1'} and \eqref{eq:p_E2'} lead to
	\begin{align}
	\Prob(E^\complement)
	&=O\left(n |J_\lambda| \exp \bigl(
	-\kappa T^{1/6}
	\bigr)\right)\nonumber,
	\end{align}
	where the big $O$ notation hides an absolute constant.
	\Halmos
	\endproof


	\begin{lemma}\label{lemma: prob events in HO-FR}
	Probabilities of the following events are bounded by
	\begin{align*}
	\Prob(Q_1)
	\geq {\color{black}0.0013}-\frac{0.9496}{\sqrt{T'}}, \quad \Prob(Q_2) 
	&\geq {\color{black}0.5}, \quad \Prob(Q_3) 
	\geq{\color{black}9.8531 \times 10^{-10}}-\frac{0.9496}{\sqrt{T'}}\\
	\Prob(B^\complement \cap B_1^\complement\cap  Q_1 \cap Q_2)
	&\leq {\color{black}e^{-0.0026\sqrt{T'}}},
	\end{align*}
	where the events $Q_1, Q_2, Q_3, B$ and $B_1$ are defined in \eqref{event:A1}--\eqref{event:B1} respectively.
\end{lemma}

\proof{Proof of Lemma \ref{lemma: prob events in HO-FR}.}
{\color{black} We will apply Berry-Esseen theorem (see formal statement in Lemma \ref{lem:berry-esseen}) to bound $\Prob(Q_1)$ and $\Prob(Q_3)$.}
The probability $\mathbb{P}(Q_1)$ can be written as
\begin{align}
\Prob(Q_1)&=\Prob(T'-{\color{black}4}\sqrt{T'}\leq \Lambda_1(0,T')\leq T'-{\color{black}3}\sqrt{T'})\nonumber\\
&=\Prob\left(-{\color{black}4}\leq \frac{\Lambda_1(0,T')-T'}{\sqrt{T'}}\leq -{\color{black}3}\right)\nonumber\\
&=F_{T'}(-{\color{black}3})-F_{T'}(-{\color{black}4})\nonumber\\
&=\Phi(-{\color{black}3})-\Phi(-{\color{black}4})-(\Phi(-{\color{black}3})-F_{T'}(-{\color{black}3}))-(F_{T'}(-{\color{black}4})-\Phi(-{\color{black}4}))\nonumber\\
&\geq \Phi(-{\color{black}3})-\Phi(-{\color{black}4})- |F_{T'}(-{\color{black}3})-\Phi(-{\color{black}3})|-|F_{T'}(-{\color{black}4})-\Phi(-{\color{black}4})|,\label{eq:berry4}
\end{align}
{\color{black}
	where $F_{T'}$ is a CDF of $\frac{\Lambda_1(0,T')-T'}{\sqrt{T'}}$. Recall that by the stationary and independent increment properties of Poisson processes, $\Lambda_1(0,T')-T'$ can be thought as the summation of $T'$ i.i.d. $\Lambda_1(1)-1$ random variables with $\E[\Lambda_1(1)-1]=0$, $\E[(\Lambda_1(1)-1)^2]=1$ and $\E[|\Lambda_1(1)-1|^3]=1$. Hence, the second and the third term of \eqref{eq:berry4} can be bounded by the Berry-Esseen theorem (Lemma \ref{lem:berry-esseen}).
}
 That is, for any $x$, we have
\begin{align}
|F_{T'}(x)-\Phi(x)|\leq \frac{0.4748}{\sqrt{T'}}.\label{eq:berrytheorem}
\end{align}
Combining \eqref{eq:berrytheorem} to \eqref{eq:berry4}, we can conclude that
\begin{align}
\Prob(Q_1)
&\geq {\color{black}0.0013}-\frac{0.9496}{\sqrt{T'}}.\label{eq:bound1}
\end{align}
{\color{black}
We will apply the same argument to bound the probablity $\Prob(Q_3)$. That is, we have
}
\begin{align}
\Prob(Q_3)
&=\Prob(T'+{\color{black}6}\sqrt{T'}\leq \Lambda_1(T'',T)\leq T'+{\color{black}7}\sqrt{T'})\nonumber\\
&=\mathbb{P}({\color{black}6}\leq \frac{\Lambda_1(T'',T)-T'}{\sqrt{T'}}\leq {\color{black}7})\nonumber\\
&=F_{T'}({\color{black}7})-F_{T'}({\color{black}6})\nonumber\\
&\geq\Phi({\color{black}7})-\Phi({\color{black}6})-2\frac{0.4748}{\sqrt{T'}}\label{eq:berry_A3}\\
&={\color{black}9.8531 \times 10^{-10}}-\frac{0.9496}{\sqrt{T'}},\label{eq:bound3}
\end{align}
where \eqref{eq:berry_A3} follows from the {\color{black}Lemma \ref{lem:berry-esseen}. Next, to bound the probability $\Prob(Q_2)$, we will apply Doob's maximal inequality (Lemma \ref{lem:doob}). Recall that the event $Q_2$ is defined as
\[
Q_2=\{(t-T')-{\color{black}2}\sqrt{T'}\leq \Lambda_1(T',t)\leq (t-T')+{\color{black}2}\sqrt{T'},\; \forall t \in (T',T'']\}.
\]
Equivalently, this event can be re-written as
\[
Q_2=\{\sup_{t\in(T',T'']}|\Lambda_1(T',t)-(t-T')|\leq 2\sqrt{T'}\}.
\]
Hence, we have
\begin{align}
	\Prob(Q_2)&=\Prob(\sup_{t\in(T',T'']}|\Lambda_1(T',t)-(t-T')|\leq 2\sqrt{T'})\nonumber\\
	&=1-\Prob(\sup_{t\in(T',T'']}|\Lambda_1(T',t)-(t-T')|> 2\sqrt{T'}).\label{eq:Q2'}
\end{align}
To bound the last term of \eqref{eq:Q2'}, we first observe that $\Lambda_1(T',t)-(t-T')$ is a martingale with paths that are right continuous with left limits, so we can apply  Doob's maximal inequality (Lemma \ref{lem:doob}) to \eqref{eq:doob}, and then use Cauchy-Schwarz inequality in \eqref{eq:cauchy after doob}. We have
\begin{align}
	\Prob(\sup_{t\in(T',T'']}|\Lambda_1(T',t)-(t-T')|> 2\sqrt{T'})
	&\leq \frac{\E[|\Lambda_1(T',T'')-(T''-T')|]}{2\sqrt{T'}} \label{eq:doob}\\
&\leq 	\frac{\sqrt{\E[(\Lambda_1(T',T'')-(T''-T'))^2]}}{2\sqrt{T'}} \label{eq:cauchy after doob}\\
&=\frac{\sqrt{\Var(\Lambda_1(T',T''))}}{2\sqrt{T'}}=\frac{\sqrt{T'}}{2\sqrt{T'}}=0.5.\label{eq:bound Q2'}
\end{align}
Combining the result in \eqref{eq:bound Q2'} to \eqref{eq:Q2'}, we get
\begin{align}
	\Prob(Q_2)\geq 1-0.5=0.5.\label{eq:bound2}
\end{align}
Next, we will bound the probability $\Prob(B^\complement \cap B_1^\complement \cap  Q_1 \cap Q_2)$.
}
If the event $B^\complement \cap B_1^\complement \cap  Q_1 \cap Q_2$ happens, {\color{black} we can observe that $z_1(0,T')\leq \Lambda_1(0,T')\leq T'-3\sqrt{T'}$ from $Q_1$ and $z_2(0,T')\leq z_2(0,T)<\frac{1}{6}\sqrt{T'}$ from $B^\complement$, and hence} the remaining capacity at time $T'$ {\color{black}will} be
\begin{align*}
C(T')&=C-z_1(0,T')-z_2(0,T')\geq T-(T'-{\color{black}3}\sqrt{T'})-\frac{1}{6}\sqrt{T'}=T-T'+\frac{\color{black}17}{6}\sqrt{T'}.
\end{align*}
{\color{black}Similarly, for any time $t\in(T',T'']$, the event $B^\complement \cap B_1^\complement \cap  Q_1 \cap Q_2$ implies that  $z_1(0,t)\leq \Lambda_1(0,T')+\Lambda_1(T',t)\leq T'-3\sqrt{T'}+(t-T)+2\sqrt{T'}$ from $Q_1$ and $Q_2$ and $z_2(0,t)\leq z_2(0,T)<\frac{1}{6}\sqrt{T'}$ from $B^\complement$; therefore,}
the remaining capacity at any time $t\in(T',T'']$ will be
\begin{align*}
C(t)&=C-z_1(0,t)-z_2(0,t) \geq T-(T'-{\color{black}3}\sqrt{T'})-((t-T')+{\color{black}2}\sqrt{T'})-\frac{1}{6}\sqrt{T'}=T-t+\frac{5}{6}\sqrt{T'}.
\end{align*}
So the average capacity per period at time $T'$ is
\begin{align}
b(T')&\geq \frac{T-T'+\frac{11}{6}\sqrt{T'}}{T-T'} =1+\frac{11\sqrt{T'}/6}{2T'} \geq 1+\frac{11}{12\sqrt{T'}},\label{eq:lb_b at T'}
\end{align}
and similarly the average capacity per period at any time $t\in(T',T'']$ is given by
\begin{align}
b(t)&\geq \frac{T-t+\frac{5}{6}\sqrt{T'}}{T-t}=1+\frac{5\sqrt{T'}/6}{T-t} \geq 1+\frac{5}{12\sqrt{T'}}.\label{eq:lb_b}
\end{align}
{\color{black}Recall that for the problem instance we consider, the admission probability of class 1 customer, which is obtained from the LP described in Algorithm \ref{algo:res}, is given by $x_1(t)=\min(b(t)/\lambda_1, 1)=\min(b(t),1)$. Thus, \eqref{eq:lb_b at T'} and \eqref{eq:lb_b} implies that}
the decision maker must accept all arrivals of class 1 customer in phase 2. {\color{black}Hence, it follows from definition of the event $Q_2$ in \eqref{event:A2}} that, for any time $t\in(T',T'']$, we have
\begin{align}
(t-T')-{\color{black}2}\sqrt{T'}\leq z_1(t,T')\leq (t-T')+{\color{black}2}\sqrt{T'}.\label{eq: arrival in 2nd phase}
\end{align}
{\color{black}Next, $B_1^\complement$ implies that} we can upper bound the remaining capacity at time $T'$ as
\[
C(T'){\color{black}= C-z_1(0,T')-z_2(0,T')} \leq T-z_1(0,T')
\leq T-(T'-{\color{black}10}\sqrt{T'})
=T-T'+{\color{black}10}\sqrt{T'},
\]
which results in
\begin{align}
b(T')&\leq \frac{T-T'+{\color{black}10}\sqrt{T'}}{T-T'} = 1 +\frac{{\color{black}10}\sqrt{T'}}{T-T'} \leq 1+\frac{\color{black}5}{\sqrt{T'}}.\label{eq:ub_b at T'}
\end{align}
Similarly, we can upper bound the remaining capacity at time $t\in (T',T'']$ using the results in \eqref{eq: arrival in 2nd phase} which yields
\begin{align}
C(t) &{\color{black}= C-z_1(0,t)-z_2(0,t)\leq T-z_1(0,T')-z_1(T',t)}\nonumber\\
&\leq T-(T'-{\color{black}10}\sqrt{T'})-((t-T')-{\color{black}2}\sqrt{T'}) =T-t+{\color{black}12}\sqrt{T'}.\nonumber
\end{align}
It follows that the upper bound of the average capacity per period at time $t\in (T',T'']$ is given by
\begin{align}
b(t) \leq \frac{T-t+{\color{black}12}\sqrt{T'}}{T-t} = 1 +\frac{{\color{black}12}\sqrt{T'}}{T-t} \leq 1+\frac{\color{black}12}{\sqrt{T'}}.\label{eq:ub_b}
\end{align}
Combining the results from \eqref{eq:lb_b at T'}, \eqref{eq:lb_b}, \eqref{eq:ub_b at T'} and \eqref{eq:ub_b}, we obtain the bound of the average capacity per period at $t\in[T',T'']$, that is,
\begin{align}
1+\frac{5}{12\sqrt{T'}}\leq b(t)\leq 1+\frac{\color{black}12}{\sqrt{T'}},
\end{align}
which also implies that the solution to the LP at time $t\in[T',T'']$ satisfies
$\frac{5}{12\sqrt{T'}}\leq x_2(t)\leq \frac{\color{black}12}{\sqrt{T'}}.$
Therefore, the probability of the event $B^\complement \cap B_1^\complement\cap  Q_1 \cap Q_2$ can be written as
\begin{align}
\mathbb{P}(B^\complement \cap B_1^\complement\cap  Q_1 \cap Q_2) &\leq \mathbb{P}\left(\frac{5}{12\sqrt{T'}}\leq x_2(t)\leq \frac{\color{black}12}{\sqrt{T'}}{\color{black},\forall t \in [T',T'']},\; z_2(T',T'')<\frac{1}{6}\sqrt{T'}\right).\label{eq:x2 bound}
\end{align}
We will use Freedman's Inequality {\color{black}(Lemma~\ref{lem:freedman})} to bound \eqref{eq:x2 bound}. {\color{black}For $i\in[T'],$ we let $\xi_0=0$ and $\xi_i=x_2(T'+i-1)-z_2(T'+i-1,T'+i).$ We observe that $\xi_i$ is $\mathcal{F}_{T'+i-1}$-measurable and $\E[\xi_i|\mathcal{F}_{T'+i-1}]=x_2(T'+i-1)-x_2(T'+i-1)=0$. Thus, $\xi_i$ is a martingale difference. Set
	\[
	S_{T'}=\sum_{i=1}^{T'}\xi_i=\sum_{i=1}^{T'}x_2(T'+i-1)-z_2(T',T'').
	\]
Then, conditions in \eqref{eq:x2 bound} imply that
\[
S_{T'}>\frac{5}{12\sqrt{T'}}T'-\frac{1}{6}\sqrt{T'}=\frac{1}{4}\sqrt{T'}.
\]
Moreover,  let $\langle S\rangle_k=\sum_{i=1}^k\E[\xi_i^2|\mathcal{F}_{i-1}]$ for $k\geq 1$. We have
\begin{align}
	\langle S\rangle_{T'}&=\sum_{i=1}^{T'}\E[\xi_i^2|\mathcal{F}_{T'+i-1}]=\sum_{i=1}^{T'}\E[(x_2(T'+i-1)-z_2(T'+i-1,T'+i))^2|\mathcal{F}_{T'+i-1}]\nonumber\\
	&=\sum_{i=1}^{T'}\Var(z_2(T'+i-1,T'+i)|\mathcal{F}_{T'+i-1})=T'x_2(T'+i-1).\nonumber
\end{align}
Then, conditions in \eqref{eq:x2 bound} imply that
\[
	\langle S\rangle_{T'}\leq T'\frac{12}{\sqrt{T'}}=12\sqrt{T'}.
\] 
Therefore, we have
\[
\mathbb{P}\left(\frac{5}{12\sqrt{T'}}\leq x_2(t)\leq \frac{\color{black}12}{\sqrt{T'}}{\color{black},\forall t \in [T',T'']},\; z_2(T',T'')<\frac{1}{6}\sqrt{T'}\right)\leq \Prob\left(S_{T'}>\frac{1}{4}\sqrt{T'}, \langle S \rangle_{T'}\leq 12\sqrt{T'}\right).
\]
Since, for all $i\in[T']$, we have $\xi_i\leq x_2(T+i-1)\leq \lambda_2=1$, we can apply Freedman's inequality (Lemma \ref{lem:freedman}) and get
\begin{align*}
\Prob\left(S_{T'}>\frac{1}{4}\sqrt{T'}, \langle S \rangle_{T'}\leq 12\sqrt{T'}\right)
&\leq\exp\left(-\frac{(\sqrt{T'}/4)^2}{2(12\sqrt{T'}+\sqrt{T'}/4)}\right)
=\exp\left(-\frac{1}{392}\sqrt{T'}\right)
=e^{-0.0026\sqrt{T'}}.\Halmos
\end{align*}
\endproof
}

\begin{lemma}[{\color{black}Bound on $b_l$}]\label{lem: expectation b} We have
	$
	\E[(b_l-b_l(t))^+]\leq K_l\sqrt{\sum_{i=0}^{t-1}\frac{1}{(T-i-1)^2}},
	$
	where $K_l=\sqrt{\sum_{j=1}^na_{lj}^2 \lambda_j^2}$.
\end{lemma}

\proof{Proof of Lemma \ref{lem: expectation b}.}
{\color{black}Recall from \eqref{eq:def of resolving} that we have
	\[
	(b_l-b_l(t))^+\leq\left[\sum_{i=0}^{t-1}\left(\frac{\sum_{j=1}^n a_{lj}(\tilde{z}_j(i+1)-\tilde{z}_j(i))}{T-i-1}-\frac{\sum_{j=1}^na_{lj}x_j(i)}{T-i-1}\right)\right]^+.
	\]}
Taking expectations on both sides yields
\begin{align}
\E[(b_l-b_l(t))^+]
&\leq \E\left[\left(\sum_{i=0}^{t-1}\frac{\sum_{j=1}^n a_{lj}(\tilde{z}_j(i+1)-\tilde{z}_j(i))-\sum_{j=1}^na_{lj}x_j(i)}{T-i-1}\right)^+\right]\nonumber\\
&\leq \E\left[\left|\sum_{i=0}^{t-1}\frac{\sum_{j=1}^n a_{lj}(\tilde{z}_j(i+1)-\tilde{z}_j(i))-\sum_{j=1}^na_{lj}x_j(i)}{T-i-1}\right|\right]\label{eq:cauchy} \\
&\leq \sqrt{\E\left[\left(\sum_{i=0}^{t-1}\frac{\sum_{j=1}^n a_{lj}(\tilde{z}_j(i+1)-\tilde{z}_j(i))-\sum_{j=1}^na_{lj}x_j(i)}{T-i-1}\right)^2\right]}.\label{eq:cond_ex}
\end{align}
The last line applies the Cauchy-Schwarz inequality.

Let $\{\mathcal{F}_t, 0\leq t\leq T\}$ be the filtration generated by $\{\Lambda_j(s),0\leq s\leq t, j\in [n]\}$. By the law of total expectation, \eqref{eq:cond_ex} becomes
\begin{align}
\E[(b_l-b_l(t))^+] \leq \sqrt{\E\left[\E\left[\left(\sum_{i=0}^{t-1}\frac{\sum_{j=1}^n a_{lj}(\tilde{z}_j(i+1)-\tilde{z}_j(i))-\sum_{j=1}^na_{lj}x_j(i)}{T-i-1}\right)^2\middle| \mathcal{F}_i\right]\right]}.\label{exp sq sum}
\end{align}
{\color{black} Since the conditional independence of the arrivals of the customers in different period implies that the arrivals of the admitted customers in different period are also conditionally independent, and each of the summands has mean zero, and the cross-terms vanish. Thus, \eqref{exp sq sum} equals to}
\begin{align}
\E[(b_l-b_l(t))^+] \leq \sqrt{\E\left[\E\left[\sum_{i=0}^{t-1}\left(\frac{\sum_{j=1}^n a_{lj}(\tilde{z}_j(i+1)-\tilde{z}_j(i))-\sum_{j=1}^na_{lj}x_j(i)}{T-i-1}\right)^2\middle| \mathcal{F}_i\right]\right]}\nonumber.
\end{align}
From the description of the $\mathsf{FR}$, we know that $\tilde{z}_j(t+1)-\tilde{z}_j(t)$ conditioned on $\mathcal{F}_t$ is distributed as Poisson distribution with parameter $x_j(t)$ (from the Poisson thinning property). {\color{black} We use the definition of variance in Equation \eqref{eq:variance}, and we use the observation that $x(t)$ is bounded above by $\lambda_j$ because $x(t)$ is the solution to the {\color{black}LP} in \eqref{eq: x(t) bounded by lambda}.} It immediately follows that $\E[\tilde{z}_j(t+1)-\tilde{z}_j(t)|\mathcal{F}_t]=x_j(t)$ and $\Var(\tilde{z}_j(t+1)-\tilde{z}_j(t)|\mathcal{F}_t)=x_j(t)$.
Therefore, we can write
\begin{align}
\E[(b_l-b_l(t))^+]
&\leq \sqrt{\E\left[\sum_{i=0}^{t-1}\frac{1}{(T-i-1)^2}\E\left[\left(\sum_{j=1}^n a_{lj}(\tilde{z}_j(i+1)-\tilde{z}_j(i))-\sum_{j=1}^na_{lj}x_j(i)\right)^2\middle| \mathcal{F}_i\right]\right]}\nonumber\\
&= \sqrt{\E\left[\sum_{i=0}^{t-1}\frac{1}{(T-i-1)^2} \Var\left(\sum_{j=1}^n a_{lj}(\tilde{z}_j(i+1)-\tilde{z}_j(i))\middle| \mathcal{F}_i\right) \right]}\label{eq:variance}\\
&= \sqrt{\sum_{i=0}^{t-1}\frac{1}{(T-i-1)^2}\sum_{j=1}^na_{lj}^2 \E\left[ \Var\left(\tilde{z}_j(i+1)-\tilde{z}_j(i)| \mathcal{F}_i\right) \right]}\nonumber\\
&=\sqrt{\sum_{i=0}^{t-1}\frac{1}{(T-i-1)^2}\sum_{j=1}^na_{lj}^2 \E\left[ x_j(t)^2\right]}\nonumber\\
&\leq \sqrt{\sum_{i=0}^{t-1}\frac{1}{(T-i-1)^2}\sum_{j=1}^na_{lj}^2 \lambda_j^2}\label{eq: x(t) bounded by lambda}\\
&=K_l\sqrt{\sum_{i=0}^{t-1}\frac{1}{(T-i-1)^2}}, \nonumber
\end{align}
where $K_l=\sqrt{\sum_{j=1}^na_{lj}^2 \lambda_j^2}$.
\Halmos
\endproof

}

%
\end{APPENDIX}

%

\end{document}